\documentclass[12pt,a4paper]{amsart}
\usepackage{amsmath,amsfonts,amssymb,amscd,amsthm}

\newtheorem{theorem}[equation]{Theorem}
\newtheorem{lemma}[equation]{Lemma}
\newtheorem{proposition}[equation]{Proposition}
\newtheorem{definition}[equation]{Definition}
\newtheorem{corollary}[equation]{Corollary}
\newtheorem{conjecture}[equation]{Conjecture}
\theoremstyle{remark}
\newtheorem{remark}[equation]{Remark}

\numberwithin{equation}{section}

\DeclareMathOperator{\proj}{proj}

\DeclareMathOperator{\End}{End}
\renewcommand\Re{{\rm Re}}
\renewcommand\Im{{\rm Im}}
\DeclareMathOperator{\Hom}{Hom}
\DeclareMathOperator{\Ext}{Ext}
\DeclareMathOperator{\Gal}{Gal}
\DeclareMathOperator{\Id}{Id}
\DeclareMathOperator{\Ind}{Ind}
\DeclareMathOperator{\Irr}{Irr}
\DeclareMathOperator{\Res}{Res}

\DeclareMathOperator{\Sh}{Sh}

\DeclareMathOperator{\St}{St}

\DeclareMathOperator{\Trace}{Trace}
\newcommand{\bbC}{{\mathbb C}}
\newcommand{\bbF}{{\mathbb F}}

\newcommand{\bbN}{{\mathbb N}}

\newcommand{\bbZ}{{\mathbb Z}}

\newcommand{\bB}{{\bf B}}

\newcommand{\bG}{{\bf G}}

\newcommand{\bI}{{\bf I}}

\newcommand{\bL}{{\bf L}}

\newcommand{\bP}{{\bf P}}

\newcommand{\bS}{{\bf S}}
\newcommand{\bT}{{\bf T}}
\newcommand{\bU}{{\bf U}}
\newcommand{\bV}{{\bf V}}
\newcommand{\bW}{{\bf W}}
\newcommand{\bX}{{\bf X}}
\newcommand{\bY}{{\bf Y}}
\newcommand{\bZ}{{\bf Z}}
\newcommand{\ba}{{\bf a}}
\newcommand{\bb}{{\bf b}}
\newcommand{\bc}{{\bf c}}

\newcommand{\bs}{{\bf s}}
\newcommand{\bt}{{\bf t}}
\newcommand{\bu}{{\bf u}}
\newcommand{\bv}{{\bf v}}
\newcommand{\bw}{{\bf w}}
\newcommand{\bx}{{\bf x}}
\newcommand{\by}{{\bf y}}
\newcommand{\bz}{{\bf z}}
\def\bpi{{\boldsymbol\pi}}
\def\bphi{{\boldsymbol\phi}}
\def\bpsi{{\boldsymbol\psi}}
\def\bsigma{{\boldsymbol\sigma}}
\newcommand{\CA}{{\mathcal A}}
\newcommand{\CB}{{\mathcal B}}
\newcommand{\CC}{{\mathcal C}}
\newcommand{\cD}{{\mathcal D}} 

\newcommand{\CF}{{\mathcal F}}
\newcommand{\CH}{{\mathcal H}}

\newcommand{\CO}{{\mathcal O}}

\newcommand{\CR}{{\mathcal R}}

\newcommand{\dz}{{\dot z}}
\newcommand{\Sgot}{{\mathfrak S}}
\newcommand{\BF}{{\bB^F}}
\newcommand{\GF}{{\bG^F}}
\newcommand{\LF}{{\bL^F}}
\newcommand{\PF}{{\bP^F}}

\newcommand{\UF}{{\bU^F}}

\newcommand\Fqbar{{\overline\bbF_q}}
\newcommand{\Qlbar}{{\overline{\mathbb Q}_\ell}}
\newcommand{\semi}{{\rtimes}}
\newcommand{\GL}{{\text{GL}}}
\newcommand{\Bred}{{\bW}}
\newcommand{\BW}{{B^+}}

\newcommand{\inv}{^{-1}}

\newcommand{\cf}{{\it cf.}}
\newcommand{\eg}{{\it e.g.}}
\newcommand{\ie}{{\it i.e.}}
\newcommand{\loccit}{{\it loc.\ cit.}}

\newcommand{\TrH}[2]{{\sum_i (-1)^i\Trace(#1|H^i_c(#2,\Qlbar))}}
\newcommand{\lexp}[2]{\kern\scriptspace\vphantom{#2}^{#1}\kern-\scriptspace#2}
\newcommand{\scal}[3]{{\langle\,#1,#2\,\rangle_{#3}}}
\newcommand{\genby}[1]{\mathopen<#1\mathclose>}
\newcommand{\sdp}[1]{\semi\genby{#1}}
\begin{document}
\author{F.~Digne and J.~Michel}
\title[Endomorphisms of Deligne-Lusztig varieties]
{Endomorphisms of Deligne-Lusztig varieties}
\begin{abstract}
We study some conjectures on the endomorphism algebras of the cohomology
of  Deligne-Lusztig  varieties  which  are  a  refinement  of  those  of
\cite{Sydney}.
\end{abstract}
\maketitle
\tableofcontents
\vfill\eject
\section{Introduction}\label{intro}
Let  $\bG$ be  a connected  reductive algebraic  group, defined  over an
algebraic closure  $\bbF$ of a  finite field of characteristic  $p$. Let
$F$  be an  isogeny on  $\bG$  such that  some power  $F^\delta$ is  the
Frobenius endomorphism  attached to a  split $\bbF_{q^\delta}$-structure
on $\bG$ (where $q$ is a real  number such that $q^\delta$ is a power of
$p$).  The finite  group $\GF$  of fixed  points under  $F$ is  called a
finite group of Lie type. When considering a simple group which is not a
Ree  or  Suzuki  groups me  may  take  $F$  to  be already  a  Frobenius
endomorphism.  

Let $W$  be the  Weyl group  of $\bG$  and let  $B$ (resp
$B^+$) be the corresponding braid group (resp. monoid).
The canonical morphism of monoids $\beta:B^+\to W$ has a section that we
denote by  $w\to\bw$: it sends  an element of  $W$ to the  only positive
braid  $\bw$  such that  $\beta(\bw)=w$  and  such  that the  length  of
$\bw$  in $B^+$  is the  same as  the Coxeter  length of  $w$; we  write
$\bW=\{\bw\mid w\in W\}$ and $\bS=\{\bs\mid s\in S\}$ where $S$ is a set
of Coxeter generators for $W$.

Let  us recall  how in  \cite{Sydney} a  ``Deligne-Lusztig variety''  is
attached to  each element  of $B^+$.
Let  $\CB$ be  the  variety  of Borel  subgroups  of  $\bG$. The  orbits
of  $\CB\times\CB$  are in  natural  bijection  with $W$.  Let  $\CO(w)$
be  the orbit  corresponding  to  $w\in W$.  Let  $\bb\in  B^+$ and  let
$\bb=\bw_1\ldots\bw_n$  be a  decomposition  of $\bb$  as  a product  of
elements  of  $\bW$. To  such  a  decomposition  we attach  the  variety
$\{(\bB_1,\ldots,\bB_{n+1})\mid (\bB_i,\bB_{i+1}) \in \CO(w_i)\text{ and
} \bB_{n+1}=F(\bB_1)\}$. It is shown in \cite[p. 163]{Deligne} that the varieties
attached to two such decompositions of $\bb$ are canonically isomorphic.
The projective limit of this system  of isomorphisms defines what we call
the  Deligne-Lusztig variety  $\bX(\bb)$ attached  to $\bb$;  it is  the
``usual'' Deligne-Lusztig variety $\bX(w)$ when we take $\bb=\bw\in\bW$.

When $W$ is an  irreducible Coxeter group, the center of  the pure braid group
is cyclic.  We denote by  $\bpi$ its positive  generator; we define  $\bpi$ in
general  as the  product of  the  corresponding elements  for the  irreducible
components of  $W$. Another way of  constructing $\bpi$ is as  $\bw_0^2$ where
$w_0$  is  the  longest  element  of  $W$.  The  Frobenius  endomorphism  acts
naturally as a diagram automorphism,  \ie, an automorphism which preserves $S$
(resp.  $\bS$), on  $W$ (resp  $B^+$); we  still denote  by $F$  these diagram
automorphisms. We call ``$F$-root of order $d$ of $\bpi$'', an element $\bb\in
B^+$  such that  $(\bb  F)^d=\bpi F^d$;  in \cite{Sydney}  it  is proved  that
$\beta(\bb)$ is  then a  regular element of  the coset $WF$  (in the  sense of
\cite{Springer}) for the eigenvalue $e^{2i\pi/d}$;  when $F$ acts trivially we
just have a  root of order $d$,  \ie, $\bb^d=\bpi$. It is  conjectured in 
\loccit\ that  the $\GF$-endomorphisms  of the  $\ell$-adic cohomology 
complex of
$\bX(\bb)$  form  a  ``cyclotomic  Hecke algebra''  attached  to  the  complex
reflection group $C_W(\beta(\bb)F)$. We will  show below that for any $F$-root
$\bb$ of $\bpi$ except for $\bpi$  itself, there is another $F$-root $\bw\in
\bW$  of  $\bpi$  of the  same  order  and  an  equivalence of  \'etale  sites
$\bX(\bb)\simeq\bX(\bw)$, thus the conjecture is about the variety $\bX(\bpi)$
and some ``ordinary'' Deligne-Lusztig varieties.

We  shall make  more specific  this conjecture  by replacing  it by  a set  of
conjectures  that we  shall study,  and prove  in some  specific examples.  In
\cite{DMR}, of  which this paper is  a continuation, we already  obtained some
general results on  some of the conjectures. We will  get here further results
for the element $\bpi$, and for all roots  of $\bpi$ in the case
of split groups  of type $A$. We  also study powers of  Coxeter elements in
type $B$ and fourth roots of $\bpi$ in split type $D$.
\vfill\eject
\section{Conjectures}\label{conjectures}

First, we  should state that  a guide  for the following  conjectures is
that, using  Lusztig's results in  \cite{LuCox}, we will show  in section
\ref{coxeter} that they all hold in the case of Coxeter elements.

We recall from \cite[2.1.1]{DMR} that a possible presentation of $B$ is
$$\langle \bw\in\bW\mid \bw_1\bw_2=\bw_3\text{ when } w_1w_2=w_3\text{ and }
l(w_1)+l(w_2)=l(w_3)\rangle.$$
We  recall  the  action  defined  in  \cite{DMR}  of  a
submonoid  of  $B^+$  on  $\bX(\bw)$,  and of  the  group  it  generates
(which  will  be  equal  to  $C_B(\bw F)$  in  every  case  we  study)  on
$H^*_c(\bX(\bw))$.  

First,     we    recall     the    definition     of    the     morphism
$D_\bt:\bX(\bb)\to\bX(\bt\inv   \bb   F(\bt))$    defined   when   $\bt$
is    a    left    divisor    of   $\bb$:    if    $\bb=\bt\bt'$,    and
if    $\bt=\bw_1\ldots\bw_n$   and    $\bt'=\bw'_1\ldots\bw'_{n'}$   are
decompositions   as   products   of   elements  of   $\bW$,   it   sends
the  element  $(\bB_1,\ldots,\bB_{n+n'+1})\in\bX(\bb)$  to  the  element
$(\bB_{n+1},\ldots,\bB_{n+n'},F(\bB_1),\ldots,F(\bB_{n+1}))$.

Then   we   introduce   categories    as   in   \cite{DMR}:   $\CB$   is
a   category   with   objects   the    elements   of   $B$,   and   such
that    $\Hom_{\CB}(\bb,\bb')=\{\by\in   B|\bb'=\by\inv\bb    F(\by)\}$;
composition  of  maps  corresponds  to  the  product  in  $B$;  one  has
$\End_\CB(\bb)=C_B(\bb F)$.

$\cD^+$  is the  smallest  subcategory  of $\CB$  which  contains the  objects
in  $B^+$  and  such  that  $$\{\by\in  B^+\mid\by\preccurlyeq\bb,  \by\inv\bb
F(\by)=\bb'\}  \subset\Hom_{\cD^+}(\bb,\bb'),$$  where $\preccurlyeq$  denotes
left divisibility in  the braid monoid, and $\cD$ is  the smallest subcategory
of $\CB$ containing $\cD^+$ and where all maps are invertible.

$\CC^+$ is the category of quasi-projective varieties on $\bbF$, together with
proper  morphisms.  $\CC$ is  the  localized  category by  morphisms  inducing
equivalences  of \'etale  sites.  An  isomorphism in  $\CC$  induces a  linear
isomorphism in $l$-adic cohomology.

It  is shown  in \cite{DMR}  that  the map  which  sends the  object $\bb$  to
$\bX(\bb)$ and the map $\bt$ to  $D_\bt$ extends to a functor $\cD^+\to\CC^+$,
which itself extends to a functor $\cD\to\CC$.

In  the following  we note  $H^*_c(\bX)$ for  the $\ell$-adic  cohomology with
compact support of the quasi-projective variety $\bX$. With this notation, the
monoid $\End_{\cD^+}(\bb)$  acts on $\bX(\bb)$  as a monoid  of endomorphisms,
and the group $\End_\cD(\bb)$ acts linearly on $H^*_c(\bX(\bb))$.

\begin{conjecture}\label{A}
When     $\bb$      is     an      $F$-root     of     $\bpi$      we     have
$\End_\cD(\bb)=\End_\CB(\bb)=C_B(\bb F)$.
\end{conjecture}
 
We will show this conjecture for  $\bpi$, $\bw_0$, Coxeter elements, all roots
of $\bpi$  in types $A$ and  $B$, and 4-th roots  of $\bpi$ in type  $D_4$. We
should note that in \cite[5.2.5]{DMR} we have defined an action of $C_{B^+(\bb
F)}$ on $\bX(\bb)$ which extends the action of $\End_{\cD^+}(\bb)$, but we are
not able to determine its  image in $H^*_c(\bX(\bb))$ (except via conjecture
\ref{A}).

It is proved in \cite[6.8]{Sydney}  that, except when $\bb=\bpi$, of course,
there is  a morphism  in $\cD^+$ between  any $d$-th $F$-root $\bb$  of $\bpi$  and a
``good''  $d$-th $F$-root $\bw$, ``good'' meaning that
$(\bw F)^i\in\bW.F^i$ for $i\leq d/2$.  Thus the  variety
$\bX(\bb)$ is  isomorphic in $\CC$  to $\bX(\bw)$,  as was asserted  above. We
thus  need only  to  consider a  variety  associated to  a  good $F$-root.  We
actually need only to consider one of them, according to the

\begin{conjecture}\label{conj}
There is always a morphism  in $\cD^+$ between
any  two $F$-roots  of  $\bpi$  of the  same  order.
\end{conjecture}

This says  in particular  that two  such roots are  $F$-conjugate in  $B$. The
result  of \cite{Sydney}  shows that  it  is sufficient  to consider  ``good''
$F$-roots of $\bpi$ in the above conjecture.

We will show this conjecture  for $\bw_0$, for Coxeter elements in split
groups and $n$-th roots of $\bpi$ in split type $A_n$. It has now been
proved in split type $A$ as a  consequence of a recent result of Birman,
Gebhardt and Gonzales-Meneses (personal communication) which states that
in general there  is a morphism in $\cD^+$ between  two conjugate roots,
and of  the theorem of  Eilenberg \cite{Eilenberg} stating that  in type
$A$ two roots of the same order are conjugate.

Since, when  $\bw$ is  an $F$-root  of $\bpi$,  $wF=\beta(\bw)F$ is  a regular
element of $WF$ (\cf\ \cite[6.6]{Sydney}),  the group $C_W(wF)$ is naturally a
complex reflection  group. We  will denote by  $B(w)$ the  corresponding braid
group; it  is shown in \eg,  \cite{BDM} in the case $F=\Id$
and in \S\ref{regular} below in the general case  that there
is  a  natural  map  $\gamma:B(w)\to  C_B(\bw  F)$  such  that  the  image  of
$\beta\circ\gamma$ is $C_W(wF)$. We recall the following conjecture
from \cite[0.1]{BDM}
\begin{conjecture}\label{B}
$\gamma$ is an isomorphism.
\end{conjecture}

The  above conjecture  is easy  when $\bw=\bpi$  or $\bw=\bw_0$.  It has  been
proved in \cite{BDM} for split types $A$ and $B$. We prove it for Coxeter
elements in split  groups, and for $4$th roots of  $\bpi$ in type $D_4$.
For this last case we  use programs of N.~Franco and J.~Gonzales-Meneses
which compute centralizers in Garside groups.

Assuming conjecture \ref{A}, and since  the operators $D_\bt$ commute with the
action  of  $\GF$, we  get  an  action  of  $B(w)$ as  $\GF$-endomorphisms  of
$H^*_c(\bX(\bw))$. The next conjecture states that this action factors through
a cyclotomic  Hecke algebra for  $B(w)$. Let  us recall their  definition; the
braid group  $B(w)$ is generated  by so-called ``braid reflections''  (see
\ref{braid reflections}
for the definition)  which form conjugacy classes in  bijection with conjugacy
classes  of distinguished  reflections in  $C_W(wF)$  (see again
\ref{braid reflections} for  the
definition) (\cf\ \cite[2.15  and Appendix 1]{BMR}). For  a representative $s$
of a  conjugacy class of  distinguished reflections  in $C_W(wF)$ we  choose a
representative  $\bs$ of  the corresponding  class of  braid reflections,  and
we  denote  by $e_s$  the  order  of  $s$.  Let  $A=\Qlbar[u_{s,j}]_{s,j}$  where
$\{u_{s,j}\}_{\{s\in S,  j=0\ldots e_s-1\}}$  are indeterminates.  The generic
Hecke algebra  of $C_W(wF)$ over  the ring $A$ is  defined as the  quotient of
$A[B(w)]$  by  the  ideal generated  by
$(\bs-u_{s,0})\ldots(\bs-u_{s,e_{s-1}})$.
A   $d$-cyclotomic   Hecke  algebra   for   $C_W(wF)$   is  a   ``one-variable
specialization'' of the generic algebra which specializes to $\Qlbar[C_W(wF)]$
by the further  specialization of the variable to $e^{2i\pi/d}$.  To make this
precise, we need some  definitions: we choose an integer $a$  and we denote by
$e^{\frac{2i\pi}{a|W|}}$ a primitive  $a|W|$-th root of unity  in $\Qlbar$; we
choose  an indeterminate  denoted by  $x^{1/a}$. Then  a $d$-cyclotomic  Hecke
algebra  is  a  specialization  of  the generic  Hecke  algebra  of  the  form
$u_{s,j}\mapsto  e^{2i\pi j/e_s}  (e^{-2i\pi/ad} x^{1/a})^{n_{s,j}}$  for some
integers $n_{s,j}$  (it is defined  over $\Qlbar[x^{1/a}]$; it  specializes to
$\Qlbar[C_W(wF)]$ by  the further  specialization 
$x^{1/a}\mapsto e^{2i\pi/ad}$).

\begin{conjecture}\label{C}
The action  of $B(w)$  on $H^*_c(\bX(\bw))$  factors through  a specialization
$x\mapsto q$ of a cyclotomic Hecke algebra $\CH(w)$ for $C_W(wF)$.
\end{conjecture}

More  precisely  we  have  to  state  the  specialization  as  $x^{1/a}\mapsto
q^{1/a}$.   This  conjecture   is  proved   for  $\bw=\bw_0$   and  $\bw=\bpi$
in   \cite[5.4.1]{DMR}   and   \cite[2.7]{Sydney}   respectively   (see   also
\cite[5.3.4]{DMR}). We will  prove this conjecture for all roots  of $\bpi$ in
split type $A$ for roots of even order in type $B$
and for $4$-th roots of $\bpi$ in type $D_4$.

Assuming  conjecture  \ref{C},  let  $\CH_q(w)$ be  the  above  specialization
(the    specialization   for    $x^{1/a}\mapsto    q^{1/a}$   of    $\CH(w)$).
We   thus    have   a   virtual   representation    $\rho_w$   of   $\CH_q(w)$
on   $\sum_i(-1)^i    H^i_c(\bX(\bw))$.   If   we    decompose   $\sum_i(-1)^i
H^i_c(\bX(\bw))=\sum_{\lambda\in\Irr(\GF)}    a_\lambda   \lambda$    in   the
Grothendieck  group  of $\GF$,  we  get  thus  for  each $\lambda$  a  virtual
character $\chi_\lambda$ of $\CH_q(w)$ of dimension $a_\lambda$.

We call  a representation  special if  its trace  defines up  to a  scalar the
canonical symmetrizing trace form on $\CH_q(w)$ (see \cite[2.1]{spetses} for the definition
of the canonical trace form). The canonical  trace form has not been proved to
exist for all complex reflection groups; however it is known to exist for those
groups that we will encounter in the present paper.

\begin{conjecture}\label{D}
\begin{enumerate}
\item The  $\chi_\lambda$ generate the Grothendieck  group of $\CH_q(w)$
and are irreducible up to sign.
\item The representation  $\rho_w$ is special.
\end{enumerate}
\end{conjecture}

(i)  above means  that  the image  of  $\CH_q(w)$ by  $\rho_w$  is the  ``full
$\GF$-endomorphism  algebra of  $\sum_i(-1)^i H^i_c(\bX(\bw))$''.  We will  be
able to  prove conjecture \ref{D} when  $\bw=\bpi$ and $\bG$ is  split of type
$A_n$, $G_2$, $E_6$ and  some small rank cases and also for  the cases when we
can prove the next conjecture.

\begin{conjecture}\label{E}
The groups $H^i_c(\bX(\bw))$ are disjoint from each other as $\GF$-modules.
\end{conjecture}

Conjectures      \ref{D}      and       \ref{E}      thus      imply      that
$\CH_q(w)\simeq\End_\GF(H^*_c(\bX(\bw)))$. Conjecture  \ref{E} is  the hardest
in  some sense,  since  it is  very difficult  to  determine individually  the
cohomology groups of a Deligne-Lusztig  variety, except when $\bw$ is rather
short. In addition to the known case of Coxeter elements, we will show \ref{E}
for $n$-th roots  of $\bpi$ in type  $A_n$ and $4$-th roots of  $\bpi$ in type
$D_4$. Also, when $\bG$ is of rank 2, conjecture \ref{E} follows from
\cite[4.2.4, 4.2.9, 4.3.4, 4.4.3 and 4.4.4]{DMR} (with some indeterminacy left
in type split $G_2$).

It should be pointed out that  conjectures \ref{C} to \ref{E} are consequences
of  a  special  case of  the  version  for  reductive  groups of  the  Brou\'e
conjectures on blocks  with abelian defect. They already  have been formulated
in  a  very  similar form  by Brou\'e, see \cite{Broue},
\cite{BrMa}  and \cite{Sydney}.  In
particular, compared to \cite[5.7]{Sydney}, we have only inverted the order of
the assertions,  in order to present  them by order of  increasing difficulty,
and made the  connection with the braid  group a little bit  more specific via
statements \ref{A} to \ref{B}.

\vfill\eject
\newcommand\Vr{{V^{\text{reg}}}}
\newcommand\VrH{{V^{{\text{reg}}_H}}}
\newcommand\Vzr{{V^{\text reg}_\zeta}}
\newcommand\VzrH{{V^{{\text reg}_H}_\zeta}}

\section{Regular elements in braid groups}\label{regular}

In this section we provide the needed background on braid monoids and groups.

We fix a  complex  vector  space $V$ of finite dimension.
A {\it complex reflection} is
an element of finite order of $\GL(V)$ whose fixed point space is a hyperplane.
Let  $W\in  \GL(V)$  be a  group generated by complex reflections, let $\CA$
be the set of reflecting hyperplanes for reflections of $W$, and let
$\Vr=V-\bigcup_{H\in\CA} H$.
Choose $x\in \Vr$ and let $\overline x\in\Vr/W$ be its image.
The group $\Pi_1(\Vr,x)$ is called the pure braid group of $W$ and
the group $B=\Pi_1(\Vr/W,\overline x)$ is called the braid group of $W$. 
The  map $\Vr\to\Vr/W$ is an \'etale 
covering, so  that it gives  rise to the exact  sequence $$1\rightarrow 
\Pi_1(\Vr,x)\rightarrow B\xrightarrow\beta W\rightarrow 1.$$
We denote by $N_W(E)$ (resp. $C_W(E)$) the stabilizer (resp. pointwise
stabilizer) of a subset $E\subset V$ in $W$.

We choose distinguished generators of $W$ and $B$ as follows
(see \cite[2.15]{BMR}):
\begin{definition}\label{braid reflections}
\begin{itemize}
\item[(i)]
A reflection $s\in W$ of hyperplane $H$ is distinguished if its only non
trivial eigenvalue is $e^{2i\pi/e_H}$ where $e_H=|C_W(H)|$.
\item[(ii)]
Let $H\in\CA$ and let $s_H$ be a distinguished reflection
of hyperplane $H$.
We call ``braid reflection'' associated to $s_H$ an element of $B$ of
the form $\gamma\circ\overline\lambda\circ\gamma\inv$, where:
$\gamma$ is a path from $\overline x$ to a point
$\overline{x_H}$ which is the image of a point $x_H\in\Vr$ ``close to $H$''
in the sense that there is a ball around
$x_H$ which meets $H$ and no other
hyperplane, and contains the path 
$\lambda:t\mapsto \proj_H(x_H)+e^{2i\pi t/e_H}\proj_{H^\perp}(x_H)$
where $\proj$ means orthogonal projection; and where $\overline\lambda$ is the
image of $\lambda$.
\end{itemize}
\end{definition}
It is clear from the definition that if $\bs_H$ is a braid reflection 
associated to $s_H$, then $\beta(\bs_H)=s_H$. Moreover it is proved in
\cite[2.8]{BMR} that braid reflections generate $B$ and in 
\cite[2.14]{BMR} that $\beta$
induces a bijection from the conjugacy classes of
braid reflections in $B$ to the conjugacy classes of distinguished reflections
in $W$.

We  assume now given $\phi\in\GL(V)$
normalizing  $W$;  thus $\Vr$  is  $\phi$-stable. 
We fix a  $d$-regular  element $w\phi\in W.\phi$
\ie,  an element  which has  an eigenvector  in $\Vr$ for the eigenvalue
$\zeta=e^{2i\pi/d}$. We refer to  \cite{Springer} \and \cite[\S B]{Sydney}
for the properties of such elements.
If $V_\zeta$ is the $\zeta$-eigenspace of $w\phi$, then 
$C_W(w\phi)=N_W(V_\zeta)$; we denote this group by  $W_\zeta$.
Its representation on $V_\zeta$ is faithful and makes it into a complex
reflection group, with reflecting hyperplanes the traces
on $V_\zeta$ of the reflecting hyperplanes of $W$.

To construct the braid  group $B(w)$
of $W_\zeta$ we choose a base point  $x_\zeta$ in $\Vzr=\Vr\cap
V_\zeta$. Then $B(w)=\Pi_1(\Vzr/W_\zeta,\overline x_\zeta)$, where  
$\overline  x_\zeta$  is the  image  of  $x_\zeta$  in  $\Vzr/W_\zeta$.
Let $\gamma:B(w)\to B$ be the morphism induced by the
injection $\Vzr/W_\zeta\hookrightarrow \Vr/W_\zeta$ composed with 
the quotient $\Vr/W_\zeta\to\Vr/W$.
It is shown in \cite[1.2 (ii)]{Bessis} that $\gamma$ factorizes through
a (unique) homeomorphism $\Vzr/W_\zeta\simeq (\Vr/W)^\zeta$,
where $(\Vr/W)^\zeta$
is the set of fixed points under the multiplication by $\zeta$.

Let $\delta$ be the path $t\mapsto e^{2i\pi  t/d}x_\zeta$ from $x_\zeta$
to $\zeta x_\zeta$ in $\Vzr$ and
let $\overline\delta$ be its image in $\Vr/W$.
Then the map $\bw\bphi:\lambda\mapsto
\overline\delta\circ (w\phi)(\lambda)\circ \overline\delta\inv$ 
is a lift to $B$ of the action of $w\phi$ on $W$.
\begin{remark}\label{bw}
Note that we have not defined independently
$\bw$ and $\bphi$; if $\phi$ is 1-regular this can be done in the following way:
let $x$ be a fixed point of $\phi$, and choose a path
$\eta_0$ from $x$  to $x_\zeta$.
Let $\eta_1$ be the path $\phi(\eta_0)\inv\circ\eta_0$,
from $\phi(x_\zeta)$ to $x_\zeta$.
The image in $\Vr/W$ of the path $\delta\circ w(\eta_1)$
from $x_\zeta$ to $w(x_\zeta)$ is an element $\bw$ of
$B$ which by definition lifts $w\in W$.

Then  if $\overline\eta_1$ is the image of $\eta_1$ in $\Vr/W$,
the map $\bphi:\lambda\mapsto\overline\eta_1\inv\circ\phi(\lambda)\circ
\overline\eta_1$ is an automorphism of $B$ which lifts
the action of $\phi$ on $W$.
\end{remark}
Note also the following:
\begin{remark}
Let $\bphi$ and $\bw$ be as in \ref{bw}, and
let $\bpi$ be the loop $t\mapsto e^{2i\pi t}x_\zeta$ in $\Vr$;
it is a generator of the
center of the pure braid group.
Then $(\bw\bphi)^d=\bpi\bphi^d$ in the semi-direct product
$B\rtimes<\bphi>$. Indeed
the element $(\bw\bphi)^d\bphi^{-d}$ is represented by the path
$\delta\circ w\phi(\delta)\circ(w\phi)^2
(\delta)\circ\cdots\circ(w\phi)^{d-1}(\delta)$
which is equal to $\bpi$.
\end{remark}

\begin{lemma}\label{C_B(w.phi)}
We have $\gamma(B(w))\subset C_B(\bw\bphi)$.  
\end{lemma}
\begin{proof}
It   is   easily   checked    that   for   any $\lambda\in B$ 
the   path    $\overline\delta\circ\lambda\circ\overline\delta\inv$    is 
homotopic  to   $\zeta\inv\lambda$,   so $(\bw\bphi)(\lambda)$
is  homotopic  to $\zeta\inv (w\phi)(\lambda)$. 
If $\lambda$ is the image of a path in $\Vzr$, we have 
$\zeta\inv(w\phi)(\lambda)=\lambda$; so the action of $\bw\bphi$ on 
the image of $\lambda$ in $B$ is trivial, as claimed.                    
\end{proof}

We now  assume that $W$  is a Coxeter group.
The  space $V$  is  the complexification of the  real reflection
representation of  $W$ and the real hyperplanes define chambers. We choose
a fundamental chamber, which defines a Coxeter generating set $S$ for $W$.

When $\phi$ is trivial, and $W$ not of type $F_4$ or $E_n$,
the  morphism $\gamma$ has been proved injective
\cite[4.1]{Bessis}; if moreover $W$ is of type 
$A$(split)  or $B$ it is proved to be an isomorphism onto
$C_B(\bw\bphi)$ in \cite{BDM}.

We assume now that $\phi$  induces a diagram  automorphism  of $W$,
\ie, that it stabilizes the fundamental chamber. Then
$\phi$ is 1-regular: in fact it has a fixed point
in the fundamental chamber.
We recall the following result from \cite{VdL}:
\begin{proposition}\label{VdL}
Assume that we have fixed a base point whose real part is in the fundamental
chamber.
Let $\bW$ be the set of elements of $B$ which
can be represented 
by paths $\lambda$  in $\Vr$, starting from the base point and
satisfying the two following properties:
\begin{itemize}
\item[(i)] The real part of $\lambda$
meets each element of $\CA$ at most once.
\item[(ii)]
When the  real part of $\lambda$ 
meets  $H\in\CA$,
the imaginary  part of $\lambda$ is  on the 
same side of $H$ as the fundamental chamber;
\end{itemize}
then $\bW$ is in bijection with $W$ via the map $B\xrightarrow\beta
W$. Moreover if $\bS\subset \bW$ is such that $\beta(\bS)=S$
then $B$ has a presentation with generators $\bS$ and relations the
braid relations given by the Coxeter diagram of $W$.
\end{proposition}
The elements of $\bS$ are braid reflections.
\begin{corollary}\label{canonical isomorphism}
Let $x$ and $y$ be two points with real parts in the fundamental chamber.
If $\gamma$ is a path from $x$ to $y$ with real part
in the fundamental chamber, the isomorphism
$\Pi_1(\Vr/W,x)\to \Pi_1(\Vr/W,y)$ which it defines
is independent of $\gamma$.
\end{corollary}
\begin{proof}
Two such isomorphisms differ by the inner automorphism of $\Pi_1(\Vr/W,x)$
defined by a loop with real part in the fundamental chamber. An element
defined by a loop is in the pure
braid group. But by the proposition \ref{VdL}
this element is  also in $\bW$,  so it is trivial. 
\end{proof}
\begin{definition}
The braid  monoid $B^+$  is  defined to  be the  submonoid of  $B$ 
generated by $\bW$.  Its elements are called positive  braids.
\end{definition}

Recall that a $d$-regular element in $W.\phi$ is called a
Springer element if it
has a $\zeta$-eigenvector with real part
in the fundamental  chamber (\cf\ \cite[3.10]{Sydney}).
We choose as base point a fixed point of $\phi$ in the fundamental chamber.
If $w\phi$ is a Springer element we can choose $\eta_0$ in \ref{bw} such that
$\Re(\eta_0)$ is  in the fundamental chamber.
\begin{proposition}\label{bw in bW}
Assume that $d>1$ and  that $w\phi$ is a Springer element. 
If $\eta_0$ in \ref{bw}
is chosen with real part in the fundamental
chamber then the element $\bw$ is independent of the
choice of $\eta_0$ and is in $\bW$.
\end{proposition}
\begin{proof}
The element 
$\bw$ is  represented by  $\delta\circ w(\eta_1)$.  As $\Re(w(\eta_1))$ 
does not meet  any reflecting hyperplane, the only  possibility for the 
real part  of this  path to  meet $H\in\CA$  is that 
for  some $t\in]0,1[$  we have  $l_H(\Re( \delta(t)))=0$, where
$l_H$ is a real linear form defining  $H$.
Write $x_\zeta=a+ib$ 
with  $a$ and  $b$ real  and let  $\theta=2\pi t/d$; since $d>1$ we have
$\theta\in]0,\pi[$.  We have 
$l_H(\Re(\delta(t)))=(\cos\theta) l_H(a)-(\sin\theta) l_H(b)=0$.
As   this  equation  has    only   one    solution 
in $]0,\pi[$, property \ref{VdL}  (i) above is satisfied.
Moreover if $\theta$ is such a solution we have
$l_H(\Im(\delta(t)))=l_H(a)(\sin\theta+\dfrac{\cos^2\theta}{\sin\theta})
=\dfrac{l_H(a)}{\sin\theta}$  which   has the same sign as $l_H(a)$
since $\sin\theta>0$.   So  property \ref{VdL}
(ii) is also satisfied.
\end{proof}

When $d=1$  we still have that  similarly $\bpi$ is independent  of the 
choice of $\eta_0$  and is in $B^+$; a  way to see this is  to use that 
$\bpi=\bw_0^2$ where $\bw_0$ is the element $\bw$ obtained for $d=2$.   

We show now how in some cases we  can lift to $B$ a distinguished reflection of 
$C_W(w\phi)$. We still assume that $w\phi$ is a Springer element

First,  to  $H\in\CA$  we  associate
the hyperplane $H\cap V_\zeta$ of $V_\zeta$  and the 
distinguished reflection $t_H$
of  $W_\zeta$  with  reflecting  hyperplane  $H\cap  V_\zeta$.
Let $W_H=C_W(H\cap V_\zeta)$; it is
a  parabolic  subgroup of  $W$, thus a reflection subgroup.  
The  element $w\phi$  normalizes  $W_H$
as  it  acts  by  $\zeta$  on  $H\cap V_\zeta$;  it  is  a  regular
element  of  $W_H.w\phi$ and  we have $C_{W_H}(w\phi)=<t_H>$.
We   can  apply  the constructions of this section
to  $W_H$: 
let  $\VrH=V-\bigcup_{\{H'\in\CA\mid H'\supset H\cap V_\zeta\}}H'$  and  
let  $\VzrH=V_\zeta\cap\VrH$.
We  get  a  morphism
$\Pi_1(\VzrH/C_{W_H}(w\phi),\overline  x_\zeta)\to \Pi_1(\VrH/W_H,\overline
x_\zeta)$,
whose image centralizes
$\bw\bphi$ as in lemma \ref{C_B(w.phi)}.
As  $x_\zeta$  and $x$  are  both  in  the
fundamental chamber of $W_H$, by \ref{canonical isomorphism}
any path from $x$ to $x_\zeta$ whose real part
stays in this fundamental chamber defines a canonical isomorphism between
$\Pi_1(\VrH/W_H,\overline  x_\zeta)$ and $\Pi_1(\VrH/W_H,\overline  x)$ which
commutes with $\bw\bphi$. Let us denote by $B_H$ this group.
By composition with this isomorphism we get a morphism
$\Pi_1(\VzrH/C_{W_H}(w\phi),\overline x_\zeta)\to B_H$,
whose  image  centralizes $\bw\bphi$.  Let  $\bt_H$  be  the generator of the 
infinite  cyclic group  $\Pi_1(\VzrH/C_{W_H}(w\phi))$  such that  its image  in 
$C_W(w\phi)$ is $t_H$; then $\bt_H$ is a braid reflection in $B(w)$ and, if
$e_H$ is the order of $t_H$, then 
$\bt_H^{e_H}$ is the loop $t\mapsto e^{2i\pi t}x_\zeta$ \ie, the element $\pi_H$
of $B_H$.
If  $H$ is  such that  $W_H$ is  a 
standard parabolic  subgroup, \ie,  is generated  by a subset  $I$ of  $S$ then 
$B_H=\Pi_1(\VrH/W_H,\overline x)$ is canonically embedded in $B$
as the subgroup  of $B$ generated by the lift 
$\bI\subset\bS$ of $I$ and if  $\bs_H\in  B_H$ is
the  image of $\bt_H$  by  the above  morphism
we have $\bs_H\in B_H\cap C_B(\bw\bphi)$.                    

\vfill\eject
\section{The case of Coxeter elements}\label{coxeter}

We prove here,  using the results of \cite{BDM},  \cite{Dual} and \cite{LuCox}
that all our conjectures hold in the case of Coxeter elements for an untwisted
quasi-simple reductive group $\bG$ (the assumption of $\bG$ being quasi-simple
is equivalent  to $W$ being irreducible).  Even though the results  of Lusztig
cover them, we are unable to handle twisted groups because the construction of
a dual braid monoid when $F$ is not trivial has not yet been carried out.

Let $h$ be the  Coxeter number of $\bG$ and denote by  $n$ the semisimple rank
of $\bG$, which is also the Coxeter rank of $W$. We begin with

\begin{proposition}\label{h-th root of pi}
The $h$-th roots of  $\bpi$ are the lift to $\bW$ of  Coxeter elements of $W$.
Conjecture \ref{conj}  holds for  $h$-th roots  of $\bpi$,  that is,  there is
always a morphism in $\cD^+$ between two such roots.
\end{proposition}
\begin{proof}
By \eg,  \cite[3.11]{Sydney} $h$-th roots of  $\bpi$ exist. Such a  root is an
element of $B^+$ of  length $n$, whose image in $W$ is  in the conjugacy class
of Coxeter  elements by  \cite[3.12]{Sydney}. Since the  minimal length  of an
element in this conjugacy class is  $n$, which is attained exactly for Coxeter
elements, we conclude that an $h$-th root of $\bpi$ is in $\bW$, and its image
is a Coxeter element.

Now, by  \cite[chap. V \S 6, Proposition 1]{Bbki}, any two  such elements are
connected  by a
morphism in  $\cD^+$ (it is easy  to identify the conjugating  process used in
\loccit\ to morphisms in $\cD^+$).
\end{proof}

We now show
\begin{proposition} \label{centralizercox} Let $\bc$ be the lift in $\bW$ of a
Coxeter element. Then $C_B(\bc)$ is the cyclic group generated by $\bc$.
\end{proposition}
\begin{proof}
Here we use the results  of \cite{BDM} and \cite{Dual}. By \cite[2.3.2]{Dual},
$B$  admits  a  Garside  structure  where  $\bc$  is  a  Garside  element.  By
\cite[2.26]{BDM}, the fixed points of $\bc$ are generated by the lcm of orbits
of atoms under $\bc$. Such an element  is a $\bc$-stable simple element of the
dual  braid monoid.  By \cite[1.4.3]{Dual},  it  is the  lift to  $\bW$ of  an
element $c_1$  in the Coxeter  class of a parabolic  subgroup of $W$.  But the
centralizer  of  $c$ in  $W$  is  the cyclic  subgroup  generated  by $c$,  in
particular a simple  element centralizes $\bc$ only  if its image in  $W$ is a
power of $c$. Thus, we have to show that no power $c^k$ of $c$ with $1<k<h$ is
the image of a simple element, or  equivalently that no such power divides $c$
for the  reflection length. By  \cite[1.2.1]{Dual} this is equivalent  to
showing
that the equality $\dim\ker(c^k-1)+\dim\ker(c^{1-k}-1)=n$ cannot hold for such
$k$. But, by \cite[4.2]{Springer}, the  eigenvalues of $c$ are $\zeta^{1-d_i}$
where  $\zeta=e^{2i\pi/h}$  and where  $d_i$  are  the reflection  degrees  of
$W$.  Thus  the equality  to  study  becomes $|\{i\mid  (1-d_i)k\equiv  0\pmod
h\}|+  |\{i\mid  (1-d_i)(1-k)\equiv  0\pmod h\}|=n$.  Both  conditions  cannot
occur simultaneously  since this  would imply $d_i\equiv  1\pmod h$,  which is
impossible since the irreducibility of $W$  implies that $1<d_i\le h$. Thus it
is  sufficient to  exhibit  a  $d_i$ which  satisfies  neither condition.  But
$d_i=h$ itself is such a $d_i$, whence the result.
\end{proof}

It follows immediately from \ref{centralizercox} that conjecture \ref{A} holds
for Coxeter  elements, that  is $\End_\cD(\bc)=C_B(\bc)=\genby \bc$,  since by
definition $\bc\in\End_\cD(\bc)$.

Let us now prove that conjecture \ref{B} holds.
The space $V_\zeta$ is one-dimensional, and for any $x_\zeta
\in\Vzr=V_\zeta-\{0\}$ the group $\Pi_1(V_\zeta/C_W(c),x_\zeta)$ is cyclic,
generated by the loop $\bb=t\mapsto e^{2i\pi t/h}$. Doing if necessary a
conjugation in $B$, we may take any $h$-th root of $\bpi$ to prove \ref{B}.
We will choose a Springer element, so we may assume that $\Re(x_\zeta)$ is
in the real fundamental chamber of $W$. Then, by \ref{bw in bW}, the image
$\gamma(\bb)$ is in $\bW$, and since its image in $W$ is $c$ it is equal to
$\bc$. We have shown that $\gamma$ is an isomorphism.

Conjectures \ref{D}  and \ref{E}  will follow  from the  following proposition
from \cite{LuCox}:

\begin{proposition}\label{lucox}
Let $c$ be a Coxeter element. Then
\begin{enumerate}
\item  $F$  is  a  semisimple automorphism  of  $\oplus_i  H^i_c(\bX(c))$;  it
has  $h$  distinct eigenvalues;  the  corresponding  eigenspaces are  mutually
non-isomorphic irreducible $\GF$-modules.
\item  For $s=1,\ldots,h-1$,  the endomorphism  $F^s$ has  no fixed  points on
$\bX(\bc)$.
\item   The    eigenvalues   of    $F$   are    monomials   in    $q$   which,
under    the   specialization    $q\mapsto    e^{2i\pi/h}$,   specialize    to
$1,\zeta,\zeta^2,\ldots,\zeta^{h-1}$ where $\zeta=e^{2i\pi/h}$.
\end{enumerate}
\end{proposition}
\begin{proof} (i) is \cite[6.1 (i)]{LuCox}; (ii) is 6.1.2 of \loccit;
(iii) results from the tables pages 146--147 of \loccit
\end{proof}

We  show now  how this  implies  conjectures \ref{D}  and \ref{E}.  Conjecture
\ref{E} is immediate  from \ref{lucox}(i). The generic  Hecke algebra $\CH(c)$
of the  cyclic group  $C_W(c)$ of order  $h$ is generated  by one  element $T$
with  the  relation  $(T-u_0)\ldots  (T-u_{h-1})=0$. The  map  which  sends  $\bc$
to  $D_\bc=F$  is  thus  a  representation of  this  algebra,  specialized  to
$u_i\mapsto \lambda_i$ where  $\lambda_0,\ldots,\lambda_{h-1}$ are the eigenvalues
of  $F$ on  $\oplus_i H^i_c(\bX(c))$.  By \ref{lucox}(iii)  this is  indeed an
$h$-cyclotomic  algebra for  $C_W(c)$.  It  remains to  see  that the  virtual
representation $\sum_i  (-1)^i H^i_c(\bX(c))$  of $\CH(c)$  is special.  But, by
\eg, \cite[2.2]{BrMa}, the symmetrizing trace  on $\CH(c)$ is characterized by
its vanishing on $T^i$ for  $i=1,\ldots,h-1$. By the Lefschetz trace formula,
one has  $\TrH{F^s}{\bX(c)}  =|\bX(c)^{F^s}|$, so this
vanishing is a consequence of \ref{lucox}(ii).

\vfill\eject
\section{Regular elements in type $A$}\label{reg in A}
We prove conjecture \ref{A} for roots of $\bpi$
when $W$ is of type $A$ and $F$ acts trivially on $W$. Here we assume that
$d>1$. The case of $\bpi$ will be treated in section \ref{secpi}.

Let $W$ be a Coxeter group of type $A_{n-1}$ and let $B$ be the associated braid
group.
There  are $d$-regular  elements  in $W$  for  any $d$
dividing $n$ or $n-1$. To handle the case $d|n-1$, it will be
simpler to  embed $W$ as  a parabolic subgroup of  a group $W'$  of type
$A_n$ and $B$ as a parabolic subgroup of the associated braid group $B'$
and to consider $d$-regular elements for $d|n$ in $W'$.

We denote by $\bS=\{\bsigma_1,\ldots,\bsigma_n\}$  the set of generators
of $B'$,  the generators of $B$  being $\bsigma_1,\ldots,\bsigma_{n-1}$.
We denote by $\bpi$ and  $\bpi'$ respectively the positive generators of
the centers of the pure braid  groups respectively associated to $W$ and
$W'$.

Let $\bc=\bsigma_1\bsigma_2\ldots\bsigma_{n-1}$, the
lift in $\Bred$ of a Coxeter element.
Let $r$ and $d$ be two integers such that $rd=n$ and let $\bw=\bc^r$; it is a
$d$-th root of $\bpi$ (\cf\ \ref{h-th root of pi}).
The image $w$ of $\bw$ in $W$ is a regular element of
order $d$ and its centralizer $C_W(w)$ is isomorphic
to the complex reflection group $G(d,1,r)$ which has a presentation
given by the diagram:
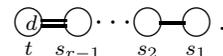
\def\nnode#1{{\kern 4.2pt\hbox to
0pt{\hss{$\mathop\bigcirc\limits_{#1}$}\hss}\kern 4.2pt}}
\def\ncnode#1#2{{\kern 4.2pt\hbox to
0pt{\hss{$\mathop\bigcirc\limits_{#1}\kern-9pt{\scriptstyle#2}
$}\hss}\kern 5.8pt}}
\def\bar{{\vrule width10pt height3pt depth-2pt}}
\def\dbar{{\rlap{\vrule width10pt height2pt depth-1pt} 
                 \vrule width10pt height4pt depth-3pt}}
$\ncnode td\dbar\nnode{s_{r-1}}\cdots\nnode{s_2}\bar\nnode{s_1}\kern 4pt$.

We denote by $B(d,1,r)$ the braid group associated to $G(d,1,r)$;
it has a presentation given by the same diagram (deleting the relations giving
the orders of the generators).

Conjecture \ref{B} holds in our case. Indeed,
Bessis  \cite[4.1]{Bessis} has proved that the morphism
$\gamma:B(w)\to  C_B(\bw)$  is injective and
in  \cite{BDM},   this morphism is proved to be bijective.
More precisely Bessis proves that
$\gamma(\bs_i)= \prod_{j=0}^{d-1}\bsigma_{i+rj}$
and  $\gamma(\prod_{i=1}^{r-1}\bs_i\bt)=\bc$, where
$\bt,\bs_1,\ldots,\bs_{r-1}$ are the generators of $B(w)$ given
in \cite[3.6]{BMR}: they are braid reflections which satisfy the braid
relations given by the above diagram for $G(d,1,r)$.
We shall identify $B(w)$ with its  image, so that we shall identify
$\bs_i$  and  $\bt$  with the elements given by the above formulas.

Let now $\bc'=\bsigma_1\ldots\bsigma_n\bsigma_n$, a $n$-th root of
$\bpi'$ in $B'$, \cf\ \cite[A1.1]{Sydney}.  Let $\bw'=\bc^{\prime r}$;
it is a $d$-th root of $\bpi'$ and if $w'$ is its image in $W'$, the
centralizer $C_{W'}(w')$ is also isomorphic to
$G(d,1,r)$. It has been proved in
\cite[5.2]{BDM} that  $C_B(\bw)\simeq C_{B'}(\bw')$.
More precisely, let $X_n$ be the configuration space of $n$ distinct
points in $\bbC$, let $\mu_n$ be the set of $n$-th roots of 1,
and let $\nu_{n+1}=\mu_n\cup\{0\}$.
We have $B'=\Pi_1(X_{n+1},\mu_{n+1})$. 
Let
$X^*_n$ be the configuration space of $n$ non-zero distinct points in
$\bbC$;
we have morphisms
$\Pi_1(X_{n+1},\nu_{n+1})\xleftarrow{\Psi}\Pi_1(X^*_n,\mu_n)
\xrightarrow{\Theta}\Pi_1(X_n,\mu_n)$,
(these morphisms are called $A$ and $B$ in \loccit).
One gets the map $\Psi$ 
by adjoining to a braid a constant string at $0$.

If we choose an isotopy from $\nu_{n+1}$  to $\mu_{n+1}$
we get an isomorphism of $\Pi_1(X_{n+1},\nu_{n+1})$ with $B'$;
this can be done by bringing $0$ along a
path ending at $e^{2i\pi\frac n{n+1}}$.

Similarly,  if we choose an
isotopy mapping the $n$ first $n+1$-th roots of 1 to $\mu_n$,
we get an isomorphism $\alpha: \Pi_1(X_n,\mu_n)\xrightarrow\sim B$.

It is shown in \loccit\ 
that the map $\alpha\circ\Theta$ has a section $\Theta'$ above $C_B(\bw)$.

Let
$\psi$ be the restriction to $C_B(\bw)$ of $\Psi\circ\Theta'$;
it is an isomorphism from $C_B(\bw)$ to $C_{B'}(\bw')$.
We have $\psi(\bc)=  \bc'$ and
$\psi(\bs_i)=\bs_i$. Let  $\bt'=\psi(\bt)$; it  satisfies
$\prod_{i=1}^{r-1}\bs_i\bt'=\bc'$.

The following theorem proves conjecture \ref{A} in type $A$. Note that if we have a
parabolic subgroup $W_1$ of a Coxeter group $W_2$, and if $B_1$ and $B_2$ are
the corresponding braid groups, the category $\cD_1$ (resp. $\cD_1^+)$
associated to $\bB_1$ as in section \ref{intro}
is a full subcategory of the category $\cD_2$ (resp.
$\cD_2^+$) associated to $B_2$. This allows us in the following theorem to
state the results in term of the categories associated to $B'$. We will denote
these categories by $\cD$ and $\cD^+$.

\begin{theorem}\label{c^i et c'^i dans A}

One has $\bs_i\in\End_{\cD^+}(\bw)\cap\End_{\cD^+}(\bw')$,
$\bt\in\End_\cD(\bw)$ and $\bt'\in\End_\cD(\bw')$, so that
$\End_\cD(\bw)=C_B(\bw)\simeq B(d,1,r)$ and
$\End_\cD(\bw')=C_{B'}(\bw')\simeq B(d,1,r)$.
\end{theorem}

The end of this section is devoted to the proof of that theorem.

In the next lemma $\preccurlyeq$
(resp. $\succcurlyeq$) denotes left divisibility (resp. right divisibility)
in the braid monoid.
\begin{lemma}\label{xy-1z} Let $B$ be the braid group of an arbitrary finite
Coxeter group;
let $\bx,\by, \bz\in\BW$ be such that
$\bx\by\inv\bz\in\BW$. Then there exist $\bx_1,\bz_1\in\BW$ such that
$\bx\succcurlyeq\bx_1$, $\bz_1\preccurlyeq\bz$ and $\by=\bz_1\bx_1$.
\end{lemma}
\begin{proof}
Let $\bb\in\BW$ be such that $\bx\by\inv\bz=\bb$, \ie, $\by\inv\bz=\bx\inv\bb$.
By \cite[3.2]{michel}, if we denote by $\bz_1$ the left gcd of
$\by$ and $\bz$, and by $\bx_2$ the left gcd of
$\bx$ and $\bb$, we have $\bz_1\inv\by=\bx_2\inv\bx$, whence the result,
putting $\bx_1=\bx_2\inv\bx$.
\end{proof}
\begin{lemma}\label{facts} For $i=1,\ldots,n$,
let $\bc_i=\bsigma_1\ldots\bsigma_i$.
\begin{itemize}
\item[(i)] We have
$\bc_k\bsigma_i=\bsigma_{i+1}\bc_k$ for $i<k$.
\item[(i')] We have
$\bc'\bsigma_i=\bsigma_{i+1}\bc'$ for $i<n-1$.
\item[(ii)]
We have $\lexp{\bc^2}\bsigma_{n-1}=\bsigma_1$.
\item[(ii')]
We have $\lexp{\bc^{\prime 2}}\bsigma_{n-1}=\bsigma_1$.
\item[(iii)]
For $\bx\in\BW$, one has
$\bsigma_{i+1}\preccurlyeq\bc_k\bx\Leftrightarrow\bsigma_i\preccurlyeq\bx$ for $i<k$.
\item[(iv)]  For $j\le k$, we have $\{i\mid\bsigma_i\preccurlyeq\bc_k^j\}=\{1,\ldots,j\}$.
\item[(iv')] For $j\le n$ we have
$\{i\mid\bsigma_i\preccurlyeq\bc^{\prime j}\}=\{1,\ldots,j\}$.
\item[(v)]
$\bc^{\prime j}=\bc_n^j\bsigma_{n-j+1}\ldots\bsigma_n$ for $1\leq j\leq n$.
\end{itemize}
\end{lemma}
\begin{proof}

Let us prove (i) and (i'). We get $\bc_k\bsigma_i =\bsigma_{i+1}
\bc_k$ by commuting $\bsigma_i$ (resp. $\bsigma_{i+1}$)
with the factors of  $\bc_k$ and applying once the braid relation between
$\bsigma_i$ and $\bsigma_{i+1}$.
Moreover, as $\bc'=\bc_n\bsigma_n$ and $\bsigma_n$
commutes with $\bsigma_i$ for $i<n-1$, we get also (i').

From (i), by induction on $j$ we get (v).

Let us prove  (ii).
For  proving $\lexp{\bc^2}\bsigma_{n-1}=\bsigma_1$,  we use  (i) to  get
$\bc\bsigma_1\ldots\bsigma_{n-2}=\bsigma_2\ldots\bsigma_{n-1}\bc$,  then
we multiply both  sides on the right by $\bsigma_{n-1}$  and on the left
by $\bsigma_1$ to get $\bsigma_1\bc^2=\bc^2\bsigma_{n-1}$.

We deduce (ii'): we have
$\lexp{\bc'}\bsigma_{n-1}=\lexp{\bc\bsigma_n}\bsigma_{n-1}
=\lexp{\bc\bsigma_{n-1}\inv}\bsigma_n=\lexp{\bsigma_n\inv\bc}\bsigma_n$,
so that $\lexp{\bc^{\prime
2}}\bsigma_{n-1}=\lexp{\bc\bsigma_n\bsigma_n\inv\bc}\bsigma_n=\bsigma_1$
by (ii).

Let us prove (iii). We have
$\bsigma_{i+1}\preccurlyeq\bc_k\bx\Leftrightarrow\bsigma_{i+1}\inv \bc_k\bx\in\BW$;
by (i) $\bsigma_{i+1}\inv\bc_k \bx=\bc_k
\bsigma_i\inv\bx$ and by lemma \ref{xy-1z}
this implies that either $\bsigma_i\preccurlyeq\bx$
or $\bc_k\succcurlyeq\bsigma_i$.
But, as no braid relation can be applied in $\bc_k$, the only
$j$ such that $\bc_k\succcurlyeq\bsigma_j$ is $k$. So we are in the case
$\bsigma_i\preccurlyeq\bx$, whence the implication from left to right.
The converse implication comes from the fact that
$\bsigma_i\preccurlyeq\bx\Rightarrow \bc_k\bsigma_i\inv\bx\in\BW$.

Let us prove (iv). 
If $\bsigma_i\preccurlyeq\bc_k^j$
we cannot have $i>j$  otherwise
applying  $j$ times  (iii), we get
$\bsigma_{i-j}\preccurlyeq  1$ which is false.
Conversely, if  $i\le  j$, by applying $i-1$ times
(iii) we get  $\bsigma_i\preccurlyeq\bc^j_k\Leftrightarrow \bsigma_1\preccurlyeq
\bc_k^{j-i+1}$ which is true.

(v) is a direct computation using (i).

Let us prove (iv').
By (iv) and (v), we have
$\bsigma_i\preccurlyeq\bc^{\prime j}$ for $i\leq j$. On the other hand
by (iii) and (v)
we see that if $\bsigma_i\preccurlyeq \bc^{\prime j}$  with $i>j$, then
$\bsigma_{i-j}\preccurlyeq\bsigma_{n-j+1}\ldots\bsigma_n$ which is impossible.
\end{proof}

In the next lemma, as in \cite[5.2.7]{DMR}, for $I\subset S$ we denote
by $B^+_I$ the submonoid of $B^+$ generated by $\bI=\{\bs\in\bS\mid s\in I\}$,
and by $\cD^+_I$ the ``parabolic'' subcategory of $\cD^+$ where we keep only
the maps coming from elements of $B^+_I$.

\begin{lemma}\label{s_i ok} 
Assume $1\le i<r$ and let
$\bI_i=\{\bsigma_i,\bsigma_{i+r},\ldots,\bsigma_{i+(d-1)r}\}$. Then
\begin{itemize}
\item[(i)] $\bs_i\in \End_{\cD^+_{I_i}}(\bw)\cap \End_{\cD^+_{I_i}}(\bw')$.
In particular $\bs_i\in\End_{\cD^+}(\bw) \cap\End_{\cD^+}(\bw')$.
\item[(ii)]  The conjugation by either $\bw$ or $\bw'$ stabilizes 
$\bI_i$ and induces on this set the cyclic permutation
$\bsigma_{i+jr}\mapsto\bsigma_{i+(j+1)r\pmod n}$.
\end{itemize}
\end{lemma}
\begin{proof}
We can assume $r\ge 2$ otherwise there is nothing to prove.
Let $\by_j=\bsigma_{i+r(j-1)}$ for $j=1,\ldots,d$, and,
following the notation of \cite[5.1.1 (i)]{DMR} let
$\bw_1=\bw$ and $\bw_{j+1}=\by_1\inv\bw\by_j$ for
$j=1,\ldots,d$. We have $\bw_j\in\BW$ as $\by_1=\bsigma_i$ divides
$\bw=\bc_{n-1}^r$ by \ref{facts}(iv).
We have $\bw_{j+1}=\by_j\inv\bw_j\by_j$ by using
$\lexp\bw\by_{j-1}=\by_j$, which is a consequence of
\ref{facts}(i), and the fact that
$\by_j$ and $\by_1$  commute as $r\ge 2$; 
from \ref{facts}(ii) we have $\lexp\bw\by_d=\by_1$,
whence $\bw_{d+1}=\bw$ so that (i) is proved for $\bw$. We get (i) for $\bw'$
by the same computation, replacing $\bw$ by $\bw'$ and using \ref{facts}
(i'), (ii') and (iv') instead of \ref{facts} (i), (ii) and (iv).
We have also got (ii) along the way.
\end{proof}

To prove  that $\bt\in\End_\cD(\bw)$, we shall find
$\by\in B^+$ such that $\by\bw\by\inv\in B^+$,
$\by\in\Hom_{\cD^+}(\by\bw\by\inv,\bw)$ and
$\by\bt\by\inv\in\End_{\cD^+}(\by\bw\by\inv)$. Then, as
$\cD$ is a groupoid, we will get $\bt\in\End_\cD(\bw)$. 

We will follow the same lines to prove that $\bt'\in\End_\cD(\bw')$.

\begin{lemma}\label{y in B<-w}
When $j\le k$, let
$\bsigma_{j,k}= \bsigma_j\bsigma_{j+1}\ldots\bsigma_k$
and let $\bx_i=\bsigma_{i,i+r-2}$.
With this notation, let
$\by=\prod_{i=1}^{d-1}\by_i$ where
$\by_i=\prod_{j=d}^{i+1}\bx_{i(r-1)+j}$. Then
$\by\bw\by\inv\in B^+$ and $\by\in \Hom_{\cD^+}(\by\bw\by\inv,\bw)$.
\end{lemma}
\begin{proof}
We set $\bw_d=\bw$, and then by decreasing induction on $i$ 
we define $\by_i\inv\bw_i\by_i=\bw_{i+1}$.
It is enough to show that
$\by_i\inv\bw_i\in\BW$; we will get this by proving by induction that
$\bw_i=\by_i\bc^{r-1}\bc_{i(r-1)+d-1}$.
This formula is true for $i=d$ (here $\by_d=1$).
Let us assume it true for $i+1$  and let us prove it for $i$. We have
\begin{align*}\bw_{i+1}&=\by_{i+1}\bc^{r-1}\bc_{(i+1)(r-1)+d-1}\cr
&=\bx_{(i+1)(r-1)+d}\ldots \bx_{(i+1)(r-1)+i+1}\bc^{r-1}\bc_{(i+1)(r-1)+d-1}\cr
&=\bc^{r-1}\bx_{i(r-1)+d}\ldots \bx_{i(r-1)+i+1}\bc_{(i+1)(r-1)+d-1}
\quad\text{(by \ref{facts}(i))}\cr
&=\bc^{r-1}\bc_{(i+1)(r-1)+d-1}\bx_{i(r-1)+d-1}\ldots \bx_{i(r-1)+i}
\quad\text{(by \ref{facts}(i))}\cr
&=\bc^{r-1}\bc_{i(r-1)+d-1}\bx_{i(r-1)+d}\bx_{i(r-1)+d-1}\ldots \bx_{i(r-1)+i}\cr
&=\bc^{r-1}\bc_{i(r-1)+d-1}\by_i\cr
\end{align*}
which, conjugating by $\by_i$, gives the stated value for $\bw_i$ .
\end{proof}

We note that in particular we have
$\by\bw\by\inv=\bw_1=\bc^{r-1}\bc_{d-1}\by_0=
\bsigma_{r,r+d-2}\bc^{r-1}\bx_d\ldots \bx_1$.

We prove now the analogous lemma for $\bw'$.
Let $\bc'_i=\bc_i\bsigma_i$.

\begin{lemma}\label{y' in B<-w'}
Let $\by'=\prod_{i=1}^{d-1}\by'_i$ where
$\by'_i=\prod_{j=d+1}^{i+1}\bx_{i(r-1)+j}$. Then
$\by'\bw'\by^{\prime-1}\in B^+$ and $\by'\in
\Hom_{\cD^+}(\by'\bw'\by^{\prime-1},\bw')$.
\end{lemma}
\begin{proof} We can assume $r\geq 2$ since for $r=1$ we have $\by'=1$.
We set $\bw'_d=\bw'$, and then by decreasing induction on $i$, 
we define  $\bw'_i=\by'_i\bw'_{i+1}\by_i^{\prime-1}$. We will show
that   $\by_i^{\prime-1}\bw'_i\in\BW$ by proving by decreasing induction on
$i$ that   $\bw'_i=\by'_i\bc_{rd}^{r-1}\bc'_{(r-1)i+d}$.
This equality has no meaning for $i=d$, as
$\bx_{dr+1}$ does not make sense. So we rewrite it as
$\bw'_i=\bc_{rd}^{r-1}\delta_{r-1}(\by'_i)\bc'_{(r-1)i+d}$, where
we define formally $\delta_i(\bx_j)=\bx_{j-i}$ which makes sense
when $j-i\le n$.
Now the formula for $i=d$ becomes
$\bw'=\bc_{rd}^r\bsigma_{rd-r+1,rd}$, which holds by \ref{facts} (v).
This is the starting point of the induction.
Let us assume the formula true for $i$ and let us prove it for $i-1$.
We have
\begin{align*}\bw'_i&=\bc^{r-1}_{rd}\delta_{r-1}(\by'_i)\bc'_{(r-1)i+d}\cr
&=\bc^{r-1}_{rd}\delta_{r-1}(\by'_i)\bc_{(r-1)i+d}\bsigma_{(r-1)i+d}\cr
&=\bc^{r-1}_{rd}\bc_{(r-1)i+d}\delta_{r-2}(\by'_i)\bsigma_{(r-1)i+d}\cr
\intertext{the last equality as $r\ge 2$ implies that the largest $j$ such
that $\bsigma_j$ occurs in $\delta_{r-1}(\by'_i)$ is at most $(r-1)i+d$,}
&=\bc^{r-1}_{rd}\bc_{(r-1)(i-1)+d}\by'_{i-1}\bsigma_{(r-1)i+d}\cr
&=\bc^{r-1}_{rd}\bc'_{(r-1)(i-1)+d}\by'_{i-1},
\end{align*}
where we get the last line by commuting
$\bsigma_{(r-1)i+d}$ and
$\bx_{(i-1)(r-1)+j}$ for $j<d$, and applying formula
$\bx_{a+1}\bx_a\bsigma_{a+r-1}=\bsigma_a\bx_{a+1}\bx_a$ to the two first terms
of $\by'_{i-1}$.
\end{proof}

From \ref{y' in B<-w'} we get
$\by'\bw'\by^{\prime-1}=
\bw'_1=\bc^{r-1}_{rd}\bc'_d\by'_0=
\bsigma_{r,r+d-1}\bsigma_{r+d-1}\bc_{rd}^{r-1}\bx_{d+1}\ldots \bx_1$.

Let us now prove

\begin{lemma}\label{yty-1}
\begin{itemize}
\item[(i)]
We have $\by\bt\by\inv=\bsigma_{r,r+d-2}$
\item[(ii)]
We have $\by'\bt'\by^{\prime-1}=\bsigma_{r,r+d-1}\bsigma_{r+d-1}$.
\end{itemize}
\end{lemma}
\begin{proof}
Let us prove (i).
We show that $\bt=\prod_{j=1}^{d-1}\ba_{jr,(j+1)r-1}$ where
$\ba_{k,l}=\bsigma_{k,l}\bsigma_{k,l-1}\inv$.
For this, replace each $\ba_{k,l}$ by its value to get
$\prod_{i=1}^{d-1}\bsigma_{ir,(i+1)r-1}
\prod_{i=1}^{d-1}\bsigma\inv_{ir,(i+1)r-2}$,
which is equal to
$\bsigma_{r,dr-1}\prod_{i=1}^{d-1}\bsigma\inv_{ir,(i+1)r-2}$,
in turn equal to $\bc_{r-1}\inv \bc\prod_{i=1}^{d-1}\bsigma\inv_{ir,(i+1)r-2}$.
We can put $\bc$ on the right of this product
if we replace $\bsigma_{k,l}$ by $\bsigma_{k+1,l+1}$ for all $k,l$.
We get $\prod_{i=0}^{d-1}\bsigma\inv_{ir+1,(i+1)r-1}\bc$, which is equal to
$\bt=(\prod_{i=1}^{r-1} \bs_i)\inv \bc$, as wanted, since
$\prod_{i=1}^{r-1}\bs_i= \prod_{j=0}^{d-1}\bx_{jr+1}$ where all factors
commute.

By decreasing induction on $i$, we find
$\prod_{j=i}^{d-1}\by_j=\prod_{j=d}^{i+1}\bsigma_{ir+j-i,jr-1}$,
so, setting $\bz_i=\bsigma_{i+r,(i+1)r-1}$, we have
$\by=\prod_{i=d-1}^1\bz_i$.
But $\bz_j$ conjugates
$\ba_{jr,(j+1)r-1}$ into $\ba_{j+r,(j+1)r-1}$ and
$\ba_{j+r-1,(j+1)r-1}$ into $\bsigma_{j+r-1}$. It commutes with $\bsigma_i$
for $i<j+r-1$ and with $\ba_{ir,(i+1)r-1}$ for $i>j$.
So by induction on $j$ we see that
$\prod_{i=j}^1 \bz_i$ conjugates $\bt$ into
$(\prod_{i=1}^j \bsigma_{i+r-1})\ba_{j+r,(j+1)r-1}
(\prod_{i=j+1}^{d-1}\ba_{ir,(i+1)r-1})$,
whence (i).

We now prove (ii). We claim that
$\by'=\bsigma_{d+r,dr}\by$: indeed
$\by'_i=\bx_{i(r-1)+d+1}\by_i$ and $\bx_{i(r-1)+d+1}$ commutes with
$\by_k$ for $k<i$; this gives the claim as
$\prod_{i=1}^{d-1}\bx_{i(r-1)+d+1}=\bsigma_{d+r,dr}$.

Let us now conjugate $\bt'=\bt\bsigma_{dr}^2$
by $\by'=\bsigma_{d+r,dr}\by$. By (i),
$\by'$ conjugates $\bt$ into $\bsigma_{r,r+d-2}$ as 
$\bsigma_{d+r,dr}$ commutes with $\bsigma_{r,r+d-2}$. 
On the other hand conjugation by $\by'$
has the same effect on $\bsigma_{dr}$ as conjugation by
$\bsigma_{d+r,dr}\bz_{d-1}=
\bsigma_{d+r,dr}\bsigma_{d+r-1,dr-1}$ as
$\bz_i$ for  $i<d-1$ commutes with $\bsigma_{dr}$. It remains to see that
$\bsigma_{d+r,dr}\bsigma_{d+r-1,dr-1}$ conjugates $\bsigma_{dr}$ into
$\bsigma_{d+r-1}$ \ie,
$\bsigma_{d+r,dr}\bsigma_{d+r-1,dr-1}\bsigma_{dr}=
\bsigma_{d+r-1}\bsigma_{d+r,dr}\bsigma_{d+r-1,dr-1}$; this can be written
$\bsigma_{d+r,dr}\bsigma_{d+r-1,dr}=
\bsigma_{d+r-1,dr}\bsigma_{d+r-1,dr-1}$, which is true.
\end{proof}

\begin{corollary}\label{t in End_D} 
We have $\bt\in\End_\cD(\bw)$ and $\bt'\in\End_\cD(\bw')$.
\end{corollary}
\begin{proof}
If $\bI=\{\bsigma_r,\ldots,\bsigma_{r+d-2}\}$,
then $\bsigma_{r,r+d-2}\in \End_{\cD^+_\bI}(\by\bw\by\inv)$,
as $\bsigma_{r,r+d-2}\preccurlyeq\by\bw\by\inv$.
By \ref{y in B<-w} and the remarks made above that lemma, we get the first
assertion.

Similarly, if $\bI'=\{\bsigma_r,\ldots,\bsigma_{r+d-1}\}$,
then $\bsigma_{r,r+d-1}\bsigma_{r+d-1}\in
\End_{\cD^+_{\bI'}}(\by'\bw'\by^{\prime -1})$
as $\bsigma_{r,r+d-1}\bsigma_{r+d-1}\prec\by\bw\by\inv$.
By \ref{y' in B<-w'}, we get the second assertion.
\end{proof}

\section{Regular elements in type $B$}\label{reg in B}
We now prove conjecture \ref{A} for roots of $\bpi$
when $W$ is of type $B_n$. We will see $W$ as the centralizer in a Coxeter
group $W'$ of type $A_{2n-1}$ of the longest element $w_0$. Let
$V'$ be the vector space which affords the reflection representation of $W'$. As $w_0$ is a
2-regular element, the eigenspace $V'_{-1}$ of $w_0$ 
affords the reflection representation of $W$,
and if we choose a base point in
$V^{\prime\text{reg}}_{-1}$ we get an embedding 
$\Pi_1(\Vr/W)\hookrightarrow\Pi_1(V^{\prime\text{reg}}/W')$ 
of the braid group $B$ of type
$B_n$ into the braid group $B'$ of type $A_{2n-1}$. If 
$\bS'=\{\bsigma'_1,\ldots\bsigma'_{2n-1}\}$
is the generating set of $B'$ then it is known (see \cite[4.4]{michel})
that $\bsigma_1=\bsigma'_n$ and
$\bsigma_{n+1-i}=\bsigma'_i\bsigma'_{2n-i}$ for $1\leq i <n$ are generators of
$B$ such that the relations are given by the Coxeter diagram
$\nnode \bsigma_1\dbar\nnode{\bsigma_2}\bar\nnode{\bsigma_3}
\cdots\nnode{\bsigma_n}\kern 4pt$.

Let $w$ be a $d$-regular element of $W$ for
some $d$, and let $\zeta=e^{2i\pi/d}$. We have
$\Vzr/C_W(w)\simeq (\Vr/W)^\zeta$ and
$\Vr/W=V^{\prime\text{reg}}_{-1}/C_{W'}(w_0)\simeq
(V^{\prime\text{reg}}/W')^{(-1)}$, so that
$\Vzr/C_W(w)\simeq(V^{\prime\text{reg}}/W')^{(-1,\zeta)}$. This is equal to
$(V^{\prime\text{reg}}/W')^\zeta$ if $d$ is even and to 
$(V^{\prime\text{reg}}/W')^{\zeta'}$ with $\zeta'=e^{i\pi/d}$ if $d$ is odd.
As $C_W(w)=C_{W'}(w,w_0)$, we see that if $w_0$ is a power of $w$,
which implies that $d$ is even, then $C_W(w)=C_{W'}(w)$.
In this case, as $\Vzr=V^{\prime\text{reg}}_\zeta$ we see that 
the map $\gamma$ from $B(w)$ to  $C_{B'}(\bw)$ is the composition of the map
which we still denote by $\gamma$ from $B(w)$ to
$C_B(\bw)$ and of the embedding $B\to B'$.

We make a specific choice of a regular element.
Let $\bc=\bsigma_1\ldots\bsigma_n$. It is a $2n$-th root of $\bpi$
(\cf\ \ref{h-th root of pi}).
Let $r$ and $d$ be two integers such that $rd=2n$ and let $\bw=\bc^r$;
it is a $d$-th root of $\bpi$ and its image $w\in W$ is a good regular
element.

If $d$ is odd, we have
$C_B(\bw)=C_B(\bc^r)=C_B(\bc^{pgcd(n,r)})=C_B(\bc^{r/2})$ and
$C_W(w)=C_W(c^{r/2})$ since $w_0=c^n$ is central.
So we are reduced to study $C_B(\bw)$ when $d$ is even and $C_B(\bw^2)$ when
$d$ is even and $d/2$ odd (see \ref{c^i dans  B_n} (ii), below).

When $d$ is even we have seen above that $C_W(w)=C_{W'}(w)$, and
this group is a complex reflection group of type $G(d,1,r)$
(see also \cite[A.1.2]{Sydney}).

We have $\bc=\bsigma'_n\bsigma'_{n-1}\bsigma'_{n+1}
\ldots\bsigma'_1\bsigma'_{2n-1}$. In order to apply the results of
\ref{reg in A}, we use a conjugation by $\bv\inv$ where
$\bv$ is the canonical lift of the longest element of the parabolic
subgroup of $W'$
generated by $\sigma'_1,\ldots,\sigma'_{n-1}$: indeed we have
$\bc=\bv\inv\bsigma'_1\ldots\bsigma'_{2n-1}\bv$.

This proves that the generators of $C_B(\bw)$ are
the elements $\bs_i=\bv\inv(\prod_{j=0}^{d-1}
\bsigma'_{i-1+rj})\bv=\prod_{k=0}^{d/2-1}\bsigma_{i+kr}$
for $i=2\ldots r$ and the element
$\bt$  such that $\bt\prod_{i=2}^r\bs_i=\bc$; note that we have not chosen for
the element $\bt$ the conjugate by $\bv\inv$ of the corresponding element of
\ref{reg in A} but we have applied a further conjugation by $\bc$ in order to
simplify the computations.
We have $\bt=(\hat\bsigma_2\bsigma_3\ldots\bsigma_{r+1}\hat\bsigma_{r+2}
\ldots\bsigma_{(d/2-1)r+1})\inv\bsigma_1\bsigma_2\ldots\bsigma_{(d/2-1)r+1}$,
where $\hat\bsigma_i$ means deleting $\bsigma_i$ from the product:
we have deleted all $\bsigma_i$ such that $i\equiv 2\pmod r$.

We first prove a lemma analogous to \ref{facts}.

\begin{lemma}\label{factsBn}
\begin{itemize}
\item[(i)] We have
$\bc\bsigma_i=\bsigma_{i+1}\bc$ for $2\le i<n-1$.
\item[(ii)]
We have $\lexp{\bc^2}\bsigma_n=\bsigma_2$.
\item[(iii)] For $\bx\in B^+$ and $2\le i<n$ we have $\bsigma_i\prec\bx
\Leftrightarrow \bsigma_{i+1}\prec\bc\bx$.
\item[(iv)] We have $\{i\mid\bsigma_i\prec\bc^j\}=\{1,\ldots,j\}$ for $j\le n$.
\end{itemize}
\end{lemma}
\begin{proof} Statements (i)  and (iii) have the same proof as
the corresponding statements in \ref{facts}.

Let us prove (ii). We have
\begin{multline*}
\bc^2=\bsigma_1\bsigma_2(\bsigma_3\ldots\bsigma_n)\bc=
\bsigma_1\bsigma_2\bc(\bsigma_2\ldots\bsigma_{n-1})=
\bsigma_1\bsigma_2\bsigma_1\bsigma_2(\bsigma_3\ldots\bsigma_n)
(\bsigma_2\ldots\bsigma_{n-1})
=\\\bsigma_2\bsigma_1\bsigma_2\bsigma_1(\bsigma_3\ldots\bsigma_n)
(\bsigma_2\ldots\bsigma_{n-1})
=\bsigma_2\bsigma_1\bsigma_2(\bsigma_3\ldots\bsigma_n)
\bsigma_1(\bsigma_2\ldots\bsigma_{n-1})
=\bsigma_2\bc\bsigma_1\ldots\bsigma_{n-1},
\end{multline*}

whence $\bc^2\bsigma_n=\bsigma_2\bc^2$.

The proof of (iv) is similar to that of \ref{facts} (iv): it uses similarly the
fact that $\bsigma_i\prec\bc$ if and only if $i=1$, but it also needs the fact
that
$\bsigma_2\prec\bc^2$, which we have seen in the proof of (ii).
\end{proof}

\begin{lemma}\label{ts2ts2}
The group $C_B(\bw)$ has a presentation with generators
$\bt$, $\bs_2$, \dots, $\bs_r$, the relations  being the braid relations given
by the diagram
$\nnode \bt\dbar\nnode{\bs_2}\bar\nnode{\bs_3}\cdots\nnode{\bs_r}\kern 4pt$.
\end{lemma}
\begin{proof}
We already know that $C_B(\bw)$ has a presentation with generators
$\bs_2$, \dots, $\bs_r$, $\bt'$ and relations the braid
relations given by
$\nnode{\bs_2}\bar\nnode{\bs_3}\cdots\nnode{\bs_r}\dbar\nnode{\bt'}\kern 4pt$,
where $\bt'$ is the element conjugate by
$\bv\inv$ of the element $\bt$ of section \ref{reg in A}. 
We have
$\bt=\bc\bt'\bc\inv=(\bs_2\dots\bs_r)\bt'(\bs_2\dots\bs_r)\inv$. The
commutation relation $\bt'$ with $\bs_i$ for
$i<r$ is equivalent by \ref {factsBn} (i) to the commutation relation
of $\bt$ with $\bs_i$ for $i>2$. It remains to see that the braid relation
between $\bt'$ and $\bs_r$ is equivalent to the braid relation between
$\bt$ and $\bs_2$. This is proved by decreasing induction on $i$
using the following fact that is the result of a simple computation:
if $\bs_{i-1}$, $\bs_i$ and $\bu'$ are elements of a group and if
$\bu=\bs_i\bu'\bs_i\inv$ then the
braid relations given by
$\nnode{\bs_{i-1}}\bar\nnode{\bs_i}\dbar\nnode{\bu'}\kern 4pt$
imply $\bs_{i-1}\bu\bs_{i-1}\bu=\bu\bs_{i-1}\bu\bs_{i-1}$.
Conjugating by $\bs_2$ the relation we get at the end of the induction
the braid relation which we want. The converse is similar.
\end{proof}

\begin{theorem}\label{c^i dans  B_n}
Assume $d$ even; then
\begin{itemize}
\item[(i)]We have
$\bs_i\in\End_{\cD^+}(\bw)$ ($i=1\ldots,r-1$)
and  $\bt\in\End_\cD(\bw)$,   so that    $\End_\cD(\bw)=C_B(\bw)\simeq
B(d,1,r)$.

\item[(ii)]
If $d/2$ odd, we have
$\bs_i\in\End_{\cD^+}(\bw^2)$ ($i=1\ldots,r-1$)
and  $\bt\in\End_\cD(\bw^2)$,   so that    $\End_\cD(\bw^2)=C_B(\bw^2)\simeq
B(d,1,r)$.
\end{itemize}
\end{theorem}
\begin{proof}
In the following lemma, the statement about $\bw^2$ assumes that $d/2$ is odd.
The proof follows the same lines as that of
\ref{s_i ok}, using \ref{factsBn} instead of \ref{facts}.
\begin{lemma}\label{s_i okB} 
Assume $2\le i\le r$ and let
$\bI_i=\{\bsigma_i,\bsigma_{i+r},\ldots,\bsigma_{i+(d/2-1)r}\}$. Then
\begin{itemize}
\item[(i)] $\bs_i\in \End_{\cD^+_{I_i}}(\bw)$ and
$\bs_i\in \End_{\cD^+_{I_i}}(\bw^2)$. In particular
$\bs_i\in\End_{\cD^+}(\bw)$ and
$\bs_i\in\End_{\cD^+}(\bw^2)$.
\item[(ii)] The conjugation
by $\bw$ (resp. $\bw^2$) stabilizes $\bI_i$ and induces the cyclic permutation
$\bsigma_{i+jr}\mapsto\bsigma_{i+(j+1)r\pmod n}$
(resp. $\bsigma_{i+jr}\mapsto\bsigma_{i+(j+2)r\pmod n}$).
\end{itemize}
\end{lemma}

\begin{lemma}\label{y in B<-wBn}
For $i\le j$ we set $\bsigma_{i,j}=\bsigma_i\bsigma_{i+1}\ldots\bsigma_j$,
and we set $\bx_i=\bsigma_{i+1,i+r-1}$.
Let $\by=\prod_{i=1}^{d/2-1}\by_i$ where
$\by_i=\prod_{k=1}^{d/2-i}\bx_{(i-1)(r-1)+d/2-k+1}$. Then
$\by\bw\by\inv\in B^+$ and $\by\in \Hom_{\cD^+}(\by\bw\by\inv,\bw)$.
\end{lemma}
\begin{proof}
We set $\bc_i=\bsigma_{1,i}$.
Let $\bw_{d/2}=\bw$, and by decreasing induction on $i$ define 
$\by_i\inv\bw_i\by_i=\bw_{i+1}$.
We claim that $\by_i\inv\bw_i\in\BW$, which implies the result: in fact
we prove by induction that
$\bw_i=\by_i\bc_{i(r-1)+d/2}\bc^{r-1}$.
This equality is clearly true for $i=d/2$ (with $\by_{d/2}=1$).
Let us assume it to be true for $i+1$ and let us prove it for $i$.
We have
\begin{align*}\bw_{i+1}&=\by_{i+1}\bc_{(i+1)(r-1)+d/2}\bc^{r-1}\cr
&=\prod_{k=1}^{d/2-i-1}\bx_{i(r-1)+d/2-k+1}\bc_{(i+1)(r-1)+d/2}\bc^{r-1}\cr
&=\bc_{(i+1)(r-1)+d/2}\prod_{k=1}^{d/2-i-1}\bx_{i(r-1)+d/2-k}\bc^{r-1}
\quad\text{(by \ref{factsBn}(i))}\cr
&=\bc_{i(r-1)+d/2}\prod_{k=0}^{d/2-i-1}\bx_{i(r-1)+d/2-k}\bc^{r-1}\cr
&=\bc_{i(r-1)+d/2}\bc^{r-1}\prod_{k=0}^{d/2-i-1}\bx_{(i-1)(r-1)+d/2-k}
\quad\text{(by \ref{factsBn}(i))}\cr
&=\bc_{i(r-1)+d/2}\bc^{r-1}\by_i\cr
\end{align*}
which, conjugating by $\by_i$ gives the equality for $\bw_i$.
\end{proof}
We note that the above proof shows that
$$\by\bw\by\inv=\by_1\bc_{r-1+d/2}\bc^{r-1}=\prod_{i=d/2}^2\bx_i
\prod_{i=1}^{r-1+d/2}\bsigma_i\bc^{r-1}=\bsigma_{1,r-1+d/2}
\prod_{i=d/2-1}^1\bx_i\bc^{r-1}$$
the last equality by \ref{factsBn} (i), so
$\by\bw\by\inv=\bsigma_{1,d/2}\prod_{i=d/2}^1\bx_i\bc^{r-1}$. 

\begin{lemma}\label{yty-1Bn}
We have $\by\bt\by\inv=\bsigma_{1,d/2}$.
\end{lemma}
\begin{proof}
For $i=1,\ldots,d/2-1$, let
$\bt_i= (\hat\bsigma_2\bsigma_3\ldots\bsigma_{r+1}\hat\bsigma_{r+2}
\ldots\bsigma_{ir+1})\inv\bsigma_{1,i(r-1)+d/2}$.
We have $\bt=\bt_{d/2-1}$ and $\bt_0=\bsigma_{1,d/2}$.
We prove by induction  that $\by_i$ (\cf\ \ref{y in B<-wBn})
conjugates $\bt_i$ into $\bt_{i-1}$, which proves the lemma.
Keeping the notation of \ref{y in B<-wBn}, we have
$\bt_i=(\bx_2\bx_{r+2}\ldots\bx_{(i-1)r+2})\inv\bsigma_{1,i(r-1)+d/2}$.
By definition
$\by_i=\prod_{k=1}^{d/2-i}\bx_{(i-1)(r-1)+d/2-k+1}$. This product commutes
with $\bx_2\bx_{r+2}\ldots\bx_{(i-2)r+2}$. As the factor indexed by
$k=d/2-i$ in $\by_i$ is equal to $\bx_{(i-1)r+2}$, we get
$$\by_i\bt_i\by_i\inv=(\bx_2\bx_{r+2}\ldots\bx_{(i-2)r+2})\inv
\prod_{k=1}^{d/2-i-1}\bx_{(i-1)(r-1)+d/2-k+1}
\bsigma_{1,i(r-1)+d/2}\by_i\inv.$$
We use the fact that conjugation by
$\bsigma_{1,i(r-1)+d/2}$ of the factor
$\bx_{(i-1)(r-1)+d/2-k+1}$ in $\by_i\inv$ for $k>1$
changes $k$ into $k-1$:  this allows to simplify the product and we get
\begin{align*}
\by_i\bt_i\by_i\inv&=(\bx_2\bx_{r+2}\ldots\bx_{(i-2)r+2})\inv
\bsigma_{1,i(r-1)+d/2}\bx_{(i-1)(r-1)+d/2}\inv\cr
&= (\bx_2\bx_{r+2}\ldots\bx_{(i-2)r+2})\inv
\bsigma_{1,(i-1)(r-1)+d/2}=\bt_{i-1}.
\end{align*}
\end{proof}

Let $\bI=\{\bsigma_1,\ldots,\bsigma_{d/2}\}$;  we have
$\bsigma_{i,d/2}\in\End_{\cD^+_\bI}(\by\bw\by\inv)$, by the remark
following the proof of \ref{y in B<-wBn}.
This, together with lemmas \ref{y in B<-wBn} and \ref{yty-1Bn}
proves the statements about $\bt$ in the theorem.
\end{proof}

\vfill\eject
\newcommand{\opp}{{\mathrm opp}}

\section{The elements $\bpi$ and $\bw_0$}\label{secpi}

We  consider  here the  order  1  $F$-root  of  $\bpi$, given  which  is
$\by=\bpi$, and the  order 2 $F$-root $\by=\bw_0$. For  $\bw_0$, we will
just show how conjectures \ref{A}  to \ref{C} follow from known results.
For $\bpi$,  we will in addition  prove conjecture \ref{D} in  a certain
number of cases, including split type $A$ in general.

We recall  that (\cf\ \cite[proposition 2.1.6]{DMR})  the group $C_W(F)$
(resp.  $C_W(w_0F)$) is  a  Coxeter group  with  Coxeter generators  the
elements $w_0^I$  for $I$ an element  of the set of  orbits $S/F$ (resp.
$I\in  S/w_0F$).  The corresponding  braid  groups  $C_B(\bpi F)$  (resp
$C_B(\bw_0 F)$) have as generators the corresponding elements $\bw_0^I$.

Since the  generators $\bw_0^I$ divide  $\bpi$ (resp $\bw_0$)  in $B^+$,
conjecture \ref{A} is trivial. For $\bpi$, conjecture \ref{conj} is also
trivial. For $\bw_0$ it results from the remark below \ref{conj} and the
fact that $\bw_0$  is the only ``good'' square $F$-root  of $\bpi$ since
it is the only element of $\bW$ of its length.

Since  $\bpi  F$ (resp.  $\bw_0  F$)  acts  as a  diagram  automorphism,
conjecture \ref{B} is \cite[corollary 4.4]{michel}.

Conjecture \ref{C} holds for the cases $\by=\bpi$ and $\by=\bw_0$ by the
following results:

\begin{proposition}\label{pi} \cite[5.3.4]{DMR}
The map $\bt\mapsto  D_\bt$ from $C_B(F)$ to  the $\GF$-endomorphisms of
$H^*_c(\bX(\bpi))$ factors through the  specialization $x\mapsto q$ of a
1-cyclotomic  Hecke algebra  for  $C_W(F)$ which  is the  specialization
$ u_{\bw_0^I,0}\mapsto  x^{l(w_0^I)}$, $u_{\bw_0^I,1}\mapsto  -1$   of  the
generic Hecke algebra of $C_W(F)$.
\end{proposition}

\begin{proposition}\label{w0} \cite[5.4.1]{DMR}
The  map   $\bt\mapsto  D_\bt$  from  $C_B(w_0   F)$  to
the   $\GF$-endomor\-phisms  of   $H^*_c(\bX(\bw_0))$  factors   through
the  specialization  $x\mapsto  q$   of  a  2-cyclotomic  Hecke  algebra
for  $C_W(w_0  F)$  which is  the  specialization 
$u_{\bw_0^I,\varepsilon}\mapsto x^{l(w_0^I)}$, $u_{\bw_0^I,1-\varepsilon}
\mapsto (-1)^{1+l(\bw_0^I)}$, where $\varepsilon=0$ if $l(w_0^I)$ is even and
$1$ otherwise, of the generic Hecke algebra of $C_W(w_0 F)$.
\end{proposition}

We will now consider conjecture \ref{D} for $\by=\bpi$. Let $\CH$ be the
cyclotomic algebra of  \ref{pi}. Since the characters of  $\CH$ are only
defined over  $\Qlbar[x^{1/2}]$ (for $W$ irreducible,  this happens only
for the characters of degree 512 of $W(E_7)$ and those of degree 4096 of
$W(E_8)$),  we  need to  take  the  integer  $a$ defined  above  \ref{C}
equal  to  $2$  and  thus consider  the  specialization  $f:x^{1/2}\mapsto
q^{1/2}$  of  the  algebra  with  parameters
$u_{\bw_0^I,0}\mapsto ((-x^{1/2})^{2l(w_0^I)})$
and $u_{\bw_0^I,1}\mapsto  -1$.  In  the terms  of
\cite[5.3]{DMR} this  corresponds to the  specialization $x^{1/2}\mapsto
-q^{1/2}$   of   $\CH_x(W,F)$.    We   recall   from   \cite[5.3.2]{DMR}
that   if   $f'$   is  the   specialization   $x^{1/2}\mapsto  -q^{1/2}$
(which corresponds to the specialization   $x^{1/2}\mapsto q^{1/2}$
of $\CH_x(W,F)$),
and   if   we   fix   an   $F$-stable   Borel   subgroup   $\bB$,   then
$\CH\otimes_{f'}\Qlbar\simeq\End_{\Qlbar\GF}(\Ind_\BF^\GF\Id)$.

Let   $\sigma$   be   the    semi-linear   automorphism, coming   from
$\Gal(\Qlbar(x^{1/2})/\Qlbar(x)$, given by  $x^{1/2}\mapsto -x^{1/2}$ of
$\CH$; thus $f\circ\sigma=f'$. 
Let   $\CH_q=\CH\otimes_{f'}\Qlbar$,  let
$\chi\mapsto\chi_q$  be the  bijection between  characters of  $\CH$ and
$\CH_q$ obtained via $f'$, and
let $\chi_q\mapsto\rho_\chi$ be the bijection
between characters of
$\CH_q$  and  characters of  $\GF$  occurring  in $\Ind_\BF^\GF\Id$  coming
from \cite[5.3.2]{DMR}. 

Let   us   recall   that   the   representation   $\Ind_\BF^\GF\Id$   of
$\CH\otimes_{f'}\Qlbar$ is special (\cf\ \ref{D} (ii));
this follows from the
fact that  the image of  any non-trivial  $\bw\in\bW$ has zero  trace in
this representation,  which characterizes  the canonical trace  form for
Hecke algebras of Coxeter groups.

It follows that conjecture \ref{D} is implied by the

\begin{conjecture}\label{conjecture Xpi}

$\sum_i(-1)^i H^i_c(\bX(\bpi))=
\sum_{\chi\in\Irr(W^F)}\rho_\chi\otimes \sigma(\chi)_q$.
\end{conjecture}

\smallskip
In the remaining part of this section we will prove the following theorem:
\begin{theorem}\label{Endomorphismes de Xpi}
Conjecture \ref{conjecture  Xpi} holds  if the characteristic  is almost
good for $\bG$ and if $(W,F)$  is irreducible of type untwisted
$A_n$, $B_2$, $B_3$, $B_4$, $D_4$, $D_5$, $D_6$, $D_7$, $G_2$ or $E_6$.
\end{theorem}
Recall that  the characteristic is almost  good for $\bG$ if  it is good
for each simple component of exceptional type of $\bG$.

\begin{proof}
We have to prove that the virtual character of $\CH_q$
appearing in the $\rho_\chi$-isotypic component
of $\sum_i(-1)^iH^i_c(\bX(\bpi),\Qlbar)$ is equal to
$\sigma(\chi)_q$. 
This is equivalent to proving that for any $\bx\in  C_{B^+}(F)$ and any
$\chi\in\Irr(W^F)$, we have:

\begin{equation}
\scal{g\mapsto\TrH{gD_\bx}{\bX(\bpi)}}
{\rho_\chi}\GF=\sigma(\chi)_q(T_\bx),\tag{1}\end{equation}  

where $T_\bx$ denotes the image of $\bx$ in $\CH_q$.
We  will   prove  this  equality   for  sufficiently  many   elements  of
$C_{B^+}(F)$ to deduce it for all elements for groups in the list
of \ref{Endomorphismes de Xpi}. In the next two lemmas we need not assume
$\bG$ split or irreducible.

\begin{lemma}\label{(1) pour x parabolique}
If conjecture  \ref{conjecture Xpi} holds  for any reductive  group with
semi-simple rank less than that of $\bG$ then (1) holds for any $\bx\in
C_{B^+_\bI}(F)$ for any $F$-stable proper subset $\bI$ of $S$.
\end{lemma}
\begin{proof}
If $\bL_I$ is the standard Levi subgroup of $\bG$ corresponding to
$I$, then by
\cite[th\'eor\`eme 5.2.10]{DMR}, for $\bx\in C_{B^+_\bI}(F)$ we have:
$$\displaylines{\scal{g\mapsto\TrH{gD_\bx}{\bX(\bpi)}}
{\rho_\chi}\GF=\hfill\cr \hfill\scal{l\mapsto\TrH{lD_\bx}
{\bX_{\bL_I}(\bpi_\bI)}} {\lexp*R^\GF_{\bL_I^F}\rho_\chi}{\bL_I^F},\cr}$$
which by assumption is equal to
$$\sum_{\varphi\in\Irr(W^F_I)}\sigma(\varphi)_q(T_\bx) \scal{\rho_\varphi}
{\lexp*R^\GF_{\bL_I^F}\rho_\chi}{\bL_I^F}.$$
As $\scal{\rho_\varphi} {\lexp*R^\GF_{\bL_I^F}\rho_\chi}{\bL_I^F}=
\scal\varphi{\Res^{W^F}_{W_I^F}\chi}{W_I^F}$, \cf\ \cite[Theorem 70.24]
{Curtis-Reiner}, we get (1) for $\bx$.
\end{proof}

\begin{lemma}\label{(1) pour pi^n}
Equality (1) holds if $\bx=\bpi^n$, with $n$ multiple of $\delta$.
\end{lemma}
\begin{proof}
As $\bpi^n$ acts by $F^n$ on $\bX_\bpi$, we have
$\TrH{gD_{\bpi^n}}{\bX(\bpi)}=|\bX(\bpi)^{gF^n}|$ by the Lefschetz trace formula.
We shall use the same methods and notation as in \cite[\S 2.B et \S 6.D]{Sydney}.
Proposition \cite[3.3.7]{DMR} shows that
$$|\bX(\bpi)^{gF^n}|=\sum_{\rho\in\Irr(\GF)}\rho(g)\sum_{\chi\in\Irr(W)^F}
\langle \rho,R_{\tilde{\chi}}\rangle\tilde{\chi}_q(T_\pi^n F).$$
We have $\tilde{\chi}_q(T_\pi^n F)=q^{n(2N-a_{\chi}-A_{\chi})}
\tilde{\chi}(F)$ \cite[proposition 6.11]{Sydney} whence we get as in the proof of
\cite[Proposition 2.5]{Sydney}
$$|\bX(\bpi)^{gF^n}|=\sum_{\rho\in\Irr(\GF)}\rho(g)
\langle \rho,\Ind_{\bB^F}^{\bG^F}\Id\rangle q^{n(2N-a_{\rho}-A_{\rho})}=
\sum_{\chi\in\Irr(W^F)}
\chi(1)q^{n(2N-a_{\chi}-A_{\chi})}\rho_\chi(g).$$
We now use $\chi_q(T_\bpi^n)=\sigma(\chi)_q(T_\bpi^n)=\chi(1)
q^{n(2N-a_\chi-A_\chi)}$
\cite[corollaire 4.21]{Sydney}, which gives the result.
\end{proof}

\begin{lemma}\label{(1) pour racines de pi}
If the characteristic is almost good for the split irreducible group $\bG$,
equality (1) holds when $\bx$ is a root of $\bpi$.
\end{lemma}
\begin{proof} \ref{(1) pour pi^n} shows the result for $\bx=\bpi$; we thus
assume that $\bx$ is a $d$-th root of $\bpi$ with $d\geq 2$.
By  \cite[5.2.2 (i)]{DMR}
the endomorphism $D_\bw$ of $\bX(\bpi)$ satisfies the trace formula so that
$$\TrH{gD_\bw}{\bX(\bpi)}=|\bX(\bpi)^{gD_\bw}|.$$ Moreover, by
\cite[5.2.2 (ii)]{DMR}, we have
$(g\mapsto|\bX(\bpi)^{gD_\bw}|)=\Sh^d(g\mapsto\Trace(gT_\bw\mid\Ind_\BF^\GF\Id))$.
So we have to prove
$$\Sh^d(g\mapsto\Trace(gT_\bw\mid\Ind_\BF^\GF\Id))=
\sum_{\chi_q\in\Irr(\CH_q(W,F))}\sigma(\chi)_q(T_\bw)\rho_\chi,$$
which is equivalent to
\begin{equation}\sum_{\chi_q\in\Irr(\CH_q(W,F))}\chi_q(T_\bw)\Sh^d\rho_\chi=
\sum_{\chi_q\in\Irr(\CH_q(W,F))}\sigma(\chi)_q(T_\bw)\rho_\chi.\tag{$1'$}
\end{equation}

To prove this, we may replace $\bw$ by a conjugate in $B$ so we may
assume that $\bw$ is a ``good'' root, in particular that $\bw\in\bW$.
As usual we set $w=\beta(\bw)$.

We have
$\chi_q(T_\bw)=\chi(w)q^{\frac{2N-a_{\rho_\chi}-A_{\rho_\chi}}d}$:
one gets this by applying \cite[6.15(2)]{spetses} to $\CH$;
it is a principal algebra (see \loccit\ 6.3), with
$ \theta_0(\bw_0^\bI)=x^{l(w_0^I)}$; we have $D_0=2N$ and if we take
$P(q)=|(\bG/\bB)^F|$ the degree $\hbox{Deg}^{(P)}_\chi$ identifies with the
generic degree of $\rho_\chi$.

As  $\CH$  is split over $\bbZ[x^{1/2},x^{-1/2}]$,
there exists a sign  $\varepsilon_{d,\chi}$  depending only on
$(a_{\rho_\chi}+A_{\rho_\chi})/d$ such that
$\sigma(\chi)_q(T_\bw)=\varepsilon_{d,\chi} \chi_q(T_\bw)$.
This sign is equal to $-1$ if and only if $(a_{\rho_\chi}+A_{\rho_\chi})/d\in
\bbZ+1/2$ and $\chi(w)\ne  0$.
Equation ($1'$) becomes then

\begin{equation}
\sum_{\chi_q\in\Irr(\CH_q(W,F))}q^{\frac{2N-a_{\rho_\chi}-A_{\rho_\chi}}d}
\chi(w)\Sh^d\rho_\chi=\\
\sum_{\chi_q\in\Irr(\CH_q(W,F))}\varepsilon_{d,\chi}
q^{\frac{2N-a_{\rho_\chi}-A_{\rho_\chi}} d}\chi(w)\rho_\chi.\tag{$1''$}
\end{equation}

For   computing  $\Sh^d$,   we  shall   use  Shoji's   results  on   the
identification of  character sheaves with  almost characters. Here  we need
the assumption that  the characteristic is almost good.  We recall these
results: unipotent characters of $\GF$ have been divided by Lusztig into
families.  Unipotent  character  sheaves  have also  been  divided  into
families  which  are  in  one-to-one correspondence  with  the  families
of  unipotent  characters.  In  \cf\  \cite[3.2  et  4.1]{Shoji2}  Shoji
proves  that the  transition  matrix from  the  unipotent characters  to
the  characteristic  functions of  the  unipotent  character sheaves  is
block diagonal  according to the  families, and in  \cite[3.3]{Shoji} he
proves that  the characteristic functions  of the character  sheaves are
eigenvectors of  $\Sh$. From this we  see that ($1''$) is  equivalent to
the  set  of its  projections  on  each  family. Moreover  $a_\rho$  and
$A_\rho$  are constant  when  $\rho$  runs over  a  family of  unipotent
characters. So ($1''$) is equivalent to the set of equations

\begin{equation}\sum_{\rho_\chi\in\CF}\chi(w)\Sh^d\rho_\chi=
\varepsilon_{d,\CF}\sum_{\rho_\chi\in\CF}
\chi(w)\rho_\chi,\tag{$1'''$}\end{equation}
where $\CF$ runs over the families.
We have written $\varepsilon_{d,\CF}$ instead of $\varepsilon_{d,\chi}$
because this sign depends only on the family of $\rho_\chi$.

We can  also assume that  $\bG$ is  adjoint as the  unipotent characters
factorize through  the adjoint group  and $\Sh$ is compatible  with this
factorization. 

Lusztig  defined  in  \cite[4.24.1]{Lubook} almost  characters  $R_\rho$
indexed  by  unipotent  characters.  If  $R_w$  is  the Deligne-Lusztig
character   given by the virtual representation
$\sum_{i\ge    0}(-1)^i    H^i_c(\bX(\bw))$,
we     have     $R_w=\sum_{\chi\in\Irr(W)}\chi(w)R_{\rho_\chi}$;     for
any      unipotent      character      $\rho$     we      also      have
$\scal\rho{R_{\rho_\chi}}\GF=\Delta_\rho\scal{R_\rho}{\rho_\chi}\GF$ for
a  sign  $\Delta_\rho$  defined   in  \cite{Lubook}.  Almost  characters
being    an   orthonormal    basis   of    the   space    of   unipotent
class   functions,   we   get   $\sum_{\rho_\chi\in\CF}\chi(w)\rho_\chi=
\sum_{\rho\in\CF}\scal{R_w}\rho\GF \Delta_\rho R_\rho$.

In  \cite[23.1]{LuCS}  Lusztig  has  defined  a  bijection  $\rho\mapsto
A_\rho$ from  the set of  unipotent characters  to the set  of unipotent
character  sheaves,   compatible  with  the  partition   into  families.
Shoji, in  (\cite{Shoji} and \cite{Shoji2})  proved   that  the   almost
character  $R_\rho$ is  a  multiple of  the  characteristic function  of
$\chi_{A_\rho}$  relative  to  the  Frobenius endomorphism  $F$  of  the
character  sheaf  $A_\rho$  and  that (\cf\  \cite[3.6  et  3.8]{Shoji})
$\Sh(\chi_{A_\rho})=\lambda_\rho \chi_{A_\rho}$  where $\lambda_\rho$ is
as in \cite[3.3.4]{DMR}.

Using this, we see that $(1''')$ is equivalent to:

\begin{equation}\text{if } 
\scal{R_w}\rho\GF\ne  0\text{  then  }\lambda_\rho^d=\varepsilon_{d,\CF}
\tag{$1''''$}\end{equation}

This  would  be  a consequence of  conjecture \cite[5.13]{Sydney}.  
We  prove it  by a case by case analysis.

If $\bG$ is classical, we have  always $\varepsilon_{d,\CF}=1$ and
$\lambda_\rho=\pm 1$, so $(1'''')$ holds
if $d$ is even. Assume $d$  odd; one checks that in a Coxeter group
of type $A_n$, $B_n$ or $D_n$, any odd order element lies in a parabolic
subgroup of  type $A$.  Let us  denote by  $\bL$ the  corresponding Levi
subgroup of $\bG$,  which  is  an  $F$-stable  Levi  subgroup  of  an 
$F$-stable
parabolic  subgroup. We  have  $R_w=R_\bL^\bG(R_w^\bL)$ where  $R_w^\bL$
is  the  Deligne-Lusztig  character  of  $\LF$  associated  to  $w$.  As
$\lambda_\rho$ is constant in a Harish-Chandra  series and is equal to 1
for a group of type $A$, we get the result in this case.

If $\bG$  is of exceptional  type
we  can  check  the  result, using  the  explicit  description  of
the  coefficients $\scal{R_w}\rho{\bG^{F'}}$  and  of $\lambda_\rho$  in
\cite{Lubook}.
The most complicated case to check is when for some
$d$ we have  $\varepsilon_{d,\CF}=-1$. In  type $E_7$
there is  exactly one such  family; it contains 4  unipotent characters.
Two of them  are some $\rho_\chi$ for a $\chi$ such that  $a_\chi+A_\chi=63$. 
In type
$E_8$ there are two such families,  each with 4 unipotent characters. In
each of  these families there are  two $\rho_\chi$
with respectively $a_\chi+A_\chi=105$ and $a_\chi+A_\chi=135$. So in all
cases we  have $\varepsilon_{d,\CF}=-1$ if  and only if  $d\equiv 2\pmod
4$. In each case for the two other unipotent characters of the family one
has $\lambda_\rho=\pm i$. One checks that if $\rho\in\CF$
and  $\scal\rho{R_w}\GF\ne 0$  then if $d\not\equiv 2\pmod 4$
one has $\lambda_\rho^d\ne1$ whence the result in this case; and if
$d\equiv  2\pmod 4$  we have  $\lambda_\rho=\pm i$, thus
$\lambda_\rho^d=-1$ and we also get the result in that case.
\end{proof}

Let  $\Phi$ be  the class  function on  $\CH_q$ with values in $\CR(\GF)$
given by
$$\Phi(T_\bx)=(g\mapsto\TrH{gD_\bx}{\bX(\bpi)})
-\sum_\chi\sigma(\chi)_q(T_\bx)\rho_\chi.$$  To  prove  theorem
\ref{Endomorphismes de Xpi} we have to prove that $\Phi=0$. By
\ref{(1) pour  x parabolique}, \ref{(1) pour pi^n}
and \ref{(1) pour  racines de pi} respectively 
we know that
\begin{itemize}
\item[(a)] $\Phi(T_\bx)=0$
for $\bx\in  B_\bI$ for any  proper subset $I$ of $S$.
\item[(b)] $\Phi(T_\bpi^n)=0$ for $n>0$.
\item[(c)] We have $\Phi(T_\bx)=0$
if  $\bx$ is a root of $\bpi$ and the characteristic is almost good.
\end{itemize}

We shall prove that in any of the cases considered in \ref{Endomorphismes
de Xpi}  a class function on $\CH_q$ which satisfies these three properties is 
zero. Such a class function, can be written
$\sum_\chi\lambda_\chi\chi_q$. We show that the three above properties imply
$\lambda_\chi=0$ for all $\chi$. Let us translate each of these properties into
a property of $(\lambda_\chi)_\chi$.

\begin{lemma}\label{(a)} Property (a) means that $(\lambda_\chi)_\chi$ is linearly spanned by
vectors $(\chi(w))_\chi$ with $w\in W$
{\it cuspidal} (\ie, the conjugacy class of $w$ has no representative in a 
proper parabolic subgroup of $W$).
\end{lemma}
\begin{proof}
Consider  the   scalar  product   on  $\CR(\CH_q)$  such   that  the
$\chi_q$  form an  orthonormal  basis (which  corresponds  to the  usual
scalar  product   on  the   vectors  $(\lambda_\chi)_\chi$);   then  the
$\chi_{c,q}=\sum_{\chi\in\Irr(W^F)}\chi(c)\chi$ are  pairwise orthogonal
when $c$ runs over a set  of representatives of the conjugacy classes in
$W$. The statement to prove is  that a class function satisfies property
(a) if  and only if  it is orthogonal to  the $\chi_{c,q}$ with  $c$ non
cuspidal.

With our  choice of  scalar product,  restriction and  induction satisfy
Frobenius  reciprocity, as  the scalar  product is  compatible with  the
specialization  to  $W$,  as  are  restriction  and  induction.  So  for
$I\subset S$,  a class function is  zero on $\CH_q(W_I)$ if  and only if
it  is orthogonal  to  any  $\Ind_{\CH_q(W_I)}^{\CH_q}\phi$; but  the
$\chi_{c,q}$  with  $c$ non  cuspidal  span  the  same subspace  as  the
$\Ind_{\CH_I}^{\CH}\phi$  with $I\subsetneq  S$, so  we get  the result.
\end{proof}

\begin{lemma}\label{(c)} If $\bx$ is a $d$-th root of $\bpi$,
property (c) is equivalent to
$\sum_\chi\lambda_\chi\chi(x)q^{\frac{2N-a_{\rho_\chi}
-A_{\rho_\chi}} d}=0$.
\end{lemma}
\begin{proof}
This is a simple translation of (c), using the value of
$\chi_q(T_\bx)$.
\end{proof}
\def\cox{{\text{cox}}}
\def\bcox{{\text{\bf cox}}}
We now prove the theorem when $\bG$ is split of type $A_n$.
The only cuspidal class is the class of a Coxeter element
$c$. So by \ref{(a)} $(\lambda_\chi)_\chi$ has to be equal to $a(\chi(c))_\chi$
for some $a\in\Qlbar$.
Lemma \ref{(c)} then gives $a\sum_\chi\chi(c)^2
q^{\frac{2N-a_{\rho_\chi}-A_{\rho_\chi}} d}=0$, so that $a=0$, as
all summands are non negative and at least one is non zero.

For the other types we need property (b).
\begin{lemma}\label{(b)} 
Property (b) is equivalent to the fact that
for all $i$, we have
$\sum_{\{\chi\mid a_{\rho_\chi}+A_{\rho_\chi}=i\}}
\lambda_\chi\chi(1)=0$.
\end{lemma}
\begin{proof}
Using the value of $\chi_q(T_\bpi^n)$ property (b) is equivalent to the fact that
for all $n$ we have $$\sum_i q^{n(2N-i)}
\sum_{\{\chi\mid a_{\rho_\chi}+A_{\rho_\chi}=i\}}\lambda_\chi\chi(1)=0.$$
We get the result using the linear independence of the characters of
$\bbZ$.
\end{proof}

The proof of \ref{Endomorphismes de Xpi} in the remaining types is obtained
by a computer calculation which shows that the vectors
given by \ref{(b)} and \ref{(c)} span for any $q$ the space given by 
\ref{(a)} (note that only the vectors given by \ref{(c)} depend on $q$).
\end{proof}

\vfill\eject
\section{Pieces of the Deligne-Lusztig varieties}\label{pieces}
\def\cox{{\bc}}

In this section, we introduce a technique inspired by \cite{LuCox}, which will
allow us  to compute  Harish-Chandra restrictions of  the cohomology  of some
Deligne-Lusztig varieties; we will also find a criterion for irreducibility of
a generalized Deligne-Lusztig variety (see \ref{H2l(t)}).

The    technique   is    intersecting   with    Bruhat   cells.    Let   $\CB$
be   the   variety   of   Borel   subgroups   of   $\bG$.   We   recall   from
\cite[2.2.18]{DMR}  that,  given  a  decomposition  $\bt=\bw_1\ldots\bw_k$  of
$\bt\in  B^+$   as  a  product   of  elements   of  $\bW$,  and   denoting  by
$\CO(w)$   the   $\bG$-orbit   in   $\CB\times\CB$   indexed   by   $w\in
W$,   the  variety   $\CO(\bw_1,\ldots,\bw_k)=\{  (\bB_1,\ldots,\bB_{r+1})\mid
(\bB_i,\bB_{i+1})\in\CO(w_i)\}$   depends   only   on   $\bt$   and   not   on
the   chosen    decomposition;   it   affords   two    canonical   projections
$p'(\bB_1,\ldots,\bB_{r+1})=\bB_1$ and $p''(\bB_1,\ldots,\bB_{r+1})=\bB_{r+1}$
and    the    Deligne-Lusztig    variety    is    $\bX(\bt)=\{x\in\CO(\bt)\mid
p''(x)=F(p'(x))\}$.  We  fix  an $F$-stable  Borel  subgroup  $\bB\subset\bG$,
and  for  $v\in  W$,  we  define  the  piece  $\bX^v(\bt)=\{x\in  \bX(\bt)\mid
(\bB,p'(x))\in\CO(v)\}$. We  have $\bX(\bt)=\coprod_{v\in W}  \bX^v(\bt)$, and
the action  of $\GF$  on $\bX(\bt)$ restricts  to an action  of $\BF$  on each
piece.

\begin{remark}\label{modeleXvt} If  $\bt=\bw_1\ldots\bw_k$ is  a decomposition
of  $\bt\in   B$  as  a  product   of  elements  of  $\bW$,   we  recall  from
\cite[2.2.12]{DMR} that $$\bX(\bt)= \{(g_1\bB, g_2\bB, \ldots,g_k\bB)| g_i\inv
g_{i+1}\in\bB w_i\bB,  \text{ for }  i=1,\ldots,k-1 \text{ and  } g_k\inv\lexp
Fg_1\in\bB  w_k\bB\}.$$  In this  model  we  get  $\bX^v(\bt)$ by  adding  the
condition $g_1\in\bB v\bB$.
\end{remark}

Let $\CH(W)$ be the  generic Hecke algebra of $W$ over  $\bbC[x]$. This is the
quotient of the  group algebra $\bbC[x]B$ by  the relations $(\bs+1)(\bs-x)=0$
for $\bs\in\bS$. We  denote by $T_\bb$ the image of  $\bb\in B^+$ in $\CH(W)$.
The  algebra $\CH(W)$  has a  basis  $\{T_\bw\mid \bw\in\bW\}$.  We will  also
sometimes  denote  by  $T_w$  the  elements of  this  basis.  We  will  denote
$\CH_q(W)$ the specialized algebra by the specialization $x\mapsto q$ and keep
the  notation $T_\bw$  for the  basis of  this algebra  (trying to  make clear
by  the  context  which  algebra  is  meant).  Finally  we  note  $A|T_v$  the
coefficient  of  the element  $A\in\CH(W)$  on  the  basis element  $T_v$.  We
recall that the canonical symmetrizing
form is $T_v\mapsto T_v|1$.
Since  $T_v$ and  $q^{-l(v)}T_{v\inv}$  are dual  bases for this form we  have
$A|T_v=q^{-l(v)}AT_{v\inv}|1$. With these notations, we have
\begin{proposition}\label{Ytv}
Let  $\bt\in   B^+$  and  $v\in   W$;  for   any  $m$  multiple   of  $\delta$
we  have  $\Big|\big(\UF\backslash\bX^v(\bt)\big)^{F^m}\Big|=T_vT_\bt|T_{\lexp
Fv}$, where the elements on the right-hand side are taken in the Hecke algebra
$\CH_{q^m}(W)$.
\end{proposition}
\begin{proof}
We     may    assume     $\bt\in\bW$.     Indeed,    by     \cite[2.3.3]{DMR},
if    $\bt=\bw_1\ldots\bw_k$    is    a    decomposition    as    a    product
of    elements    of     $\bW$,    and    if    $F_1$     is    the    isogeny
on    $\bG^k$   defined    by   $F_1(g_1,\ldots,g_k)=(g_2,\ldots,g_k,F(g_1))$,
then     $\bX(\bt)\simeq    \bX_{\bG^k}((w_1,\ldots,w_k),F_1)$     and    this
isomorphism   restricts  to   $$\bX^v(\bt)\simeq\coprod_{v_2,\ldots,v_k\in  W}
(\bX_{\bG^k}^{(v,v_2,\ldots,v_k)}((w_1,\ldots,w_k),F_1)).$$
Thus $(\UF\backslash\bX^v(\bt))^{F^m}\simeq
\coprod_{v_2,\ldots,v_k}\left((\bU^k)^{F_1}\backslash
(\bX_{\bG^k}((w_1,\ldots,w_k),F_1))^{(v,v_2,\ldots,v_k)}\right)^{F_1^{km}}$,
and it is also clear that:
$$\displaylines{
\sum_{v_2,\ldots,v_k}T_{(v,v_2,\ldots,v_k)}T_{(w_1,\ldots,w_k)}
|T_{\lexp{F_1}(v,v_2,\ldots,v_k)}=\hfill\cr
\hfill\sum_{v_2,\ldots,v_k}(T_vT_{w_1}|T_{v_2})(T_{v_2}T_{w_2}|
T_{v_3})\ldots(T_{v_k}T_{w_k}|T_{\lexp Fv})=T_v T_\bt|
T_{\lexp Fv}.\cr}$$
As  $F_1^{k\delta}$  is  the  smallest  power   of  $F_1$  which  is  a  split
Frobenius, we  are indeed reduced  to the  same statement for  $\bG^k$, $F_1$,
$(w_1,\ldots,w_k)$, $(v,v_2,\ldots,v_k)$.

We   then  assume   $\bt\in\bW$.  Thus   $\bX^v(t)=\{g\bB\mid  g\in\bB   v\bB,
g\inv\lexp  Fg\in   \bB  t\bB\}$.  The   map  $u\mapsto  uv\bB$   then  fibers
$\{u\in\bU\mid(uv)\inv  \lexp F(uv)\in\bB  t\bB\}$ on  $\bX^v(t)$ with  fibers
isomorphic to $\bU\cap\lexp v\bU$, since $\bB v\bB=\bU_vv\bB$,
where for $v\in W$ we set $\bU_v=\bU\cap\lexp v\bU^-$.
This fibration  is $\UF$-equivariant for the action of
$\UF$ by left multiplication on both  spaces. The  quotient by  $\UF$ is  thus obtained  by
$u\mapsto  u\inv.\lexp Fu$  which maps  the above  variety to  $\bU\cap v  \bB
t\bB\lexp F v\inv$.  As the fibers $\bU\cap\lexp v\bU$ are  connected and have
$q^{ml(w_0v)}$  fixed points  under $F^m$,  the cardinality  we seek  is thus
$q^{-ml(w_0v)}|\bU^{F^m}\cap v\bB  t\bB\lexp Fv\inv|$. Thus  the proposition
results from the  following lemma, applied with $F$ replaced  by $F^m$ and $w$
by $\lexp Fv$:
\begin{lemma}\label{TvTt|Tv1} Assume $F$  split. For $v, t, w\in  W$ and $T_v,
T_t,  T_w\in\CH_q(W)$  we  have $$T_vT_t|T_w=q^{-l(w_0v)}|(\bU\cap  v\bB  t\bB
w\inv)^F|.$$
\end{lemma}
\begin{proof}
We have $T_v T_t|T_w=T_{t\inv}T_{v\inv}|T_{w\inv}=
q^{-l(w)} T_{t\inv}T_{v\inv}T_w|1=q^{l(t)-l(w)}T_{v\inv}T_w|T_t$.
We recall that $\CH_q(W)$ may be realized as a subalgebra of $\bbC[\GF]$ via
the isomorphism $\CH_q(W)\simeq\End_\GF\Ind_\BF^\GF\bbC$. By this isomorphism,
$T_w$ corresponds to $q^{l(w)}e_\bB \dot w e_\bB$ where $\dot w$ is a
representative of $w$ in $N(\bT)^F$ and where $e_\bB$ is the idempotent
$|\BF|\inv\sum_{b\in\BF} b$. Thus:
\begin{align*} T_{v\inv} T_w|T_t&=
q^{l(v)+l(w)-l(t)} e_\bB v\inv e_\bB w e_\bB\mid e_\bB t e_\bB\\
&=|\BF|\inv q^{l(v)+l(w)-l(t)} |v\inv \BF w\cap \BF t\BF|\\
&=q^{-l(w_0)+l(v)+l(w)-l(t)}|\UF\cap v\BF t \BF w\inv|.\\
\end{align*}
The lemma follows, since $(\bU\cap  v\bB t\bB w\inv)^F=\bU^F\cap v\bB^F t\bB^F
w\inv$  which may  be  seen by  using  the uniqueness  properties  of the  Bruhat
decomposition.
\end{proof}
\end{proof}
For $\bt\in B^+$, we call support of $\bt$ the set of $\bs\in\bS$ which appear
in a decomposition of $\bt$ as a product of elements of $\bS$: this set does
not depend on the decomposition as it is not changed by a braid relation.
With this notation, we have

\begin{proposition}\label{H2l(t)}  Let $\bt\in  B^+$.  The variety  $\bX(\bt)$
is  irreducible  (in particular,  with  the  convention of  \cite[3.3.5]{DMR},
$H^{2l(\bt)}_c(\bX(\bt))$ is  given by  $h^{2l(t)}t^{l(\bt)}\Id$) if  and only
all the support of $\bt$ meets every orbit  of $F$ on $\bS$ (\ie, if the group
$\bL$ of \cite[2.3.8]{DMR} cannot be taken different from $\bG$).
\end{proposition}
\begin{proof}
We will adapt the proof of \cite[3.10(d)]{LuMa} to our case. From \ref{Ytv} we
get that if $m$ is a multiple of $\delta$, $|(\UF\backslash\bX^v(\bt))^{F^m}|$
is  the coefficient  $T_v  T_\bt| T_{\lexp  Fv}$  in $\CH_{q^m}(W)$,
this coefficient is equal to $T_\bt T_{\lexp  Fv}|T_v$ thus
$|(\UF\backslash\bX(\bt))^{F^m}|$ is  the trace of the  endomorphism $x\mapsto
T_\bt\lexp  Fx$  of $\CH_{q^m}(W)$.

\begin{lemma}\label{TtTv|Tz}  Let  $\bt\in  B^+,  v,z\in  W$;  then  the  coefficient  $T_\bt
T_v|  T_z$  is  a  polynomial  in  $q^m$ of  degree  less  than  or  equal  to
$\inf(l(\bt),l(\bt)+l(v)-l(z))$.
\end{lemma}
\begin{proof}[Proof  of lemma]  We first  show by  induction on  $l(\bt)$ that
$T_\bt T_v| T_z$ is of degree $\le l(\bt)+l(v)-l(z)$. Let $\bt''$ be a maximal
right divisor of $\bt$ in $\bW$ such that $\bt''\bv\in\bW$.

If  $\bt''=\bt$,  then  $T_\bt  T_v=T_{t''v}$ and  the  result  is  immediate:
$T_{t''v}|T_z=0$ (of degree  $-\infty$) except when $t''v=z$,  in which case
the coefficient is $1$, of degree $0=l(\bt)+l(v)-l(z)$.

Otherwise,  let   $\bt=\bt'  \bs  \bt''$  with   $\bs\in\bS$.  By  assumption
on    $\bt''$,    we    have    $\bt''\bv\in\bW$    and    $l(st''v)<l(t''v)$.
Then    $T_\bt    T_v=(q^m-1)T_{\bt'}T_{t''v}+q^m    T_{\bt'}T_{st''v}$.    As
$1+l(\bt')+l(st''v)<l(\bt)+l(v)$   and    $1+l(\bt')+l(t''v)=l(\bt)+l(v)$   we
conclude by induction, using the result  for $\bt'$ which is of smaller length
than $\bt$.

Let    $\bt\mapsto\tilde\bt$    be   the    ``reversing''    anti-automorphism
of   $B$, \ie,   the   anti-automorphism    which   extends   $\bs\mapsto\bs$   for
$\bs\in\bS$.  We  have  $T_\bt  T_v|  T_z=T_{v\inv}T_{\tilde\bt}|  T_{z\inv}=
q^{-ml(z)}(T_{v\inv}T_{\tilde\bt}T_z|  1)=q^{m(l(v)-l(z))}  (T_{\tilde\bt}T_z|
T_v)$. By the first part of the proof $T_{\tilde\bt}T_z|T_v$ is of degree $\le
l(\bt)+l(z)-l(v)$, thus $q^{m(l(v)-l(z))}(T_{\tilde\bt}T_z|T_v)$  is of degree
$\le l(\bt)$, which concludes the proof of the lemma.
\end{proof}

\begin{lemma}\label{pol. unitaire}
Let $\bt\in B^+$; then $\Trace(x\mapsto T_\bt \lexp Fx\mid\CH_{q^m}(W))$
is a polynomial in $q^m$ of  degree $l(\bt)$ and the coefficient of
$q^{ml(\bt)}$ in  this  polynomial is  the  number  of  $\bv\in\bW^F$ who are
right multiples of all elements of the support of $\bt$.
\end{lemma}
\begin{proof}
To show the lemma,  it is enough to show that  $T_\bt T_{\lexp Fv}| T_v$
is a polynomial in $q^m$ of degree  $< l(\bt)$ except if $v=\lexp Fv$ and
all  $\bs$  in the  support  of  $\bt$ divide  $\bv$  on  the left,  and  that
in  this last  case it  is  a  unitary polynomial  of  degree  $l(\bt)$.  Let
us  write  $\bt=\bt'\bs$ where  $\bs\in\bS$.  If  $l(s\lexp Fv)>l(\lexp  Fv)$,
then  by  \ref{TtTv|Tz} the  degree  of $T_\bt  T_{\lexp  Fv}|T_v=T_{\bt'}T_{s\lexp  Fv}|T_v$ is
less  than $l(\bt)$. Otherwise,  $T_\bt T_{\lexp Fv}|T_v=q^m(T_{\bt'}T_{\lexp  Fv}| T_v)+
(q^m-1)(T_{\bt'}T_{s\lexp  Fv}| T_v)$.  By \ref{TtTv|Tz},  we see  that only  the
leftmost term can contribute to  $q^{ml(\bt)}$; and by induction on $l(\bt')$,
we see that  the contribution to $q^{ml(\bt)}$ of $T_\bt  T_{\lexp Fv}|T_v$ is
$T_{\lexp Fv}|T_v$  if all $\bs$  in the support of  $\bt$ divide on  the left
$\lexp Fv$, and is $0$ otherwise. The result follows.
\end{proof}

From  the  last  lemma,   $|(\UF\backslash\bX(\bt))^{F^m}|$  is  a  polynomial
of  degree  $l(\bt)$  in  $q^m$,  unitary  if  and  only  if  the  support  of
$\bt$  meets  every   $F$-orbit  in  $\bS$.  As   all  irreducible  components
of  $\UF\backslash\bX(\bt)$   have  the   same  dimension,   since  $\bX(\bt)$
is  the  transverse  intersection  of  the   graph  of  $F$  with  the  smooth
irreducible variety  $\CO(\bt)$, this  variety is irreducible  if and  only if
$|(\UF\backslash\bX(\bt))^{F^m}|$ is a unitary polynomial in $q^m$.

To prove the proposition
it  remains  to   check  that  $\bX(\bt)$  is  irreducible  if   and  only  if
$\UF\backslash\bX(\bt)$  is.  The``only  if''  part   is  clear;  to  see  the
``only''  part,  we   may  follow  the  arguments   of  \cite[4.8]{LuCox}:  if
$\UF\backslash\bX(\bt)$  is   irreducible,  the   set  $\pi$   of  irreducible
components   of  $\bX(\bt)$   is  a   single   orbit  under   $\UF$,  so   its
cardinality  is  a  power  of  $p$.   The  set  $\pi$  is  in  bijection  with
the  set of  irreducible  components of  the (smooth) compactification  $\bX(\underline
s_1,\ldots,\underline s_r)$ (see  \cite[2.3.4]{DMR}). But the $\GF$-stabilizer
of $(\bB,\ldots,\bB)\in\bX(\underline s_1,\ldots,\underline  s_r)$ is $\bB^F$,
thus the orbit of $(\bB,\ldots,\bB)$ (and a fortiori the number of irreducible
components of  $\bX(\underline s_1,\ldots,\underline s_r)$) has  cardinality a
divisor of $|\bG^F/\bB^F|$, which is prime to $p$, whence the result.
\end{proof}
By \ref{Ytv}, we see that the variety $\bX^v(\bw)$ is non-empty if and only if
$T_v T_\bw|T_{\lexp Fv}\ne 0$. We  shall study this condition, especially when
$v$ is $F$-stable. In  what follows, we will denote by  $\le$ the Bruhat order
on $W$.
\begin{proposition}\label{Yv non vide}
Assume  that  $\bw\in  B^+$  is  of the  form  $\bw=\bw_1\ldots  \bw_k$  where
$\bw_i\in\bW$ have  mutually disjoint support.  Then $T_v T_\bw| T_v\ne  0$ is
equivalent to $T_v T_{\bw_i}| T_v\ne 0$ for all $i$.
\end{proposition}
\begin{proof}
By induction on $k$,  it is enough to show the case  $k=2$ of the proposition.
By the isomorphism of the Hecke algebra with a subalgebra of the group algebra
of  $\GF$, we  have $T_v  T_\bw| T_v\ne  0$ if  and only  if $\bB  v\bB w_1\bB
w_2\bB\supset \bB v \bB$. We then use the following lemma:
\begin{lemma}\label{BwBw'B} For  $w,w'\in W$  we have: $\bB  w\bB w'\bB\subset
(\coprod_{v'\le w'} \bB wv'\bB)\cap(\coprod_{v\le w} \bB vw'\bB)$.
\end{lemma}
\begin{proof}
The inclusion in  \eg, the left union  is an easy induction  on $l(w')$, using
the exchange lemma.
\end{proof}
Thus  $\bB  v\bB w_1\bB  w_2\bB=\coprod  _{v_1}\bB  vv_1\bB w_2\bB$  for  some
$v_1\le w_1$, and in  turn this last union is a union of  double cosets of the
form $\bB  vv_1 v_2\bB$, where  $v_2\le w_2$;  now the assumption  on supports
implies that $vv_1  v_2=v$ if and only if $v_1=v_2=1$.  Since $v_1=1$ occurring
is equivalent to $\bB v \bB w_1\bB  \supset \bB v\bB$ and then in turn $v_2=1$
occurring  is equivalent  to $\bB  v \bB  w_2\bB\supset \bB  v\bB$, we  get the
proposition.
\end{proof}
\begin{lemma}\label{z>xy} If $v,w,x\in W$ and if $T_vT_w|T_x\ne 0$, then $x\ge
vw$.
\end{lemma}
\begin{proof}
The  condition  is  equivalent  to  $T_wT_{x\inv}|T_{v\inv}\neq  0$.  Applying
\ref{BwBw'B},   this   implies   $v\inv=wx'$   with   $x'\leq   x\inv$,   \ie,
$vw=x^{\prime-1}$ with $x^{\prime-1}\leq x$.
\end{proof}
\begin{lemma}\label{vsts}  Let $v,t\in  W$ where  $t$ is  a reflection  of root
$\alpha>0$. Then $T_vT_t|T_v\ne 0$ if and only if $v\alpha<0$.
\end{lemma}
\begin{proof}
By  \cite[1.2   and  1.12]{Dy},  $v\alpha<0$   if  and  only  if   $vt<v$.  If
$T_vT_t|  T_v\ne 0$,  then  by  \ref{z>xy} we  have  $v>vt$ thus  $v\alpha<0$.
Conversely,  if  $v\alpha<0$,  we  will  show  by  induction  on  $l(t)$  that
$T_vT_t| T_v\ne  0$. If  $l(t)=1$, then $v=v't$  with $l(v)=l(v')+1$.  We have
then  $T_vT_t=(q-1)T_v+qT_{vt}$,  thus  $T_vT_t|T_v=q-1\ne 0$.  Otherwise,  by
\cite[1.4]{Dy},  we may  write  $t=at'a$ where  $a\in  S$ and  $l(t)=l(t')+2$.
We  have:  $T_vT_{at'a}|  T_v=q^{-l(v)}T_vT_a T_{t'}  T_a  T_{v\inv}|1$;  when
$va>v$  this  is  equal  to $qT_{va}T_{t'}|T_{va}$.  Otherwise,  it  is  equal
to      $q^{-l(v)}(q^2(T_{va}T_{t'}T_{av\inv}|1)+(q-1)^2(T_vT_{t'}T_{v\inv}|1)
+q(q-1)(T_{va}T_{t'}T_{v\inv}|1)+q(q-1)(T_vT_{t'}T_{av\inv}|1))$  whose  first
term is  equal to $q^3T_{va}T_{t'}|T_{va}$.  Since the structure  constants of
the Hecke algebra  are polynomials which positive highest  coefficient, we see
in both cases that $T_vT_{t}| T_v$ will be non zero if $T_{va}T_{t'}|T_{va}\ne
0$ is non zero. Since $t'$ is a  reflection of root $a\alpha$, we see by
induction, that this  coefficient is non zero if  $va(a\alpha)<0$, \ie,
$v\alpha<0$, whence the result.
\end{proof}

Recall that an element  $w\in W$ is reduced-$I$, with $I\subset S$,  if it is of
minimal length in its coset $wW_I$.
To continue  our study, we  define $E_W(w)=\{w_0v\in W\mid T_v  T_w|T_v\ne 0\}$.
With this notation, we have
\begin{lemma}\label{calculE}
Let $I\subset S$  be $F$-stable. Assume that  $w\in W$ is of  the form $w=sw'$
with $s\in  S-I$ and and $w'\in  W_I$. Then $E_W(w)$ consists  of the products
$v_1 v_2$ where $v_2\in E_{W_I}(w')$ and  where $v_1$ is a reduced-$I$ element
such that $l(v_1 v_2 s)>l(v_1 v_2)$.
\end{lemma}
\begin{proof} By proposition \ref{Yv non vide}, $T_v T_w|T_v\ne 0$ if and only
if $T_v  T_s|T_v\ne 0$ and  $T_v T_{w'}|T_v\ne 0$.  Let $w_0 v=v_1  v_2$ where
$v_2\in W_I$ and $v_1$ is reduced-$I$.  Then we have $v=(w_0 v_1 w_0^I).(w_0^I
v_2)$.  Note that  $w_0  v_1  w_0^I$ is  still  reduced-$I$;  it follows  that
$T_v=T_{w_0  v_1  w_0^I}T_{w_0^I v_2}$  and  that  the condition  $T_{w_0  v_1
w_0^I}T_{w_0^I v_2}T_{w'}|T_{w_0  v_1 w_0^I}T_{w_0^I v_2}\ne 0$  is equivalent
to  $T_{w_0^I  v_2}T_{w'}|T_{w_0^I  v_2}\ne  0$. It  remains  to  express  the
condition  $T_{w_0  v_1  v_2}T_s|T_{w_0  v_1 v_2}\ne  0$;  this  condition  is
equivalent to $l(w_0  v_1 v_2 s)<l(w_0 v_1 v_2)$, which  in turn is equivalent
to $l(v_1 v_2 s)>l(v_1 v_2)$.
\end{proof}
Note that $v_1=1$  and $v_2$ arbitrary in $E_{W_I}(w')$ satisfy the above
condition, so that $E_W(w)\supset E_{W_I}(w')$.

We will apply \ref{calculE} in a more specific situation where the following
holds:
\begin{proposition}\label{calculE2}
Under the assumptions of \ref{calculE}, assume in addition that $S=I\cup\{s\}$,
that  there is a unique  $s'\in I$ which  does not commute with  $s$ and that
$ss's=s'ss'$. Assume also  that any $v\in E_{W_I}(w')$  whose support contains
$s'$ is  such that $s'$  is not  in the  support of  $s'v$. Then
$E_W(w)=E_{W_I}(w')\cup\{\,sv\mid v\in E_{W_I}(w')\text{ and }s'v<v\,\}$.
\end{proposition}
\begin{proof}
We take  $v_1v_2\in E_W(w)$  as in  lemma \ref{calculE};  if $v_1\ne  1$, then
$l(v_1s)<l(v_1)$ since  $v_1$ is reduced-$I$.  It follows that $v_2$  does not
commute with $s$,  thus the support of  $v_2$ must contain $s'$.  We claim that
$v_1=s$; otherwise, as  $v_1$ is reduced-$I$, it would end with  $s's$ since
$s''s=ss''$ for $s''\neq s'$;
but then $s'sv_2s$  would not be reduced since by the  assumption of the
proposition $v_2 s$ has a reduced expression starting with $s's$ as $s$ commutes
with all terms of a reduced expression for $v_2$ excepted $s'$.
\end{proof}
For  $I\subset S$,  we denote  by  $\Phi_I$ the  corresponding parabolic  root
subsystem, and we  denote by $\bL_I$ the Levi subgroup  generated by $\bT$ and
$\{\,\bU_\alpha\,\}_{\alpha\in\Phi_I}$.  We denote  $\bB_I$ (resp.  $\bB^-_I$,
$W_I$,  $\bU_I$, $\bU^-_I$)  the  intersection with  $\bL_I$  of $\bB$  (resp.
$\bB^-$, $W$,  $\bU$, $\bU^-$), by  $\bP_I$ the parabolic  subgroup $\bL_I\bB$
and  by  $\bU_{\bP_I}$  its  unipotent  radical.  We  will  use  the  following
proposition in the proof of \ref{union of pieces}:
\begin{proposition}\label{B(w) inter U-}
Let   $I_1,\ldots,I_k$  be   mutually  disjoint   subsets  of   $S$  and   let
$x_i\in\bL_{I_i}$.  Then   the  condition   $x_1\ldots  x_k\in   \bU^-\bB$  is
equivalent to $x_i\in\bU^-_{I_i}\bB_{I_i}$ for all $i$.
\end{proposition}
\begin{proof}
If $k=1$,  let us  write $x_1=u\dot vb$  with $u\in  \bU_{I_1}^-$, with $\dot
v$  a  representative  of  $v\in   W_{I_1}$  and  $b\in  \bB_{I_1}$.  As
$\bU_{I_1}^-\subset  \bU^-$ and  $\bB_{I_1}^-\subset  \bB$, the  existence
of the Bruhat decomposition with respect to the pair of Borel
subgroups $(\bB^-,\bB)$, which is obtained by multiplying on the left by $w_0$ the 
classical Bruhat decomposition, implies that $v=1$.

By  induction on  $k$ it  is enough  to prove  the statement  for $k=2$.
Let $I=I_1\cup I_2$; we have $\bU_I^-=\bU_{I_1}^-\bU_{I_2}^-$ and
$\bB_I=\bB_{I_1}\bB_{I_2}$. From the case $k=1$ we get that $x_1x_2\in\bU^-\bB$ is
equivalent to $x_1x_2\in \bU_{I_1}^-\bU_{I_2}^-\bB_{I_1}\bB_{I_2}$. As $\bB_{I_1}$
normalizes $\bU_{I_2}^-$, this is equivalent to
$x_1x_2\in (\bU_{I_1}^-\bB_{I_1}).(\bU_{I_2}^-\bB_{I_2})$.
Let us write $x_1x_2=y_1y_2$ according to this decomposition; as $\bL_{I_1}\cap\bL_{I_2}=\bT$
we get $x_i\in y_i\bT\subset\bU_{I_i}^-\bB_{I_i}$ for $i=1,2$.
\end{proof}
\begin{lemma}\label{w=v1..vk}  Let   $w=v_1\ldots  v_k\in  W$  be   such  that
$l(v_1)+\ldots+l(v_k)=l(w)$  and  let  $\dot w,\dot  v_1,\ldots,\dot  v_k$  be
representatives  in $N_\bG(\bT)$  such that  $\dot w=\dot  v_1\ldots\dot v_k$.
Then $\bU_w\dot w=\bU_{v_1}\dot v_1\ldots \bU_{v_k}\dot v_k$.
\end{lemma}
\begin{proof}
We      have       $\bU_{v_1}\dot      v_1\ldots       \bU_{v_k}\dot      v_k=
\bU_{v_1}\lexp{v_1}\bU_{v_2}\ldots\lexp{v_1\ldots   v_{k-1}}\bU_{v_k}\dot  w$,
and    we    have   $\bU_w=\prod_{\alpha\in    N(w^{-1})}\bU_\alpha$,    where
$N(w)=\{\,\alpha>0\mid\lexp  w\alpha<0\}$.  Let $\coprod$  represent  disjoint
union. The  lemma is thus a  consequence of $N(w\inv)=\coprod_i v_1\ldots
v_{i-1}(N(v_i\inv))$,  which itself  is obtained by  iterating the  well known
formula: $l(x)+l(y)=l(xy) \Leftrightarrow N(xy)=y\inv(N(x))\coprod N(y)$.
\end{proof}
\begin{corollary}\label{BwB  inter   U-}  Let  $I_1,\ldots,I_k$   be  disjoint
parts  of $S$,  and  let $v_i\in  W_{I_i}$. Then  $\bB  v_1\ldots v_k  \bB\cap
\bU^-=\prod_i((\bB_{I_i}v_i\bB_{I_i})\cap\bU^-)$.
\end{corollary}
\begin{proof}   As  in   \ref{B(w)  inter   U-},   it  is   enough  to   prove
the   result    for   $k=2$.    By   \ref{w=v1..vk}    we   have    $\bB   v_1
v_2\bB=\bU_{v_1v_2}v_1v_2\bB=\bU_{v_1}\dot  v_1   \bU_{v_2}\dot  v_2\bB$.  Let
$x\in\bB  v_1  v_2\bB\cap \bU^-$  and  write  accordingly $x=x_1  x_2b$  where
$x_1\in\bU_{v_1}\dot  v_1$,   $x_2\in\bU_{v_2}\dot  v_2$  and   $b\in\bB$.  We
have  $x_1  x_2\in \bU^-\bB$  thus  by  \ref{B(w)  inter U-}  $x_1=u_1  b_1\in
\bU^-_{I_1}\bB$. We have $u_1=x_1b_1\inv\in \bU_{v_1}\dot v_1\bB\cap \bU^-=\bB
v_1\bB\cap \bU^-$.  As $x\in \bU^-$  we have also  $b_1 x_2 b\in  \bU^-$, thus
$b_1x_2 b\in \bB v_2\bB\cap \bU^-$.
\end{proof}
\begin{proposition}\label{piece isolee}
Let $\bw\in  B^+$, and  let $I$ be an $F$-stable  subset of $S$.  Then for  any $v\in
W$,  the left  multiplication  action of  $\bP_I^F$  on $\bX(\bw)$  stabilizes
$\coprod_{v'\in W_Iv}\bX^{v'}(\bw)$.
\end{proposition}
\begin{proof}An   element    $(g_1\bB,\ldots,g_r\bB)\in   \bX(\bw)$    is   in
$\bX^v(\bw)$ if and only if $g_1\in\bB  v\bB$. If $p\in \bB w'\bB$ with $w'\in
W_I$, then by  \ref{BwBw'B} we have $pg_1\in \bB w''  v\bB$ with $w''\leq w'$,
thus $w''\in W_I$, whence the result.
\end{proof}
In the next proposition, for $I\subset S$ we denote
by $B^+_I$ the submonoid of $B^+$ generated by $\bI=\{\bs\in\bS\mid s\in I\}$.
\begin{proposition}\label{union of pieces}
Under  the  assumptions  of  \ref{piece   isolee},  assume  in  addition  that
$\lexp{\bw_0}\bw=\bs\bw'$  where  $\bs\notin  I$ and  where  $\bw'\in  B^+_I$.
Let   $\bU_{\bP_I}$  be   the   unipotent  radical   of   $\bP_I$.  Then,   if
$\bw_0^I\in\bW$  lifts  $w_0^I$,  for  any  $i$  we  have  an  isomorphism  of
$\bL_I^F\times\genby{F^\delta}$-modules
$$H^i_c(\left(\coprod_{v\in W_Iw_0}\bX^v(\bw)\right)/\bU_{\bP_I}^F)
\xrightarrow\sim
H^{i-2}_c(\bX_{\bL_I}(\lexp{\bw_0^I}\bw'))(-1)\oplus
H^{i-1}_c(\bX_{\bL_I}(\lexp{\bw_0^I}\bw'))$$
\end{proposition}
\begin{proof}
To  simplify  the notation we  write  just  $\bL$,  $\bP$ for  $\bL_I$,  $\bP_I$,
and  we  set $\bY=\coprod_{v\in  W_Iw_0}\bX^v(\bw)$.  Let  us see  first  that
the  proposition  follows   from  its  special  case   where  $\bw\in\bW$.  If
$\bw=\bw_1\ldots \bw_k$ is  a decomposition of $\bw$ as a  product of elements
of $\bW$, we have $\bw'=\bw'_1\lexp{\bw_0}\bw_2\ldots\lexp{\bw_0}\bw_k$, where
$\lexp{\bw_0}\bw_1=\bs\bw'_1$.  Using \cite[2.3.3]{DMR}  as  in the  beginning
of  the proof  of  \ref{Ytv},  we apply  the  proposition  with
$\bG^k$, $F_1$, $\bL^k$,
$(\bw_1,\ldots,\bw_k)$ and $(\bs,1,\ldots,1)$ replacing respectively $\bG$, $F$,
$\bL$, $\bw$ and $\bs$. 
If we set $\bY'=\coprod_{v_1,\ldots v_k\in
W_Iw_0}\left(\bX_{\bG^k}^{(v_1,\ldots,v_k)}((\bw_1,\ldots,\bw_k),F_1)\right)$
we  obtain   that  the  cohomology  groups   of  $\bY'/(\bU_\bP^k)^{F_1}$  are
sums of  those of $\bX_{\bL^k}((\lexp{\bw_0^I}\bw'_1,\bw_2\ldots,\bw_k),F_1)$.
This  last   variety  is   isomorphic  to   $\bX_\bL(\lexp{\bw_0^I}\bw')$.  On
the  other  hand  $(\bU_\bP^k)^{F_1}\simeq\bU_\bP^F$   and  $\bY'$  is  formed
from   pieces  from   $\bY$.   Indeed,   by  the   beginning   of  the   proof
of   \ref{Ytv},  we   have   $\bY=   \coprod_{v\in  W_I,v_2,\ldots,v_k\in   W}
(\bX_{\bG^k}^{(v,v_2,\ldots,v_k)}((w_1,\ldots,w_k),F_1))$.  We show that
the  only  non-empty  pieces  of  $\bY$  are  those  such  that  $v_i\in  W_I$
for  all  $i$,  \ie,  those  of  $\bY'$:  the  $(v,v_2,\ldots,v_k)$  piece  is
non-empty if and only if $T_{(w_1,\ldots,w_k)}T_{(v_2\inv,v_3\inv,\ldots,\lexp
Fv\inv)}|T_{(v\inv,v_2\inv,\ldots,v_k\inv)}\ne   0$.   As   $w_2,\ldots,w_k\in
\lexp{w_0}W_I$ and $v\in W_Iw_0$ the non-vanishing of this coefficient implies
that $v_2,\ldots,v_k\in  W_Iw_0$; indeed, proceeding  by induction on  $i$, if
$v_{i+1}\in W_Iw_0$  with $i\geq2$, then the  product $T_{w_i}T_{v_{i+1}\inv}$
involves only $T_y$ for $y\in W_Iw_0$, thus $v_i\in W_Iw_0$.

We  thus assume  now  $\bw\in\bW$. We  use then  the  model \ref{modeleXvt}  of
$\bX^v(\bw)$ taking $k=2$, $w_1=\lexp{w_0}s$ and $w_2=\lexp{w_0}w'$:
$$\bX^v(\bw)= \{(g_1\bB, g_2\bB)|
g_1\inv g_2\in\bB \lexp{w_0}s\bB,
g_2\inv\lexp Fg_1\in\bB \lexp{w_0}w'\bB,
\text{ and } g_1\in\bB v\bB \}.$$
Let $\dot w_0$ be a rational representative of $w_0$.
Taking $g_1\dot w_0$ and $g_2\dot w_0$ as variables we get
$$\bX^v(\bw)=\{(g_1\bB^-,g_2\bB^-)\mid g_1\in\bB vw_0\bB^-,\;
g_1\inv.g_2\in\bB^-s\bB^-,\;
g_2\inv.\lexp Fg_1\in\bB^- w'\bB^-\}.$$

By arguing as in the proof of \ref{piece isolee} we see that
$\cup_{v\in W_Iw_0}\bB vw_0\bB^-=\bP\bB^-$. Thus
\begin{equation}
\bY=\{(g_1\bB^-,g_2\bB^-)\mid  g_1\in\bP\bB^-,\, g_2\in\bG,\,
g_1\inv.g_2\in\bB^- s\bB^-,\,g_2\inv.\lexp Fg_1\in\bB^-
w'\bB^-\}.\tag{1} \end{equation}
The  action of $\bP^F$ is by left multiplication. In
the  following, we fix $F^\delta$-stable representatives, denoted $\dot s$ and
$\dot w'$ of $s$ and $w'$. For $v\in W$, let $\bU_v^-=\bU^-\cap\lexp v\bU$.
\begin{lemma}\label{chgt of variable}
The  variety  $\bX=\{p\in\bP\mid   p\inv\lexp  Fp\in\bB^-  sw'\bB^-\}$  admits
natural  actions of  $\bB^-$  by right  multiplication and  of  $\PF$ by  left
multiplication. The map $p\mapsto (p\bB^-,pu_p\dot s\bB^-)$ where $u_p$ is the
unique element of $\bU^-_s$ such that $p\inv\lexp Fp\in u_p\dot s\bB^-w'\bB^-$
defines  a   $\PF$  and   $F^\delta$  equivariant  isomorphism   of  varieties
$\bX/\bB^-_I\xrightarrow\sim\bY$.
\end{lemma}
\begin{proof}
The  existence  and  uniqueness  of  $u_p$  come  from  \ref{w=v1..vk}:  we  have
$$\bB^-sw'\bB^-=\bU^-_{sw'}\dot  s\dot  w'\bB^-=  \bU^-_s\dot  s\bU^-_{w'}\dot
w'\bB^-=\bU^-_s\dot s\bB^-w'\bB^-$$  where the  $\bU^-_s$ part is  unique. The
image of the  map $p\mapsto (p\bB^-,pu_p\dot s\bB^-)$ is easily  checked to be
in  the  model  (1)  of  $\bY$. If  $b\in\bB_I^-$,  the  element  $u_{pb}$  is
determined by $b\inv u_p\dot s\bB^-=u_{pb}\dot s\bB^-$. We have thus $pu_p\dot
s\bB^-=pbu_{pb}\dot  s\bB^-$, which  shows that  $p$  and $pb$  have the  same
image in $\bY$.

Conversely, given $(g_1\bB^-,g_2\bB^-)\in\bY$,  the equality $g_1\bB^-=p\bB^-$
defines   $p\in\bX$   up to   right   translation   by   $\bB_I^-$.   We   must
check that $g_2\bB^-=pu_p\dot    s\bB^-$.    By    definition
of    $\bY$,     we    have     $g_2\in\lexp    Fg_1\bB^-w^{\prime-1}\bB^-\cap
g_1\bB^-s\bB^-=   \lexp   Fp\bB^-w^{\prime-1}\bB^-\cap  p\bB^-s\bB^-$   whence
$p\inv    g_2\in     p\inv\lexp    Fp\bB^-w^{\prime-1}\bB^-\cap    \bB^-s\bB^-
\subset    u_p\dot   s\bB^-w'\bB^-w^{\prime-1}\bB^-\cap    \bB^-s\bB^-$;   but
$\bB^-w'\bB^-w^{\prime-1}\bB^-$   is  a   union  of   double  cosets   of  the
form   $\bB^-   v\bB^-$,   where   $v\in  W_I$.   Thus   $u_p\dot   s\bB^-\dot
w'\bB^-w^{\prime-1}\bB^-\cap \bB^-s\bB^-=u_p\dot s\bB^-$, whence the result.
\end{proof}
Let us decompose $p\in\bX$  as $ul$, with  $u\in\bU_\bP$ and $l\in\bL$.
The action  of $\bB_I^-$
does  not  change  the  component  $u$  thus  the  quotient  of  $\bX/\bB_I^-$
by   $\bU_\PF$   is  realized   by   the   Lang  map   $(u,l)
\mapsto(u\inv.\lexp Fu,l)$.
If   we   take  $\lexp{l\inv}(u\inv.\lexp Fu)$ and $l$  as   variables
we   get  $$\bY/\bU_\bP^F\simeq\{u\in\bU_\bP,l\in\bL\mid   ul\inv.\lexp  Fl\in
\bB^-sw'\bB^-\}/\bB_I^-,$$ where the action of $b\in\bB_I^-$ is by conjugation
by $b\inv$ on $u$ and by right  multiplication on $l$, and where the action of
$\LF$ is by left multiplication on $l$.

\begin{lemma} For $u\in\bU_\bP$, $l\in\bL$, the condition $ul\in\bB^-sw'\bB^-$
is   equivalent  to   $u\in\bB^-_{\{s\}}   s\bB^-_{\{s\}}$  and   $l\in\bB^-_I
w'\bB^-_I$.
\end{lemma}
\begin{proof}     We      have     $\bB^-sw'\bB^-=\bU^-_s\dot     s\bU^-_w\dot
w\bB^-_I\bU_\bP^-=  \bU^-_s\dot  s\bB^-_Iw\bB^-_I\bU_\bP^-$.  Thus  using  that
$l$  normalizes  $\bU_\bP^-$,  we   see  that  there  exists  $u'\in\bU_\bP^-$
such   that   $uu'\in   \bU^-_s\dot  s\bB^-_Iw\bB^-_I   l\inv$.   Thus   there
exists   $l_s\in\bU^-_s\dot   s\subset\bL_{\{s\}}$,   $l'\in\bB^-_I   w\bB^-_I
l\inv\subset\bL$  such   that  $uu'=l_sl'$.   We  may  then   apply  \ref{B(w)
inter  U-}   with  $k=2$,   $I_1=\{s\}$  and   $I_2=I$  (exchanging the roles  of
$\bB$   and  of   $\bB^-$),  and   we  get   $l_s\in\bU_s\bB_{\{s\}}^-$.  Thus
there   exists  $u_s\in\bU_s$   such   that   $l_s\in  u_s\bB^-_{\{s\}}$.   We
have   $u_s\inv  u\in\bB^-_{\{s\}}\bL   u^{\prime-1}\subset\bL\bU^-_\bP$  thus
$u_s\inv  u\in\bU_\bP\cap\bL\bU^-_\bP=\{1\}$ thus  $u=u_s\in\bU_s$, and  since
$l_s\in\bU^-_s\dot   s$   we   even  have   $u\in\bU^-_s\dot   s\bB^-_{\{s\}}=
\bB^-_{\{s\}}s\bB^-_{\{s\}}$. The condition $uu'=l_s  l'$ becomes thus $u'=b_s
l'$ for some $b_s\in\bB^-_{\{s\}}$; as $\bT\bU^-_\bP\cap\bL=\bT$, this implies
$l'\in\bT$ thus $l\in\bB^-_I w'\bB^-_I$ q.e.d.
\end{proof}
As we have $\bU\cap\bB^-_{\{s\}}s\bB^-_{\{s\}}=\bU_s^*\subset\bU_\bP$, we get thus
$\bY/\bU_\bP^F\simeq\{u\in\bU_s^*,l\in\bL\mid    l\inv\lexp   Fl\in    \bB^-_I
w'\bB^-_I\}/\bB^-_I$ where
the action  of $\LF$  is by  left multiplication on $l$ and
the action of $b\in\bB^-_I$ is by right multiplication
on $l$ and conjugation by $b\inv$ on $u$. Note that, as $\bU_I^-$ centralizes $\bU_s$
since no root in $\Phi_I^-$ can add to the simple root corresponding to $s$, the action of
$\bB^-_I$  on $u$  is through $\bT$.

We have $\bX_{\bL_I}(\lexp{w_0^I}w')=\{l\in\bL\mid   l\inv\lexp  Fl\in   \bB^-_I
w'\bB^-_I\}/\bB^-_I$.  We  conclude  arguing  as  in  \cite[3.2.10]{DMR}. To simplify the notation we
write $w''$ for $\lexp{w_0^I}w'$. Let
$\tilde\bY$ be  the variety $\{u\in\bU_s,l\in\bL\mid l\inv\lexp  Fl\in \bB^-_I
w'\bB^-_I\}/\bB^-_I$; the projection $\pi:\tilde\bY\to\bX_{\bL_I}(w'')$ defined
by $(u,l)\mapsto  l$ is a fibration  by affine lines, and  $\pi$ restricted to
$\tilde\bY-\bY/\bU_\bP^F$ is an  isomorphism. Let $i$ be  the closed inclusion
$\tilde\bY-\bY/\bU_\bP^F\hookrightarrow\tilde\bY$   and  $j$   the  open   inclusion
$\bY/\bU_\bP^F\hookrightarrow\tilde\bY$.  If we  make $\LF$ act by left multiplication
on $\tilde\bY$
then  $i$, $j$  and  $\pi$  are $\LF$-equivariant.  Let  $\Lambda_\bX$ be  the
constant sheaf  $\Qlbar$ on a variety  $\bX$; we have an  exact sequence $0\to
j_!\Lambda_{\bY/\bU^F_\bP}\to\Lambda_{\tilde\bY}\to  i_!\Lambda_{\tilde\bY-\bY/\bU^F_\bP}\to
0$. Using that $R\pi_!\Lambda_{\tilde\bY}\simeq\Lambda_{\bX_{\bL_I}(w'')}[-2](-1)$
and that $R\pi_!i_!\Lambda_{\tilde\bY-\bY/\bU^F_\bP}\simeq\Lambda_{\bX_{\bL_I}(w'')}$,
we deduce a distinguished triangle
$R\pi_!j_!\Lambda_{\bY/\bU^F_\bP}\to\Lambda_{\bX_{\bL_I}(w'')}[-2](-1)
\xrightarrow\partial \Lambda_{\bX_{\bL_I}(w'')}\rightsquigarrow$ where
$\partial\in\Ext^2(\Lambda_{\bX_{\bL_I}(w'')},
\Lambda_{\bX_{\bL_I}(w'')})=H^2(\bX_{\bL_I}(w''),\Qlbar)$.   All    maps   being
$\GF$-equivariant, we  even have  $\partial\in H^2(\bX_{\bL_I}(w''),\Qlbar)^\GF
\simeq    H^{2l-2}_c(\bX_{\bL_I}(w''),\Qlbar)^\GF$   where    the   isomorphism
comes  from   the  smoothness  of   the  variety  $\bX_{\bL_I}(w'')$.   But  by
\cite[3.3.14]{DMR}  we  have $H^{2l-2}_c(\bX_{\bL_I}(w''),\Qlbar)^\GF=0$,  thus
$\partial=0$ and the distinguished triangle gives an isomorphism
$R\pi_!j_!\Lambda_{\bY/\bU^F_\bP}\simeq\Lambda_{\bX_{\bL_I}(w'')}[-2](-1)
\oplus\Lambda_{\bX_(\bL_I)(w'')}[-1]$ whence
$H^i_c(\bY/\bU^F_\bP)\simeq H^{i-2}_c(\bX_{\bL_I}(w''))(-1)\oplus
H^{i-1}_c(\bX_{\bL_I}(w''))$ as wanted.
\end{proof}
In the proof of proposition \ref{variety restante} we will use the next lemma.
\begin{lemma}\label{Uv1Uv2...}
Let $\dot w_1,\ldots,\dot w_k$ be representatives
in  $N_\bG(\bT)$ of $w_1,\ldots,w_k\in W$; for any $u_1,\ldots,u_k\in \bU$, there exist
unique $u'_i\in \bU_{w_i}$ such that for all $i$ we have
$u_1\dot w_1\ldots u_i\dot w_i\in u'_1\dot w_1\ldots u'_i\dot w_i\bU$. This
defines a morphism $\bU^k\to\prod_{i=1}^k\bU_{w_i}$.
\end{lemma}
\begin{proof}
It is known that for $v\in W$, the equality $u\dot v=u'\dot vu''$ with $u\in\bU$, $u'\in\bU_v$ and
$u''\in\bU\cap\lexp v\bU$ defines an isomorphism $\bU\xrightarrow\sim \bU_v\times \bU\cap\lexp v\bU$.
The lemma is a consequence of this fact by induction on $k$.
\end{proof}
The elements $w$ that  we will handle in this paper will  have $E_W(w)$ of the
form $W_I  w_0\cup \{v\}$, where $v$ satisfies the assumptions of the next proposition.
\begin{proposition}\label{variety restante} Let
$\bw=\bw_1\ldots \bw_k$ be a
decomposition of $\bw$ with $\bw_i\in \bW$; let $v\in W^F$,
and let $I$ be an $F$-stable subset of $S$; we write $\bP$ for $\bP_I$. Then
\begin{enumerate}
\item \label{restante1}  For all $i$  we      have       an      isomorphism
of       $\genby       F$-modules:       $$H^i_c(\bX^v(\bw))^{\bU_\bP^F}\simeq
H_c^{i+2l(w_0v)}(\bZ^v_\bw)(l(w_0v)),$$
where $(l(w_0v))$ is a ``Tate twist'' and
$$\bZ^v_\bw=\{(y,x,u_1,\ldots,u_k)\in
\bU_\bP\times\bU_I\times\prod_i\bU_{w_i}\mid  yx\inv\lexp Fx\in  \dot vu_1\dot
w_1\ldots u_k\dot w_k\bB v\inv\}.$$
\item\label{restante4}  Let  $\overline\bZ^v_\bw$  be  the  variety
$\{(y,x,u_1,\ldots,u_k)  \in\bZ^v_\bw\mid y\in\bU_\bP\cap\lexp  v\bU^-\}$.
The      map       $\bZ_\bw^v\to\overline\bZ_\bw^v$      given      by
$(y,x,u_1,\ldots,u_k)\mapsto(y_2,x,u'_1,\ldots,u'_k)$  where $y_2$  is defined
by $y=y_1y_2$ with $y_1\in\bU_\bP\cap\lexp  v\bU$ and $y_2\in \bU_\bP\cap\lexp
v\bU^-$ and where the $u'_i$  are defined by $\lexp{\dot v\inv}y_1\inv u_1\dot
w_1\ldots u_i\dot w_i\in u'_1\dot w_1\ldots u'_i\dot w_i\bU$ for any $i$ (\cf\
\ref{Uv1Uv2...}), is a  fibration whose  fibers are all  isomorphic to
$\lexp {v\inv}\bU_\bP\cap\bU$.
\item \label{restante2}  If in addition  $v$ is  the unique element  of $W_Iv$
such  that  $\bX^v(\bw)$  is  nonempty,   there  is  an  action  of  $\bL_I^F$
on   $\bZ^v_\bw$  such   that  the   isomorphism  of   \ref{restante1}  is
$\bL_I^F$-equivariant (for  the action of  $\bL_I^F$ on $\bX^v(\bw)$  given by
\ref{piece isolee}).
\item  \label{restante3} If  in  addition
$\lexp{v\inv}\bU_I\subset\bU^-$ and
$\proj_{\bU_I}(\bU\cap  v\bB
w_1\bB\ldots  w_k\bB v\inv)  =\prod_{s\in  I}\bU^*_s$ (where  we have  denoted
$\proj_{\bU_I}$  the  natural   projection  $\bU\to\bU_I$),  then  the
projection  $\pi:(y,x,u_1,\ldots,u_k)\mapsto x$  is  an epimorphism
$\bZ^v_\bw\to\bX_{\bL_I}(c)$, where $c$  is a Coxeter element
of  $W_I$, which  is $\bL_I^F$-equivariant  (for the  action of  $\bL_I^F$ on
$\bZ^v_\bw$ given in \ref{restante2}).
\end{enumerate}
\end{proposition}
\begin{proof}
\begin{lemma}\label{chgmt modele}
Let $\dot v,\dot w_1,\ldots,\dot w_k$ be representatives
in  $N_\bG(\bT)$ of $v, w_1,\ldots,w_k\in W$ and  let 
$$\CO^v(w_1,\ldots,w_k)=\{(x_1\bB,\ldots,x_k\bB)\in
(\bG/\bB)^k\mid
x_i\inv x_{i+1}\in \bB w_i\bB \text{ and } x_1\in \bB v\bB\};$$ there exist
unique $u\in \bU_v$ and $u_i\in \bU_{w_i}$ such that $x_1\in uv\bB$ and that for all $i$ we have
$x_{i+1}\in u\dot v u_1\dot w_1\ldots u_i\dot w_i\bB$.
This defines a morphism $\CO^v(w_1,\ldots,w_k)\to \bU_v\times(\prod_{i=1}^k\bU_{w_i})$.
\end{lemma}
\begin{proof}
The proof is similar to that of \ref{Uv1Uv2...}.
\end{proof}
Consider the map
$\Psi: (u,g_2\bB,\ldots,g_k\bB)\mapsto(uv\bB,g_2\bB,\ldots,g_k\bB)$
from the variety $$\bZ=\{(u,g_2\bB,\ldots,g_k\bB)\mid
u\in\bU,\,(uv)\inv g_2\in\bB w_1\bB,\,
g_i\inv g_{i+1}\in\bB w_i\bB,\, g_k\inv\lexp F(uv)\in\bB w_k\bB\}$$
to  $\bX^v(\bw)$. As $\bB v\bB=\bU_v v\bB$, it is a fibration whose fibers are
isomorphic  to  $\bU\cap\lexp  v\bU$, an affine space of dimension $l(w_0v)$.
This  fibration  is  $\UF$ (thus $\bU_\bP^F$)-equivariant, for the action by
left multiplication  of  all  components.
Applying \ref{chgmt modele} with $x_i=(uv)\inv g_{i+1}$ we get
$\bZ\simeq\{(u,u_1,\ldots,u_k)\mid  u\in\bU,u_i\in\bU_{w_i}, \;u\inv\lexp Fu\in
\dot vu_1\dot  w_1\ldots u_k\dot  w_k\bB v\inv\}$.  The map  $u=x_\bP x\mapsto
(y=\lexp{x\inv}(x_\bP\inv\lexp  F  x_\bP),x)$,   where  $x_\bP\in\bU_\bP$  and
$x\in\bU_I$,  defines   thus  an   isomorphism  between   $\bZ/\bU_\bP^F$  and
$\bZ^v_\bw$;  \ref{restante1} results immediately  from  this
isomorphism and the isomorphism of cohomology implied by $\Psi$.

Let   us   prove   \ref{restante4}.   It   is  clear   that   the   image   of
the  map   is  in   $\overline\bZ_\bw^v$.  Consider   now  the   fiber  of
$(y_2,x,u'_1,\ldots,u'_k)$.   Let   $y_1$   in  $\bU_\bP\cap\lexp   v\bU$   be
arbitrary; the  formulas $\lexp{\dot  v\inv}y_1\inv u_1\dot  w_1\ldots u_i\dot
w_i\in  u'_1\dot  w_1\ldots  u'_i\dot  w_i\bB$   for  any  $i$  define  unique
$u_i$  (\cf\   \ref{Uv1Uv2...}),   and  the  element   thus  obtained
$(y_1y_2,x,u_1,\ldots,u_k)$ is in the fiber and we get thus all the fiber.

We prove  now \ref{restante2}.  If $v$  is the unique  element of  $W_Iv$ such
that  $\bX^v(\bw)$ is  non  empty, we  will  define an  action  of $\bP^F$  on
$\bZ$ such  that $\Psi$  is $\bP^F$-equivariant.  Under this  assumption, by
\ref{piece  isolee} $\bX^v(\bw)$  is $\bP^F$-stable,  \ie, if  $p\in\bP^F$ and
$(uv\bB,g_2\bB,\ldots,g_k\bB)\in\bX^v(\bw)$  where $u\in\bU_v$,
then   $pu\in\bU_vv\bB$.   The  action   of  $p\in\bP^F$   on
$\bX^v(\bw)$ is thus  given by $$(uv\bB,g_2\bB,\ldots,g_k\bB)\mapsto(\overline
uv\bB,pg_2\bB,\ldots,pg_k\bB)$$  where   $\overline  u\in\bU_v$
is  defined  by  $pu=\overline  ub$   with  $b\in\lexp  v\bB$.  Let  $z\in\bZ$
have    image   $(uv\bB,g_2\bB,\ldots,g_k\bB)$   in    $\bX^v(\bw)$,   \ie,
$z=(uu',g_2\bB,\ldots,g_k\bB)$   where  $u'\in\bU\cap\lexp   v\bU$.  We   want
to define the action of $p$ by the map $$z\mapsto(\overline
u.\proj_{\bU\cap\lexp    v\bU}(bu'),pg_2\bB,\ldots,pg_k\bB)$$    where
$pu=\overline ub$  as above  and where $\proj_{\bU\cap\lexp  v\bU}$ is
the  projection  of $\lexp  v\bB$  on  $\bU\cap\lexp  v\bU$ according  to  the
decomposition $\lexp v\bB=(\bU\cap\lexp v\bU).(\bU^-\cap\lexp v\bU).\bT$. This
map clearly commutes with $\Psi$, but it is not obvious that it defines
an action.  Let $p'\in\bP^F$; let us  write $p'\overline u=\overline{\overline
u}b'$  where $\overline{\overline  u}\in\bU_v$ and  $b'\in\lexp
v\bB$. We  must check  that the  action of $pp'$  is the  composition
of  that of
$p$  and  that  of  $p'$;  this is equivalent to
$\proj_{\bU\cap\lexp v\bU}(b'bu')= \proj_{\bU\cap\lexp
v\bU}(b' \proj_{\bU\cap\lexp  v\bU}(bu'))$, and this last  equality is
easy to check. This action  of $\bP^F$ gives after quotienting $\bZ$ by
$\bU_\bP^F$ an action of $\bL_I^F$ on $\bZ^v_\bw$, for which \ref{restante2}
holds.

Let us now prove \ref{restante3}.    By    \cite[2.5]{LuCox}    the    variety
$\bX_{\bL_I}(\cox)$  has a  single piece  $\bX_{\bL_I}^{w_0^I}(\cox)$. As
$\bc\in\bW$, the model \ref{modeleXvt} gives
$\bX_{\bL_I}^{w_0^I}(\cox)\simeq\{gB_I\in\bL_I/\bB_I\mid
g\inv\lexp F g\in \bB_I\cox\bB_I\text{ and }
g\in\bB_Iw_0^I\bB_I\}$. Defining $u\in \bU_I$ by $g\in uw_0^I\bB_I$
we get $\bX_{\bL_I}^{w_0^I}(\cox)\simeq\{u\in\bU_I\mid \lexp{w_0^I}(u\inv\lexp  Fu)\in \bB_I\cox\bB_I\}$.
As $\bB_I\cox\bB_I\cap\bU^-_I=\prod_{s\in  I}\bU^{-*}_s$   (\cf\  \cite[2.6]{LuCox})
and $\lexp{w_0^I}(\prod_{s\in  I}\bU^*_s)=\prod_{s\in  I}\bU^{-*}_s$,
we get $\bX_{\bL_I}^{w_0^I}(\cox)\simeq\{u\in\bU_I\mid u\inv\lexp  Fu\in\prod_{s\in  I}\bU^*_s\}$,
on which $t\in\bT^F$ acts  by $u\mapsto\lexp tu$ and $u_+\in\bU_I^F$ 
acts by $u\mapsto u_+u$; the action of $u_-\in\bU_I^{-F}$ maps $u$ on
the element $u'$ of $\bU_I$ such that $u_-u\in u'\bB_I^-$: such
an  element exists  by uniqueness  of the  piece $\bX_{\bL_I}^{w_0^I}(\cox)$.
On    the   other    hand   the   projection    on   the
first two   components   of    $\bZ^v(\bw)$   is   surjective   on   couples
$(y,x)\in\bU_\bP\times\bU_I$   which  satisfy   $yx\inv.\lexp  F   x\in\bU\cap
v\bB  w_1\bB\ldots  w_k\bB   v\inv$  and  the  projection   of  $\bU\cap  v\bB
w_1\bB\ldots   w_k\bB   v\inv$   on   $\bU_I$   is   equal   to   $\prod_{s\in
I}\bU_s^*$  by  assumption. Thus
we  see  that  $\pi$ indeed  defines
an  epimorphism  $\bZ^v_\bw\to\bX_{\bL_I}(\cox)$.

It  remains  to  check
that  this  epimorphism  is  $\bL_I^F$-equivariant. For  this,  it  is  enough
to  check   that  the  above  action   and  that  on  the   component  $x$  of
$(y,x,u_1,\ldots,u_k)\in\bZ^v_\bw$, which results itself  of the action on
the  projection on $\bU_I$  of the first component  of  an  element of  $\bZ$,
coincide.  For this  it is  enough  to check  separately that  the actions  of
$\bB_I^F$ and $\bU_I^{-F}$ coincide.

Let $u\in\bU_v$ and $u'\in\bU\cap\lexp v\bU$;
the first component of the image of $z=(uu',g_2\bB,\ldots)\in\bZ$
by the action of $xt\in\bB_I^F$, with $t\in\bT^F$ and $x\in\bU_I^F$,
is $x\lexp{t}u\proj_{\bU\cap\lexp  v\bU}(tu')=x \lexp t u\lexp t u'$,
since $x\in\bU_v$ and $u'\in\bU\cap\lexp v\bU$.
Since $\bU_I\subset\bU_v$, the action of $xt$ on the projection $u_I$ of $uu'$ 
is thus $u_I\mapsto x\lexp t u_I$, which coincides with the action on $\bX_{\bL_I}(\cox)$.

Similarly, the action of  $y\in\bU_I^{-F}$ maps $u_I$ on
$\proj_{\bU_I}(\overline u)$ where $yu\in\overline u\lexp v\bB$ with $\overline u\in\bU_v$.
On the other hand, the image by $y$ of $u_I\in \bX_{\bL_I}^{w_0^I}(\cox)$ is
$u'_I$ such   that   $yu_I\in u'_I\bB_I^-\subset u'_I\lexp  v\bB$.
Let us write $u=u_\bP u_I$ with $u_\bP\in \bU_\bP$;
thus  $yu_\bP u_I=\lexp{y}u_\bP yu_I\in\lexp{y}u_\bP u'_I\lexp  v\bB=
u'_I(\lexp{u_I^{\prime-1}y}u_\bP)\lexp v\bB=u'_I(\lexp{u_I^{\prime-1}y}u_\bP)_-
\lexp v\bB=\lexp{u'_I}(\lexp{u_I^{\prime-1}y}u_\bP)_-u'_I\lexp v\bB$,
where we have denoted by 
$x_-$ the projection on $\bU_v$ of an element $x\in\bU=\bU_v.(\bU\cap\lexp
v\bU)$.  But we have
$\lexp{u'_I}(\lexp{u_I^{\prime-1}y}u_\bP)_-\in\bU_\bP\cap\lexp v\bU^-$ since
$u'_I\in\bU_I\cap\lexp v\bU^-$,  thus $\overline u=
\lexp{u'_I}(\lexp{u_I^{\prime-1}y}u_\bP)_-u'_I$ and its component in $\bU_I$
is $u'_I$ as required.
\end{proof}

The    map     $\pi$    of     \ref{restante3}    is    the     composition
of the    fibration     of    \ref{restante4}     and    the    projection
$$\overline\pi:(y,x,u_1,\ldots,u_k)\mapsto  x$$ which  is an 
epimorphism $\overline\bZ^v_\bw\to\bX_{\bL_I}(\cox)$.
All fibers  of $\bZ^v_\bw\to\overline\bZ^v_\bw$  are affine  spaces of
dimension $l(w_0v)$. It will be easier to compute the fibers of $\overline\pi$
than those of  $\pi$. In the case where we  will apply \ref{variety restante},
we will be in the situation of one of the next two propositions:
\begin{proposition}\label{restante   An}   Assume  under  the  assumptions  of
\ref{restante3} that  the fibers of $\overline\pi$ are affine
spaces  of  dimension  $d$.  Then  for  every  $i$  we  have an isomorphism of
$\bL_I^F\times\genby F$-modules:
$H^i_c(\bX^v(\bw))^{\bU_{\bP_I}^F}\simeq H^{i-2d}_c(\bX_{\bL_I}(\cox))(-d)$.
\end{proposition}
\begin{proof}
\begin{lemma}\label{fibers  affines} With  the same  notations as  above
\ref{B(w)  inter  U-},  let  $\bV$  be  a variety  given  with  an  action  of
$\bP_I^F$,  let $\bV'$  be  a variety  given  with an  action  of $\bL_I^F$  and
let $\pi:\bV/\bU_{\bP_I}^F\to\bV'$  be an $\bL_I^F$-equivariant  epimorphism whose
fibers  are  all  affine  spaces  of  dimension  $d$.  Then  for  all  $i$  we
have  an  isomorphism $H^i_c(\bV)^{\bU_{\bP_I}^F}\simeq  H^{i-2d}_c(\bV')(-d)$  of
$\bL_I^F\times\genby F$-modules.
\end{lemma}
\begin{proof}
This lemma  results from standard  properties of $\ell$-adic cohomology  , see
\eg, \cite[10.10 and 10.12]{DMb}.
\end{proof}
The proposition results  immediately  of  \ref{variety  restante}  and of
\ref{fibers  affines},  taking  in  account  the ``Tate twist'' induced by the
quotient $\bZ^v_\bw\to\overline\bZ^v_\bw$. 
\end{proof}

In the next proposition we assume $\bw\in\bW$ and write $\bX(w)$ for $\bX(\bw)$.
\begin{proposition}\label{restante D4} Assume that in addition to $w$, there
is another element $w'<w$ satisfying the assumptions of \ref{restante3}
with the same $v$; let $\overline\pi':\overline\bZ^v_{w'}\to \bX_{\bL_I}(\cox)$
be the epimorphism analogous to $\overline\pi$ and assume
that the fibers of $\overline\pi'$
are affine lines and that the fibers of
$$\overline\pi\coprod\overline\pi':
\overline\bZ^v_w\coprod\overline\bZ^v_{w'}\to \bX_{\bL_I}(\cox)$$
are  affine planes, the above union being taken in
$(\bU_{\bP_I}\cap\lexp v\bU^-)\times\bU_I\times\bU$.
Then for any $i$ we have an isomorphism of $\bL_I^F\times\genby F$-modules:
$$H^i_c(\bX^v(w))^{\bU_{\bP_I}^F}\simeq H^{i-3}_c(\bX_{\bL_I}(\cox))(-1)\oplus
H^{i-4}_c(\bX_{\bL_I}(\cox))(-2).$$
\end{proposition}
\begin{proof} \ref{restante An} and \ref{fibers affines} give respectively
the isomorphisms of $\bL_I^F\times\genby F$-modules
$H^i_c(\bX^v(w'))^{\bU_{\bP_I}^F}\simeq H^{i-2}_c(\bX_{\bL_I}(\cox))(-1)$ and
$H^i_c(\bX^v(w)\coprod \bX^v(w'))^{\bU_{\bP_I}^F}\simeq
H^{i-4}_c(\bX_{\bL_I}(\cox))(-2)$. The assumption $w'<w$ implies
that $\bX^v(w)$ is open in $\bX^v(w)\coprod \bX^v(w')$.
We deduce a long exact sequence
\begin{align*}
\ldots&\to H^{i-3}_c(\bX_{\bL_I}(\cox))(-1)\to
H^i_c(\bX^v(w))^{\bU_{\bP_I}^F}\to H^{i-4}_c(\bX_{\bL_I}(\cox))(-2)\\
&\to H^{i-2}_c(\bX_{\bL_I}(\cox))(-1)\to\ldots
\end{align*}
We deduce the proposition by observing that a morphism of $\bL_I^F$-modules
from  $H^{i-4}_c(\bX_{\bL_I}(\cox))(-2)$ to $H^{i-2}_c(\bX_{\bL_I}(\cox))(-1)$
must   be   $0$.   Indeed   cohomology   groups   of   different   degrees  of
$\bX_{\bL_I}(\cox)$ are disjoint as $\bL_I^F$-modules.
\end{proof}

\vfill\eject
\section{The $n$-th roots of $\bpi$ in type $A_n$\label{n-regular in An}}
In this section we compute the cohomology groups
$H^i_c(\bX(\bw))$ as $\GF\sdp F$-modules when $\bG$ is a split group of type
$A_n$ ($n\ge 1$) and $\bw$ an $n$-th root of $\bpi$ and we show conjectures
\ref{A} to \ref{E} for this case. 
The Coxeter presentation of $W$ is given by the diagram
$\nnode{s_1}\bar\nnode{s_2}\cdots\nnode{s_n}\kern 4pt$, and we denote
$\bS=\{\bs_1,\ldots,\bs_n\}$ the corresponding generating set of $B$.

Conjecture \ref{A} holds by \ref{c^i et c'^i dans A}, and
conjecture \ref{B} holds as remarked in section \ref{reg in A}.
As noticed in
section \ref{conjectures}, conjecture \ref{conj} follows from a result of
Eilenberg and a recent result of Birman, Gebhardt and Gonzales-Meneses.
However, we will give a simple proof of it in our case.
\begin{proposition}\label{ncycles conj} 
There is a morphism in $\cD^+$ between any two of $n$-th roots of $\bpi$ in 
$B^+$.
\end{proposition}
\begin{proof}
We will show the proposition by showing that any two roots
$\bb$ and $\bb'$ are equivalent by the transitive closure of the
relation on $\BW$ given by $\bx\by\sim\by\bx$.

By \cite[3.12 et A.1.1]{Sydney} the image in $w$, identified with
the symmetric group $\Sgot_{n+1}$,
of any $n$-th root of $\bpi$ is an $n$-cycle. So it is enough to see that if 
$\bb\in\BW$ has length $n+1$ and support $\bS$, and is such that $\beta(\bb)$
is an $n$-cycle then there exists a morphism from $\bb$ to
$\bs_1\ldots \bs_{n-1} \bs_n \bs_n$ in $\cD^+$.

We first show that there exists a morphism in $\cD^+$ from $\bb$
to an element of the form $\bx \bs_n\bs_n \by$.
By assumption,
in the decomposition of $\bb$ as a product of elements of $\bS$, exactly one,
say $\bs_i$ is present twice.  We write
$\bb=\bx\bs_i\by\bs_i\bz$ where
$\bx,\by,\bz\in\BW$ do not have $\bs_i$ in their support. There are three cases:
\begin{itemize}
\item[(a)] The support of $\by$ contains neither $\bs_{i+1}$ nor $\bs_{i-1}$.
Then $\bb=\bx\by\bs_i\bs_i\bz$, so that
$\beta(\bb)=\beta(\bx\by\bz)$ is an element of length $n-1$ of
support $S-\{s_i\}$.
Such an element can be an $n$-cycle only if $i=1$ or $i=n$.
In the latter case $\bb$ is as desired.
If $i=1$, as $\bs_1$ commutes with all elements of $\bS-\{\bs_2\}$,
we see that $\bb$ has the form $\bx\bs_1\bs_1\bs_2\by$ or
$\bx\bs_2\bs_1\bs_1\by$. In both cases 
$\bb$ is equivalent to $\bx\bs_1\bs_2\bs_1\by=\bx\bs_2\bs_1\bs_2\by$ and we
are reduced to case (c) below with $i=2$.
\item[(b)] The support of $\by$ contains $\bs_{i+1}$ and $\bs_{i-1}$. 
Then we use that $\bb$ is equivalent to $\by\bs_i\bz\bx\bs_i$, and
the support of $\bz\bx$ contains neither $\bs_{i+1}$ nor $\bs_{i-1}$. 
We are back to case (a).
\item[(c)] The support of $\by$ contains one of $\bs_{i+1}$ and $\bs_{i-1}$.
Then $\bs_i$ commutes with either $\bx$ or $\bz$ and
if $i=n$ then $\bb$ is equivalent to an 
element of the form we want. Otherwise replacing if needed $\bb$ by
the equivalent element $\by\bs_i\bz\bx\bs_i$, we may assume that $\by$ involves
$\bs_{i+1}$. Then $\bb$ can be written
$\bx\bs_i\bs_{i+1}\bs_i\by'\bz= \bx\bs_{i+1}\bs_i\bs_{i+1}\by'\bz$.
By induction on $i$ we are reduced to $i=n$ whence the result.
\end{itemize}

We prove now that any pair of elements
$\bx \bs_n\bs_n\by$ of length $n+1$ with support $\bS$
are connected by a morphism in $\cD^+$.
We use \cite[chap. V \S 6, lemma 1]{Bbki} which says
``If $X$ is a finite forest and if $x\mapsto g_x$
is a mapping from $X$ to a group $\Gamma$ such that $g_x$ and
$g_y$ commute if and only if  $x$ and $y$ are not connected in $X$, then 
the elements of $\Gamma$  which are the products of all the $g_x$ in some
order are conjugate by cyclic permutation'', where conjugation by cyclic
permutation  is the transitive closure of $g_{x_1}\ldots g_{x_k}\mapsto
g_{x_2}\ldots g_{x_k}g_{x_1}$
(\cite{Bbki} does not state that the conjugation is by cyclic permutation but
it is established in the proof).
We apply this result to the map from the Coxeter diagram
of type $A_n$ which maps the $i$-th vertex to $\bs_i\in B$,
with the exception of the $n$-th vertex which is mapped onto $\bs_n\bs_n$.
This gives the result.
\end{proof}

We denote by $\rho_b^{(n)}$ the unipotent representation of $\GF$ which 
corresponds to the partition $1,\ldots,1,2,b$ of $n+1$. Let $\St^{(n)}$ 
be  the  Steinberg  representation  and  $\Id^{(n)}$  be  the  identity 
representation of $\GF$. We will  deduce conjectures \ref{C} to \ref{E} 
from the following  theorem. In this theorem, we  adopt the conventions 
of \cite[3.3.5]{DMR} to describe the cohomology of a variety $\bX(\bw)$ 
as a  $\GF\sdp F$-module;  we describe the  cohomology as  a 2-variable 
polynomial  with  coefficients  in  the Grothendieck  group  of  $\GF$, 
where  the  degree  in  the  variable  $h$  represents  the  degree  of 
the  cohomology  group,  and  where  the  degree  in  $t$  encodes  the 
eigenvalues  of  $F$:  by  a  theorem of  Lusztig,  given  a  unipotent 
character $\rho$, the eigenvalues of $F$ on the $\rho$-isotypic part of 
a cohomology group $H^j_c(\bX(\bw))$  are of the form $q^i\lambda_\rho$ 
where $\lambda_\rho$  is a  complex number of  module $1$  or $q^{1/2}$ 
which depends only on $\rho$ and neither on $j$ nor on $\bw$. We encode 
such an eigenvalue by $t^i$.                                            

\begin{theorem}\label{principal} Let $\bw\in\BW$ be an $n$-th root of $\bpi$;
then we have as $\GF\genby F$-modules:
$$\sum_i h^i H^i_c(\bX(\bw))=
\St^{(n)}h^{n+1}+\sum_{b=2}^{n-1}\rho_b^{(n)}t^b h^{n+b}+\Id^{(n)}t^{n+1}
h^{2(n+1)}.$$
\end{theorem}
\begin{proof}
We prove the theorem by induction on $n$. If $n=1$ we have $\bw=\bpi$.
Then the only unipotent representations of
$\bG^F$ are $\St^{(1)}$ and $\Id^{(1)}$ and the result is given by
\cite[3.3.14]{DMR} and \cite[3.3.15]{DMR}.
If $n\ge 2$, by \ref{ncycles conj} and
\cite[3.1.6]{DMR}, it is sufficient to prove the result for a fixed root of
$\bpi$.

We choose $\bw=\bs_1\ldots\bs_{n-1}\bs_n\bs_n$.
We shall prove the theorem using the results of section \ref{pieces}.
Let $I=\{s_1,\ldots,s_{n-1}\}$.
\begin{lemma}\label{Y^v_w non vide} The variety $\bX^v(\bw)$
is not empty if and only if $v$ is the longest element in its coset
$v W_I$.
\end{lemma}
\begin{proof} We apply \ref{Ytv}: $\bX^v(\bw)$ is not empty if and only if
$T_v T_\bw| T_v$ is not equal to zero.
In the Hecke algebra $T_\bw=(q-1)T_{s_1\ldots s_n}+
qT_{s_1\ldots s_{n-1}}$, so if $T_v T_\bw| T_v\neq 0$ then
$T_v T_{s_1\ldots s_n}| T_v$
or $T_v T_{s_1\ldots s_{n-1}}| T_v$ is not zero.
By \cite[2.5]{LuCox} or \ref{Yv non vide}, the only $v$ such that
$T_v T_{s_1\ldots s_n}| T_v\ne 0$ is $w_0$.  Let us write
$v=xy$ with $x$ reduced-$I$ and $y\in W_I$; then 
$T_v T_{s_1\ldots s_{n-1}}| T_v\neq 0$ if and only if
$T_y T_{s_1\ldots s_{n-1}}| T_y\neq 0$; by the same result as above,
applied in $W_I$, this coefficient is not zero if and only if
$y=w_0^I$. So $v$ has to be the longest element in its coset, and we have
shown that this is equivalent to the non-vanishing of $T_vT_w|T_v$ except
possibly for  $v=w_0$. In this last case $T_{w_0} T_{s_1\ldots s_n}|
T_{w_0}=(q-1)^n$ and $T_{w_0} T_{s_1\ldots s_{n-1}}| T_{w_0}=(q-1)^{n-1}$
and the sum of these two coefficients is again non zero.
\end{proof}
The reduced-$I$ elements are $s_i  s_{i+1}\ldots  s_n$,  for
$i\le  n$; their number is
$n+1=|W/W_I|$. The elements $v$ of maximal length in their cosets $vW_I$
are then $s_i s_{i+1}\ldots s_n w_0^I= s_1\ldots
s_{i-1} w_0$.  They are in the coset $W_I w_0$,  except $s_1\ldots
s_n w_0=w_0^I$.  Let $\bP$ denote the parabolic subgroup $\bP_I$;
by proposition \ref{piece isolee} $\bX(\bw)$ is the union
of two $\bP^F$-stable pieces: $\bigcup_{v\in W_I}\bX^{vw_0}(\bw)$,
which is an open subvariety as
$\bigcup_{v\in   W_I}\bB   vw_0\bB$ is open in
$\bigcup_{v\in W_I}\bB vw_0\bB \bigcup \bB w_0^I\bB$,
and the closed subvariety $\bX^{w_0^I}(\bw)$.
As $H^i_c(\bX(\bw)/\bU_\bP^F)=\lexp*R^\bG_\bL(H^i_c(\bX(\bw)))$
(\cf\ \eg, \cite[10.10]{DMb}), we get, setting
$\bX_1=(\bigcup_{v\in W_I}\bX^{vw_0}(w))/\bU_\bP^F$ and
$\bX_2=\bX^{w_0^I}(w)/\bU_\bP^F$, the following long exact sequence of
$\LF\genby F$-modules, where $\bL$ denotes $\bL_I$:
\begin{equation}\tag{1}
\ldots\to H^i_c(\bX_1)\to\lexp *R^\bG_{\bL}(H^i_c(\bX(w)))\to
H^i_c(\bX_2)\to H^{i+1}_c(\bX_1)\to\ldots
\end{equation}
We now apply \ref{union of pieces} with
$s=s_1$ (indeed
$\lexp{\bw_0}\bw=\bs_n \bw'$ where
$\bw'=\bs_{n-1}\ldots\bs_2\bs_1\bs_1\in B^+_I$),
whence $$H^i_c(\bX_1)\simeq
H^{i-2}_c(\bX_{\bL}(\lexp{w_0^I}\bw'))(-1)\oplus
H^{i-1}_c(\bX_{\bL}(\lexp{w_0^I}\bw')).$$

The element $\lexp{w_0^I}\bw'=\bs_1\ldots\bs_{n-2}\bs_{n-1}\bs_{n-1}$
is the element of $\bL$ analogous to $\bw$. So by induction on $n$
we get the equality of $\bL^F\genby F$-modules:
\begin{equation}
\sum_i h^i H^i_c(\bX_1)
=(th^2+h)(\St^{(n-1)}h^n+\sum_{b=2}^{n-2}\rho_b^{(n-1)}t^b h^{n+b-1}
+\Id^{(n-1)}t^nh^{2n}).\tag{$1'$}\end{equation}
To compute the $\LF\genby F$-module $\sum_i h^i H^i_c(\bX_2)$,
we apply \ref{variety restante} and \ref{restante An} with $v=w_0^I$,
with $k=n+1$ and with $\bs_1,\bs_2,\ldots,\bs_n,\bs_n$ for $\bw_1,\ldots,
\bw_k$. Let us check the assumptions:
the assumption (iii) on $v$ holds by \ref{Y^v_w non vide};
assumption (iv) holds since
$\lexp{v\inv}\bU_I=\bU_I^-\subset\bU^-$; the other assumptions
to check are
\begin{equation}\tag{2}
{\hbox{proj}}_{\bU_I}(\bU\cap v\bB s_1\bB\ldots\bB
s_{n-1}\bB s_n\bB s_n\bB v\inv)=\prod_{\alpha\in I}\bU^*_\alpha
\end{equation}
and that
the fibers of $\overline\pi$ are affine lines.
We have
\begin{align*}
{\hbox{proj}}_{\bU_I}(\bU\cap v\bB s_1\bB\ldots\bB
s_{n-1}\bB s_n\bB s_n\bB v\inv)&=\\
\lexp v({\hbox{proj}}_{\bU_I^-}(\lexp{v\inv}\bU\cap 
\bB s_1\bB\ldots\bB s_n\bB s_n\bB))&=\\
\lexp v({\hbox{proj}}_{\bU_I^-}((\bU_I^-\cdot\bU_\bP)\cap\bB s_1\bB\ldots\bB
s_{n-1}\bB s_n\bB s_n\bB))&=\\
\lexp v(\bU_I^-\cap \bB s_1\bB\ldots\bB s_{n-1}\bB s_n\bB s_n\bB),
\end{align*}
the last equality as, since $\bU_\bP\subset\bB$ we have
$$(\bU_I^-\cdot\bU_\bP)\cap \bB s_1\bB\ldots\bB s_{n-1}\bB s_n\bB s_n\bB=
(\bU_I^-\cap \bB s_1\bB\ldots\bB s_{n-1}\bB s_n\bB s_n\bB)\cdot\bU_\bP.$$
But we have $\bB s_1\bB\ldots\bB s_{n-1}\bB s_n\bB s_n\bB=
\bB s_1\ldots s_{n-1} s_n\bB \cup \bB s_1\ldots s_{n-1}\bB$.
By \ref{BwB inter U-} we have  $\bB s_1\ldots s_{n-1}\bB\cap\bU^-=
\prod_{i=1}^{i=n-1}\bU^*_{-\alpha_i}\subset \bU_I^-$, and in the same way
(or by \cite[2.2]{LuCox}) we have $\bB   s_1\ldots
s_n\bB\cap  \bU^-=\prod_{i=1}^{i=n}\bU^*_{-\alpha_i}$ which has empty
intersection with $\bU^-_I$. So (2) is proved.

Let us compute the fibers of $\overline\pi$.
We have $\bU_\bP\cap\lexp v\bU^-=\bU_\bP\cap(\bU_I.\bU^-_\bP)=1$, so
$$\displaylines{\overline{\bZ^v_\bw}\simeq\{(x,u_1,\ldots,u_{n+1})
  \mid x\in\bU_I, u_i\in\bU_{s_i} (i=1,\ldots,n),u_{n+1}\in\bU_{s_n},
  \hfill\cr
\hfill x\inv\cdot \lexp Fx\in \dot v u_1\dot s_1\ldots u_{n-1}\dot s_{n-1}
u_n\dot s_n u_{n+1}\dot s_n\bB v\inv\}.\cr}$$
Using the above description of the projection we see that the 
fibers of $\overline\pi$ are those of the map from
$$\displaylines{\{(u_1,\ldots,u_{n+1})\mid
u_i\in\bU_{s_i} (i=1,\ldots,n),u_{n+1}\in\bU_{s_n},\hfill\cr
\hfill
u_1\dot s_1\ldots u_{n-1}\dot s_{n-1}u_n\dot s_n u_{n+1}\dot s_n\bB \in
\bB s_1\ldots s_{n-1}\bB\}\cr}$$ to $\bG/\bB$ given by
$(u_1,\ldots,u_{n+1})\mapsto u_1\dot s_1\ldots u_n\dot s_n u_{n+1}\dot
s_n\bB$. The condition on the $u_i$ implies $u_{n+1}=1$ and that the
image of an $n$-tuple $(u_1,\ldots,u_{n-1},u_n)$ does not depend on $u_n$
and is injective on $(u_1,\ldots,u_{n-1})$. So the fibers are indeed affine
lines.

So we may apply \ref{restante An}: let us write $\gamma_b^{(n)}$ for
the unipotent character corresponding to the partition $1,\ldots,1,b$
of a split group of type $A_n$; multiplying by 
$th^2$ the two variable polynomial which encodes the cohomology of the
Coxeter variety $\bL$, we get that
the cohomology of  $\bX_2$ is given by
$\sum_{b=1}^n t^bh^{n+b}\gamma_b^{(n-1)}$ as a $\bL^F\times\genby F$-module.

We now compute the $t^b$-isotypic part of the exact sequence (1).
For $2\le b\le n-1$ we get the exact sequence
$0\to\rho_b^{(n-1)}+\rho_{b-1}^{(n-1)}\to
(\lexp*R^\bG_{\bL}H^{n+b}_c(\bX(w)))_{t^b}\to\gamma_b^{(n-1)}\to 0$
which gives $\lexp*R^\bG_{\bL}H^{n+b}_c(\bX(w))=\rho_b^{(n-1)}+
\rho_{b-1}^{(n-1)}+\gamma_b^{(n-1)}$. The Littlewood-Richardson formula
allows to compute $\lexp*R^\bG_{\bL}$ of any 
unipotent character. In particular it shows that
the only characters of
$\GF$ whose $\lexp*R^\bG_{\bL}$ contains $\gamma_b^{(n-1)}$ are
$\rho^{(n)}_b$,$\gamma^{(n)}_{b+1}$ and $\gamma^{(n)}_b$; so
$H^{n+b}_c(\bX(w))$ contains one of these three characters. But
$\lexp*R^\bG_{\bL}(\gamma^{(n)}_{b+1})$ and
$\lexp*R^\bG_{\bL}(\gamma^{(n)}_b)$ contain characters different from
$\rho_b^{(n-1)}$, $\rho_{b-1}^{(n-1)}$ and $\gamma_b^{(n-1)}$.
So $H^{n+b}_c(\bX(w))$ contains $\rho^{(n)}_b$. But
$\lexp*R^\bG_{\bL}(\rho^{(n)}_b)=\rho_b^{(n-1)}+
\rho_{b-1}^{(n-1)}+\gamma_b^{(n-1)}$. So 
$H^{n+b}_c(\bX(w))=\rho^{(n)}_b$ as a $\GF$-module.

For $b=0,1,n,n+1$ we get respectively the exact sequences:
$$0\to\St^{(n-1)}\to(\lexp*R^\bG_{\bL}H^{n+1}_c(\bX(w)))_{t^0}\to 0$$
$$0\to(\lexp*R^\bG_{\bL}H^{n+1}_c(\bX(w)))_t\to\St^{(n-1)}\to\St^{(n-1)}\to 
(\lexp*R^\bG_{\bL}H^{n+2}_c(\bX(w)))_t\to 0$$
$$0\to(\lexp*R^\bG_{\bL}H^{2n}_c(\bX(w)))_{t^n}\to\Id^{(n-1)}\to\Id^{(n-1)}\to 
(\lexp*R^\bG_{\bL}H^{2n+1}_c(\bX(w)))_{t^n}\to 0$$
$$0\to\Id^{(n-1)}\to(\lexp*R^\bG_{\bL}H^{2n+2}_c(\bX(w)))_{t^{n+1}}\to 0$$

We know that the only character $\chi$ of $\GF$ such that  $\lexp
*R^\bG_{\bL}\chi$ is $\St^{(n-1)}$-isotypic is $\St^{(n)}$,
and that the only character  $\chi$ such that $\lexp
*R^\bG_{\bL}\chi$ is $\Id^{(n-1)}$-isotypic is $\Id^{(n)}$.
So we see that the
$(H^i(\bX(w)))_{t^b}$ in the above exact sequences are
$\Id^{(n)}$-isotypic or $\St^{(n)}$-isotypic.
The exact sequence for $b=0$ (resp. $b=n+1$)
gives $H^{n+1}(\bX(w))$ (resp. $H^{2n+2}(\bX(w))$).
For $b=1$ or $b=n$, applying
propositions \cite[3.3.14]{DMR} and \cite[3.3.15]{DMR}  we see that the
$(H^i(\bX(w)))_{t^b}$ in the above exact sequences must be zero
and the arrows $\St^{(n-1)}\to\St^{(n-1)}$ and $\Id^{(n-1)}\to\Id^{(n-1)}$
must be isomorphisms.
This completes the proof of the theorem.
\end{proof}

Let  us explain  now  why theorem  \ref{principal} implies  conjectures 
\ref{C}  to  \ref{E}; conjecture  \ref{E}  is  immediate. Let  us  show 
\ref{C}.  By section  \ref{reg  in A},  the  centralizer $C_B(\bw)$  is 
cyclic,  generated by  $\bw$, and  $C_W(w)=G(1,1,n)$. The  endomorphism 
$D_\bw$   acts  as   $F$  on   $\bX(\bw)$.  Thus   the  value   of  the 
eigenvalues  of $F$  given in  theorem \ref{principal}  shows that  the 
representation $\bw\mapsto  D_\bw$ of $B(1,1,n)$  on $\End_\GF(\oplus_i 
H^i_c(\bX(\bw)))$  factors  through  an  $n$-cyclotomic  Hecke  algebra 
$\CH(w)$  with  parameters  $(1,x^2,x^3,\ldots,x^{n-1},x^n)$.  To  show 
\ref{D},   it  remains   to   see  that   the  virtual   representation 
$\sum_i  (-1)^i H^i_c(\bX(\bw))$  of  $\CH(w)$  is special.  Proceeding 
as   in   section   \ref{coxeter},   it  is   enough   to   show   that 
$|\bX(\bw)^{F^i}|=0$  for $i=1,\ldots,n-1$.  But  this  is exactly  the 
statement \cite[5.2]{Sydney}.                                           

\vfill\eject
\section{Conjecture \ref C in type $A$}

We consider a group  $\bG$ of type $A_{n-1}$ and a  group $\bG'$ of type
$A_n$. We keep the notations of section \ref{reg in A}: $W$ (resp. $W'$)
is the  Weyl group of $\bG$  (resp. $\bG'$) and we  consider $\bw=\bc^r$
(resp. $\bw'=\bc^{\prime r}$),  a $d$-th root of $\bpi$ in the  braid group $B$
(resp. $B'$) of $W$ (resp. $W'$). We have $dr=n$ with $d\ge 2$.

In   section  \ref{reg   in  A}   we   had  two   incarnations  of   the
braid   group   $B(d,1,r)$:   one   as   $C_B(\bw)$,   with   generators
$\bt,\bs_1,\ldots,\bs_{r-1}$ and  one as $C_{B'}(\bw')$  with generators
$\bt',\bs_1,\ldots,\bs_{r-1}$; in both cases these generators correspond
to braid reflections of $B(d,1,r)$.

The    group     $G(d,1,r)$    has    two    orbits     of    reflecting
hyperplanes   corresponding   to   reflections    of   order   $d$   and
$2$   respectively,   so   that   for   a   choice   of   indeterminates
$\bu=(u_{\bt,0},\ldots,u_{\bt,d-1};u_{\bs_1,0},u_{\bs_1,1})$ the generic
Hecke algebra $\CH_\bu$ of $G(d,1,r)$ it is the quotient of $\Qlbar[\bu]
B(d,1,r)$  by  the relations  $(\bt-u_{\bt,0})\ldots(\bt-u_{\bt,d-1})=0$
and  $(\bs_1-u_{\bs_1,0})(\bs_1-u_{\bs_1,1})=0$;  we  will   write  the
first  relation   as  $(\bt'-u_{\bt,0})\ldots(\bt'-u_{\bt,d-1})=0$  when
considering the other incarnation of $B(d,1,r)$.

Next theorem proves that \ref C holds.

\begin{theorem} 

\begin{itemize}
\item
The map $\bx\mapsto D_\bx$ from $C_B(\bw)$ to $\End_\GF(\oplus_i
H^i_c(\bX(\bw)))$ factors through the specialization $x\mapsto q$ of a
$d$-cyclotomic Hecke algebra $\CH$ for $B(d,1,r)$ with parameters
$(1,x,x^2,\ldots,x^{d-1};x^d,-1)$.
\item
The map $\bx\mapsto D_\bx$ from $C_{B'}(\bw')$ to 
$\End_\GF(\oplus_i H^i_c(\bX(\bw')))$ factors through the specialization
$x\mapsto q$ of a
$d$-cyclotomic Hecke algebra $\CH'$ for $B(d,1,r)$ with parameters
$(1,x^2,x^3,\ldots,x^{d-1},x^{d+1};x^d,-1)$.
\end{itemize}
\end{theorem}

With the notation of section \ref{reg in A}, we have to show that the
operators induced on the cohomology by 
$D_{\bs_i}$, $D_\bt$ and $D_{\bt'}$ satisfy the expected polynomial
relations. 
The end of this section is devoted to the proof of this theorem.

The next lemma will allow us to compute by induction the eigenvalues of
$D_{\bs_i}$, using \cite[5.2.9]{DMR}.
\begin{lemma}\label{alpha(bw)}
For $i=1,\ldots,r-1$ let
$\bI_i=\{\bsigma_i,\bsigma_{i+r},\ldots,\bsigma_{i+(d-1)r}\}$
as in \ref{s_i ok}; then we have 
$\alpha_{\bI_i}(\bw)=\alpha_{\bI_i}(\bw')=\bsigma_i^2$.
\end{lemma}
\begin{proof}
As the elements of $\bI_i$ commute pairwise and as, by
\ref{facts} (iii), $\bsigma_i$ is the only divisor of $\bw$ in $\bI_i$,
we have $\alpha_{\bI_i}(\bw)=\bsigma_i^k$ for some $k$. Now
$$\bsigma_i^k\preccurlyeq\bw\Leftrightarrow\bsigma_i^{-k}\bw\in\BW\Leftrightarrow 
\bc^{i-1}\bsigma_1^{-k}\bc^{r-i+1}\in\BW.$$
But $\bsigma_1$ does not divide $\bc^{i-1}$
on the right for $i<n$ by the ``right-side version'' of \ref{facts} (iv)
so, by \ref{xy-1z},  $\bsigma_i^k\preccurlyeq\bw$ if and only if
$\bsigma_1^k\preccurlyeq\bc^{r-i+1}$. To prove the assertion about $\alpha_{\bI_i}(\bw)$
it is sufficient to show that
$\bsigma_1^2\preccurlyeq\bc^2$ and $\bsigma_1^3\not\preccurlyeq\bc^r$.
The former statement follows from \ref{facts} (ii) which gives
$\bsigma_1^{-2}\bc^2=\bsigma_1\inv\bc^2\bsigma_{n-1}\inv\in\BW$. To get the latter
we write $\bsigma_1^{-3}\bc^r=(\bsigma_1^{-2}\bc^2)\bsigma_{n-1}\inv\bc^{r-2}$
and, as $\bsigma_{n-1}\not\preccurlyeq\bc^{r-2}$, we have to see  by \ref{xy-1z} that
$\bsigma_1^{-2}\bc^2\not\succcurlyeq\bsigma_{n-1}$, which is equivalent to
$\bsigma_1^{-2}\bc^2\bsigma_{n-1}=\bsigma_1^{-3}\bc^2\not\in\BW$.
But $\bsigma_1^3\preccurlyeq\bc^2$ is impossible by \cite[4.8]{michel} as $\nu(\bsigma_1^3)=3$
and $\nu(\bc^2)\le 2$ (recall that $\nu(\bb)=\inf\{k\in \bbN|\bb\preccurlyeq\bw_0^k\}$
for $\bb\in B^+$).

The proof of the assertion about $\alpha_{\bI_i}(\bw')$ follows the same lines, using
\ref{facts} (iv') instead of \ref{facts} (iv) and \ref{facts} (ii') instead of \ref{facts} (ii).
At the end we have to see that
$\bsigma_1^3\not\preccurlyeq\bc^{\prime2}$. But by
\ref{facts} (v), we have $\bc^{\prime 2}=\bc^2\bsigma_{n-1}\bsigma_n$, which can be
written $\bc^{\prime
2}=(\bc\bsigma_1\ldots\bsigma_{n-1})(\bsigma_n\bsigma_{n-1}\bsigma_n)$.
The two factors are in $\bW$, so that $\nu(\bc^{\prime2})=2$, and we conclude as in the
$\bw$ case.
\end{proof}
\begin{proposition}\label{Tsi} For $i\in\{1,\ldots,r-1\}$, the image
$T_{\bs_i}$ of $D_{\bs_i}$ in either
$\End_\GF(\oplus_jH^j_c(\bX(\bw)))$ or
$\End_\GF(\oplus_jH^j_c(\bX(\bw')))$ satisfies
$(T_{\bs_i}+1)(T_{\bs_i}-q^d)=0$.
\end{proposition}
\begin{proof}
We prove the statement for $\bw$, the proof for $\bw'$ being exactly the same.
By \cite[5.2.9]{DMR}, which can be applied by \ref{s_i ok} (i),
the eigenvalues of
$T_{\bs_i}$ are equal to those of $D_{\bs_i}$ on
$\oplus_jH^j_c(\bX_{\bL_{\bI_i}}(\alpha_{\bI_i}(\bw),\omega_{\bI_i}(\bw)F))$.
By the above lemma we have $\alpha_{\bI_i}(\bw)=\bsigma_i^2$. The group
$\bL_{\bI_i}$ has 
type $A_1^d$ where the components are permuted cyclically by
$\omega_{\bI_i}(\bw)F$. If $\bs$ denotes the positive generator of the braid group of type 
$A_1$, through the isomorphism with $A_1^d$ we have
$\bsigma_i^2\leftrightarrow(\bs^2,1,\ldots,1)$, 
$\bs_i\leftrightarrow(\bs,\bs,\ldots,\bs)$ and $\omega_{\bI_i}(\bw)F$
corresponds to $(x_1,\ldots,x_d)\mapsto
(\lexp Fx_2,\ldots,\lexp Fx_d,\lexp Fx_1)$.
The variety we have to study is
$\bX_{\bL_{\bI_i}}(\bs^2,1,\ldots,1)$,
which can be identified with the three term sequences of $d$-tuples of Borel subgroups of
$\bL_{I_i}$ of the form
$$(\bB_1,\ldots,\bB_d)\xrightarrow
{(s,1,\ldots,1)}(\bB'_1,\ldots,\bB'_d)\xrightarrow{(s,1,\ldots,1)}
(\lexp F\bB_2,\ldots,\lexp F\bB_d,\lexp F\bB_1),$$
where for two Borel subgroups $\bB$ and $\bB'$, we write
$\bB\xrightarrow v\bB'$ to say that $(\bB,\bB')\in\CO(v)$ (we say that $\bB$
and $\bB'$ are in relative position $v$).
The conditions on the relative positions imply
$\bB'_2=\bB_2$,\dots, $\bB'_d=\bB_d$ and $\lexp F\bB_3=\bB_2$,\dots,
$\lexp F\bB_d=\bB_{d-1},\lexp F\bB_1=\bB_d$; so that
$\bX_{\bL_{\bI_i}}(\bs^2,1,\ldots,1)$ identifies with the three term sequences
$\bB_1\xrightarrow s\bB'_1 \xrightarrow s\lexp F\bB_2=\lexp{F^d}\bB_1$, \ie, to the
variety $\bX(\bs^2)$ of a group of type
$A_1$ with Frobenius endomorphism $F^d$.
Let us put $\bs_{(i)}=(\bs,1,\ldots,1,\bs,1,\ldots,1)$ where the second
$\bs$ is at the place $i$ with $i>1$.
In the same way as above we identify the variety
$\bX_{\bL_{\bI_i}}(\bs_{(i)})$ to the variety $\bX(\bs^2)$ by 
identifying a sequence $(\bB_1,\ldots,\bB_d)$ such that $(\bB_1,\ldots,\bB_d)\xrightarrow{s_{(i)}}
(\lexp F\bB_2,\ldots,\lexp F\bB_{d-1},\lexp F\bB_d)$ 
with the three term sequence $\bB_1\xrightarrow s\lexp{F^{i-1}}\bB_i\xrightarrow s\lexp{F^d}\bB_1$.
We can decompose the morphism $D_{\bs,\ldots,\bs}$ as
$D_{\bs^{(2)}}\circ\ldots\circ D_{\bs^{(d)}}\circ D_{\bs^{(1)}}$, where
$\bs^{(i)}=(1,\ldots,1,\bs,1,\ldots,1)$, with $\bs$ at the $i$th place.
Then $D_{\bs^{(1)}}$ sends $\bX(\bs^2,1,\ldots,1)$ to $\bX(\bs_{(d)})$
and $D_{\bs^{(i)}}$ sends $\bX(\bs_{(i)})$ to $\bX(\bs_{(i-1)})$ for
$i>1$.
With the above identifications, one checks that
$D_{\bs^{(1)}}$ sends
$\bB_1\xrightarrow s\bB'_1\xrightarrow s\lexp{F^d}\bB_1$ to 
$\bB'_1\xrightarrow s\lexp{F^d}\bB_1\xrightarrow s\lexp{F^d}\bB'_1$
and that $D_{\bs^{(i)}}$ sends 
$\bB_1\xrightarrow s\lexp{F^{i-1}}\bB_i\xrightarrow s\lexp{F^d}\bB_1$
to $\bB_1\xrightarrow s\lexp{F^{i-2}}(\lexp F\bB_i)\xrightarrow s\lexp{F^d}\bB_1$
for $i>1$. So that 
$D_{(\bs,\ldots,\bs)}$ identifies with the operator $D_\bs$ on the variety
$\bX(\bs^2)$ of a group of type $A_1$ with Frobenius endomorphism
$F^d$. This is a particular case of \cite[5.3.4]{DMR}, where it is proved that the operator $T_\bs$
induced by $D_\bs$ on $H_c^*(\bX(\bs^2)$ satisfies $(T_\bs-q^d)(T_\bs+1)=0$, whence the proposition.
\end{proof}

We will now prove that the operators induced on the cohomology of $\bX(\bw)$
and $\bX(\bw')$ respectively by $D_\bt$ and $D_{\bt'}$
satisfy the claimed polynomial relation.

\begin{lemma}\label{ywy-1}
Let $\bI=\{\bsigma_r,\ldots,\bsigma_{r+d-2}\}$ and
$\bI'=\{\bsigma_r,\ldots,\bsigma_{r+d-1}\}$, as in \ref{t in End_D};
we have
\begin{itemize}
\item[(i)] $\alpha_\bI(\by\bw\by\inv)=\by\bt\by\inv$,
\item[(i')]$\alpha_{\bI'}(\by'\bw'\by^{\prime -1})=\by'\bt'\by^{\prime -1}$.
\item[(ii)] $\omega_\bI(\by\bw\by\inv)$ 
commutes with $\bsigma_{r+i}$ for $0\le i\le d-2$,
\item[(ii')] $\omega_{\bI'}(\by'\bw'\by^{\prime-1})$
commutes with $\bsigma_{r+i}$ for $0\leq i\leq d-1$.
\end{itemize}
\end{lemma}
\begin{proof} Let us prove (i).
As $\by\bt\by\inv=\bsigma_{r,r+d-2}\in B_I^+$ and
$\by\bw\by\inv=\by\bt\by\inv\bc^{r-1}
\bx_d\ldots\bx_1$
by the remark which follows the proof of \ref{y in B<-w},
it is enough to see that $\alpha_\bI(\bc^{r-1}\bx_d\ldots\bx_1)=1$,
\ie, that for $i\in 0,\ldots,d-2$
we have $\bsigma_{r+i}\not\preccurlyeq\bc^{r-1}\bx_d \ldots\bx_1$.
By \ref{xy-1z}, this amounts to prove that for
$i\in 1,\ldots,d-1$  we have
$\bsigma_i\not\preccurlyeq\bx_d \ldots\bx_1$.
But $\bsigma_i$ commutes with  $\bx_j$ for $j> i+1$, so that
$\bsigma_i\inv\bx_d\ldots\bx_1= \bx_d\ldots\bx_{i+2}\bsigma_i\inv\bx_{i+1}\ldots\bx_1$.
As $\bsigma_i$ is not in the support of
$\bx_d\ldots\bx_{i+2}$, we have $\bx_d\ldots\bx_{i+2}\not\succcurlyeq\bsigma_i$, so by
\ref{xy-1z} it is enough to see that
$\bsigma_i\not\preccurlyeq\bx_{i+1}\ldots\bx_1$.
But $\bsigma_i\inv\bx_{i+1}\bx_i=
\bx_{i+1}\bx_i\bsigma_{i+r-1}\inv$: indeed we can write this equality as
$\bsigma_{i+1,i+r-1}\bsigma_{i,i+r-1}=\bsigma_{i,i+r-1}\bsigma_{i,i+r-2}$,
using \ref{facts}(i) in the parabolic subgroup generated by
$\bsigma_i,\ldots,\bsigma_{i+r-1}$. So, as
$\bsigma_{i+r-1}\not\preccurlyeq\bx_{i-1}\ldots\bx_1$
because $\bsigma_{i+r-1}$ is not in the support of this element,
again by \ref{xy-1z} it is sufficient to see that
$\bx_{i+1}\bx_i\not\succcurlyeq\bsigma_{i+r-1}$.  But this is the exact analogue in
the parabolic subgroup generated by
$\bsigma_i,\ldots,\bsigma_{i+r-1}$ of
what we have proved at the end of lemma \ref{alpha(bw)},
as $\bsigma_i^{-2}\bsigma_{i,i+r-1}^2=\bx_{i+1}\bx_i$.

The proof of (i') is along the same lines:
we start with
$\by'\bw'\by^{\prime -1}=\by'\bt'\by^{\prime
-1}\bc_{rd}^{r-1}\bx_{d+1}\ldots\bx_1$
as noticed in the remark following \ref{y' in B<-w'}.
We have to see that
$\alpha_{\bI'}(\bc_{rd}^{r-1}\bx_{d+1}\ldots\bx_1)=1$,
\ie, that
$\bsigma_{r+i} \not \preccurlyeq \bc_{rd}^{r-1}\bx_{d+1}\ldots\bx_1$ for $i\in 0,\ldots,d-1$.
By \ref{xy-1z},  this amounts to prove that for
$1,\ldots,d$ we have $\bsigma_i\not\preccurlyeq\bx_{d+1}\ldots\bx_1$.
But $\bsigma_i$ commutes with $\bx_j$ for $j> i+1$,
and as in the proof of (i) we are reduced to prove that
$\bsigma_i\not\preccurlyeq\bx_{i+1}\ldots\bx_1$. Then we finish exactly as in (i).

Let us prove (ii).
We have to see that $\bsigma_{r+i}\bc^{r-1}\bx_d\ldots\bx_1=
\bc^{r-1}\bx_d\ldots\bx_1\bsigma_{r+i}$. This can be written
$$\bc^{r-1}\bx_d\ldots \bx_{i+3}\bsigma_{i+1}\bx_{i+2}\bx_{i+1}\ldots\bx_1=
\bc^{r-1}\bx_d\ldots\bx_{i+2}\bx_{i+1}\bsigma_{r+i}\bx_i\ldots\bx_1.$$
So we have to prove that
$\bsigma_{i+1}\bx_{i+2}\bx_{i+1}=\bx_{i+2}\bx_{i+1}\bsigma_{r+i}$,
\ie,
$\bsigma_{i+1,i+r}\bsigma_{i+1,i+r-1}=\bsigma_{i+2,i+r}\bsigma_{i+1,r+i}$.
This last equality is a consequence of \ref{facts} (i) applied in the parabolic 
subgroup generated by $\{\bsigma_{i+1},\ldots,\bsigma_{i+r}\}$.

The proof of (ii') is exactly the same, replacing
$\bI$ by $\bI'$, $\bc$ by $\bc_{rd}$ and $\bx_d\ldots\bx_1$ by $\bx_{d+1}\ldots\bx_1$.
\end{proof}
\begin{corollary}\label{vp de t} 
\begin{itemize}
\item The image of $T_\bt$ of $D_\bt$ in $\End_\GF(\oplus_iH^i_c(\bX(\bw)))$
satisfies $(T_\bt-1)(T_\bt-q)\ldots(T_\bt-q^{d-1})=0$.
\item The image of $T_{\bt'}$ of $D_{\bt'}$ in 
$\End_\GF(\oplus_iH^i_c(\bX(\bw')))$ satisfies
$(T_{\bt'}-1)(T_{\bt'}-q^2)(T_{\bt'}-q^3)\ldots(T_{\bt'}-q^{d-1})(T_{\bt'}-q^{d+1})=0$.
\end{itemize}
\end{corollary}
\begin{proof}
We prove first the result for $\bt$.
Using conjugation by $D_\by$ we see that
it is equivalent to prove that the image of $D_{\by\bt\by\inv}$ in
$\End_\GF(\oplus_iH^i_c(\bX(\by\bw\by\inv)))$ satisfies the same relation.
By \cite[5.2.9]{DMR} and \ref{ywy-1} the eigenvalues of
this operator are the same as those of
$D_{\by\bt\by\inv}$ on $\oplus_iH^i_c(\bX_{\bL_I}
(\alpha_\bI(\by\bw\by\inv),\omega_\bI(\by\bw\by\inv)F))=\oplus_iH^i_c(\bX_{\bL_I}(\bsigma_{r,r+d-2}))$.
The element $\bsigma_{r,r+d-2}$ is the lift of a  Coxeter element in the braid group of type $\bI$.
The group $\bL_I$ has type 
split  $A_{d-1}$, and $D_{\by\bt\by\inv}$ acts as $F$ on
$\bX_{\bL_I}(\bsigma_{r,r+d-2})$. The eigenvalues of $F$ in the case of a Coxeter element
are given in \cite{LuCox} and are the eigenvalues which we want for $T_\bt$.

We follow the same lines for proving the assertion on $\bt'$.
We are reduced to compute the eigenvalues of
$D_{\by'\bt'\by^{\prime -1}}$ on $\oplus_iH^i_c(\bX_{\bL_{I'}}
(\alpha_{\bI'}(\by'\bw'\by^{\prime -1}),
\omega_{\bI'}(\by'\bw'\by^{\prime -1})F))=
\oplus_iH^i_c(\bX_{\bL_{I'}}(\bsigma_{r,r+d-1}\bsigma_{r+d-1}))$.
The element $\bsigma_{r,r+d-1}\bsigma_{r+d-1}$ relatively to
the group $\bL_{I'}$ is an element as studied in section
\ref{n-regular in An}.
The group $\bL_{I'}$ is split of type
$A_d$ and $D_{\by'\bt'\by^{\prime -1}}$ acts as $F$ on
$\bX_{\bL_{I'}}(\bsigma_{r,r+d-1}\bsigma_{r+d-1})$.
We end the proof by noticing, as in the end of section \ref{n-regular in An},
that the eigenvalues of $F$ given in
\ref{principal} give the eigenvalues we want for $T_{\bt'}$.
\end{proof}

\vfill\eject
\section{Conjecture \ref C in type $B$}

In this section we keep the notation of \ref{reg in B};
we prove conjecture \ref C in a group $\bG$
of type $B_n$ for a $d$-regular element with $d$ even.
So we have $\bw=\bc^r$ with $dr=n$ and $d$ even.
\begin{theorem}
The map $\bx\mapsto D_\bx$ from $C_B(\bw)$ to $\End_\GF(\oplus_i
H^i_c(\bX(\bw)))$ factors through the specialization $x\mapsto q$ of a
$d$-cyclotomic Hecke algebra $\CH$ for $B(d,1,r)$ with parameters
$(1,x,x^2,\ldots,x^{d/2-1},x^{d/2},-x,-x^2,\ldots,-x^{d/2-1};x^{d/2},-1)$.
\end{theorem}

With the notation of section \ref{reg in B}, we have to show that the
operators induced on the cohomology by 
$D_{\bs_i}$ and $D_\bt$ satisfy the expected polynomial relations. 

The end of this section is devoted to the proof of this theorem.
The next lemma is the analogue of \ref{alpha(bw)}.

\begin{lemma}\label{alpha(bw)B}
For $i=2,\ldots,r$ let
$\bI_i=\{\bsigma_i,\bsigma_{i+r},\ldots,\bsigma_{i+(d/2-1)r}\}$ as in
\ref{s_i okB}; then we have $\alpha_{\bI_i}(\bw)=\bsigma_i$.
\end{lemma}
\begin{proof}

As the elements of $\bI_i$ commute pairwise and as, by
\ref{factsBn} (iv), $\bsigma_i$ is the only divisor of $\bw$ in $\bI_i$,
we have $\alpha_{\bI_i}(\bw)=\bsigma_i^k$ for some $k$.
But $\bc$ is a good root of $\bpi$, so for any $r\leq n$ we have
$\nu(c^r)=1$. As $\nu(\bsigma_i^2)=2$ we cannot have $\bsigma_i^2\preccurlyeq
\bw$, whence the result.
\end{proof}

\begin{corollary}\label{TsiBw} For $i\in\{1,\ldots,r-1\}$ 
the image  $T_{\bs_i}$ of $D_{\bs_i}$ in
$\End_\GF(\oplus_jH^j_c(\bX(\bw)))$ satisfies
$(T_{\bs_i}+1)(T_{\bs_i}-q^{d/2})=0$.
\end{corollary}
\begin{proof}
The proof is similar to that of \ref{Tsi}, except that here
$D_{\bs_i}$ eventually identifies to
the operator $D_\bs$ on the variety
$\bX(\bs)$ for a group of type $A_1$ with Frobenius endomorphism
$F^{d/2}$, whence the result.
\end{proof}

\begin{lemma}\label{ywy-1B}
Let $\bI=\{\bsigma_1,\ldots,\bsigma_{d/2}\}$ as at the end of section
\ref{reg in B}; we have
\begin{enumerate}
\item $\alpha_\bI(\by\bw\by\inv)=\by\bt\by\inv$.
\item $\omega_\bI(\by\bw\by\inv)$ commutes with  $\bsigma_i$
for $1\le i\le d/2$.
\end{enumerate}
\end{lemma}
\begin{proof}
Let us prove (i). As $\by\bt\by\inv=\bsigma_{1,d/2}$ and
$\by\bw\by\inv=\by\bt\by\inv\prod_{i=d/2}^1\bx_i\bc^{r-1}$
by the remark which follows the proof of \ref{y in B<-wBn}, it is enough to
see that $\bsigma_i\not\preccurlyeq\prod_{i=d/2}^1\bx_i\bc^{r-1}$
pour $i=0,\ldots,d/2$.
We use the following lemma which is an immediate consequence of
\ref{xy-1z}:
\begin{lemma}\label{chaine} Assume that $\ba,\bb\in B^+$, that
$\bsigma_i\preccurlyeq\ba\bb$, that
$\bsigma_i\ba=\ba\bsigma_j$ and that $\bsigma_i\not\preccurlyeq\ba$;
then $\bsigma_j\preccurlyeq\bb$.
\end{lemma}
By this lemma, for $i>1$ we have
$\bsigma_i\preccurlyeq\prod_{j=d/2}^1\bx_j\bc^{r-1}\Leftrightarrow
\bsigma_i\preccurlyeq\prod_{j=i}^1\bx_j\bc^{r-1}\Leftrightarrow
\bsigma_{i+r-1}\preccurlyeq\prod_{j=i-2}^1\bx_j\bc^{r-1}\Leftrightarrow
\bsigma_{i+r-1}\preccurlyeq\bc^{r-1}$,
the second equivalence as, by the proof of \ref{ywy-1},
we have $\bsigma_i\inv\bx_i\bx_{i-1}= \bx_i\bx_{i-1}\bsigma_{i+r-1}\inv$.
But $\bsigma_{i+r-1}\preccurlyeq\bc^{r-1}$ is false by \ref{factsBn}(iv).

For $i=1$ we get
$\bsigma_1\preccurlyeq\prod_{j=d/2}^1\bx_j\bc^{r-1}\Leftrightarrow
\bsigma_1\preccurlyeq\bx_1\bc^{r-1}\Leftrightarrow
\bsigma_1^2\preccurlyeq\bsigma_1\ldots\bsigma_r\bc^{r-1}\Rightarrow
\bsigma_1^2\preccurlyeq\bc^r$, the last equivalence as
$\bsigma_{r+1}\ldots\bsigma_n \bc^{r-1}=\bc^{r-1}\bsigma_2\ldots\bsigma_
{n-r+1}$. But $\bsigma_1^2\preccurlyeq\bc^r$ is impossible as
$\nu(\bc^r)=1$.

Let us prove (ii). By (i) we have
$\omega_\bI(\by\bw\by\inv)=\prod_{i=d/2}^1\bx_i\bc^{r-1}$.
For $i>1$, we have
$\prod_{j=d/2}^1\bx_j\bc^{r-1}\bsigma_i=
\prod_{j=d/2}^1\bx_j\bsigma_{i+r-1}\bc^{r-1}=
\prod_{j=d/2}^{i+1}\bx_i\bx_{i-1}\bsigma_{i+r-1}\prod_{j=i-2}^1\bx_j\bc^{r-1}=
\prod_{j=d/2}^{i+1}\bsigma_i\bx_i\bx_{i-1}\prod_{j=i-2}^1\bx_j\bc^{r-1}=
\bsigma_i\prod_{j=d/2}^1\bx_j\bc^{r-1}$.

For $i=1$, we prove by induction on  $i$ that
$\bsigma_2\ldots\bsigma_{i+1}\bc^i\bsigma_1=\bsigma_1\ldots\bsigma_{i+1}\bc^i$;
this equality for $i=r-1$ gives
$\bx_1\bc^{r-1}\bsigma_1=\bsigma_1\bx_1\bc^{r-1}$; then we use that
$\bsigma_1$ commutes with $\bx_i$ for $i>1$.
\end{proof}

As for \ref{vp de t}, we deduce the following corollary which ends the proof
of the theorem.
\begin{corollary}\label{vp de tBn} 
Let $T_\bt$ be the image of $D_\bt$ in $\End_\GF(\oplus_iH^i_c(\bX(\bw)))$;
then we have $(T_\bt-1)(T_\bt-q)\ldots(T_\bt-q^{d/2})(T_\bt+q)(T_\bt+q^2)\ldots
(T_\bt+q^{d/2-1})=0$.
\end{corollary}
\begin{proof}
Using conjugation by $D_\by$, we see that it is equivalent to prove that
the image of $D_{\by\bt\by\inv}$ in
$\End_\GF(\oplus_iH^i_c(\bX(\by\bw\by\inv)))$ satisfies the same relation.
By \cite[5.2.9]{DMR} and \ref{ywy-1B} the eigenvalues of this
operator are the same as those of
$D_{\by\bt\by\inv}$ on $\oplus_iH^i_c(\bX_{\bL_I}
(\alpha_\bI(\by\bw\by\inv),\omega_\bI(\by\bw\by\inv)F))=
\oplus_iH^i_c(\bX_{\bL_I}(\prod_{i=1}^{d/2}\bsigma_i))$.
The element $\prod_{i=1}^{d/2}\bsigma_i$ is the lift of a Coxeter element
in the braid group of type $\bI$. The group $\bL_I$ has type
$B_{d/2}$ and $D_{\by\bt\by\inv}$ acts as $F$ on
$\bX_{\bL_I}(\prod_{i=1}^{d/2}\bsigma_i)$.
The eigenvalues of $F$ in the case of a Coxeter element
are given in \cite{LuCox} and are the eigenvalues which we want for $T_\bt$.
\end{proof}

\vfill\eject
\section{The $4$-th roots of $\bpi$ in $D_4$}

We  use  the  same numbering for the reflections of $W=W(D_4)$ as in Bourbaki:
we     will     study     the    $4$-th    root    of    $\bpi$    given    by
$\bw=\bs_2\bs_3\bs_1\bs_3\bs_4\bs_3$.
$$
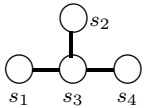
\def\node{{\kern -0.4pt\bigcirc\kern -1pt}}
\def\nnode#1{{\kern -0.6pt\mathop\bigcirc\limits_{#1}\kern -1pt}}
\def\bar#1pt{{\vrule width#1pt height3pt depth-2pt}}
\def\vertbar#1#2{\rlap{$\nnode{#1}$}
                 \rlap{\kern4pt\vrule width1pt height17.3pt depth-7.3pt}
		 \raise19.4pt\hbox{$\node$\rlap{$\kern 1pt\scriptstyle#2$}}}
\nnode{s_1}\bar10pt\vertbar{s_3}{s_2}\bar10pt\nnode{s_4}
$$

The centralizer of $w$ is the complex reflection group $G(4,2,2)$. We 
show  conjectures \ref{A} to \ref{E} for $w$.

\begin{proposition}
The   element    $\bw$   verifies   conjectures   \ref{A}    and   \ref{conj}.
In   particular,   the   map   which  sends   the   standard   generators   of
the braid group $B(4,2,2)$ to $\bb_1=(\bs_1\bs_2)^{\bs_3}$, $\bb_2=\bs_1\bs_4$
and $\bb_3=(\bs_2\bs_4)^{\bs_3\bs_4}$ is an isomorphism  $B(4,2,2)\simeq
C_B(\bw)$.
\end{proposition}
\begin{proof}
A  4-th  root  of $\pi$ is an element of length 6 whose image is a $4$-regular
element.  As  all  twelve  $4$-regular  elements of $W$ are of length $6$, all
$4$-th  roots of $\bpi$
are in $\bW$. There are thus $12$ good roots of $\bpi$, and a direct check
shows  that  there  are  morphisms  in  $\cD^+$  between any two of them, thus
proving  \ref{conj}  in this case. We may show \ref{A} using the programs
written   by   Franco  and Gonzalez-Meneses
to   compute   centralizers   in   braid   groups   .
They     give     that     $\bb_1$, $\bb_2$  and $\bb_3$ generate $C_B(\bw)$.
However,          one          sees          that         the         elements
$\bw\bb_1=\bs_1\bs_2\bs_3\bs_1\bs_2\bs_4\bs_1\bs_3$,        $\bb_2$        and
$\bw\bb_3=\bs_2\bs_3\bs_1\bs_2\bs_3\bs_4\bs_3\bs_3$           are           in
$\End_{\cD^+}(\bw)$  by using the decompositions
$\bw\bb_1=(\bs_1\bs_2\bs_3\bs_1)(\bs_2\bs_4)(\bs_1\bs_3)$,
$\bb_2=(\bs_1)(\bs_4)$ and $\bw\bb_3=(\bs_2\bs_3\bs_1)   (\bs_4)  (\bs_2\bs_3\bs_4)  (\bs_3)$.  We  get  thus  that
$\bb_1, \bb_2$ and $\bb_3$ are in $\End_\cD(\bw)$ whence \ref{A}.

Let us  show that the  elements $\bb_1, \bb_2,  \bb_3$ are the  image of
braid reflections around hyperplanes of $C_W(w)$. We start with $\bb_2$;
let  $V$ be  the  reflection representation  of $W$,  and  let $V_i$  be
the  eigenspace  of  $w$ in  $V$  for  the  eigenvalue  $i$. This  is  a
2-dimensional space, and the fixed points  of $s_1s_4$ on this space are
of  dimension $1$,  and  of the  form $H\cap  V_i$  for some  reflecting
hyperplane of  $W$. We use now  the constructions at the  end of section
\ref{regular}:  the  group  $B_H$  is  the  parabolic  subgroup  of  $B$
generated by  $\bs_1$ and $\bs_4$, and  the element $\bs_H$ is  a square
root of the element $\bpi$ of this parabolic subgroup, so it is equal to
$\bs_1\bs_4$ thus this element is indeed the image in $B$ of a the braid
reflection $\bt_H$  of $B(w)$. To handle  the other elements we  use the
triality  automorphism of  $\bG$.  This automorphism  corresponds to  an
automorphism  of  $V$  which  induces  on  $B$  the  automorphism  given
by  $\bs_2\mapsto\bs_1\mapsto\bs_4\mapsto\bs_2$  and  the  corresponding
automorphism on $W$. Let $\bpsi$ be the automorphism of $B$ which is the
triality followed by  conjugation by $\bs_2\bs_3$ and let  $\psi$ be the
corresponding automorphism  of $V$.  Then $\bpsi$  fixes $\bw$  and does
the permutation $\bb_1\mapsto\bb_2\mapsto\bb_3\mapsto\bb_1$. Thus $\psi$
induces an  automorphism of $V_i$,  and the braid reflections  which are
images of $\bt_H$ by the powers of this automorphism have images $\bb_1$
and $\bb_3$.

The  fact the images  of these braid reflections generate  $C_B(\bw)$  imply
that  they  generate  $B(w)$; thus conjecture \ref{B} holds.  They
satisfy the  defining braid relations for  the braid group of  $G(4,2,2)$, that
is $\bb_1\bb_2\bb_3=\bb_2\bb_3\bb_1=\bb_3\bb_1\bb_2$  (all three  products are
easily checked to be equal to $\bw$).
\end{proof}

The next proposition shows conjecture \ref{C}.
\def\Hcyc{{\CH(w)_q}}
\begin{proposition}
The   map   $\bb_i\mapsto   D_{\bb_i}$  factors   through   a   representation
$\Hcyc\to\oplus_i H^i_c(\bX(w))$ of the  specialization $x\mapsto q$ of
a $4$-cyclotomic Hecke algebra $\CH(w)$ for $G(4,2,2)$ with parameters $1$ and $x^2$.
\end{proposition}
We note  that the  above algebra  is indeed a  4-cyclotomic algebra  since the
parameters specialize to $1$ and  $-1$ by the specialization $x\mapsto
i$.
\begin{proof}
It   is   enough   to   prove   the   quadratic   relations
$(D_{\bb_i}-1)(D_{\bb_i}-q^2)=0$.

For    $\bb_2=\bs_1\bs_4$,   we   use   proposition   \cite[5.2.8]{DMR}   with
$\bx=\by=\bs_1$   and   $I=\{s_1,s_4\}$.  The  subgroup  $\bL_I$  is  of  type
$A_1\times  A_1$; The element $z=s_1\inv w$ normalises $\bL_I$, exchanging the
two  $A_1$  components.  The  variety  $\bX_{\bL_I}(\bs_1,\dz F)$ is mapped by
$D_{\bs_1}$, which  acts  as  $\dz  F$,  to  $\bX_{\bL_I}(\bs_4,\dz F)$. Then
$D_{\bs_4}$  maps  $\bX_{\bL_I}(\bs_4,\dz  F)$  to $\bX_{\bL_I}(\bs_1,\dz F)$,
and    $D_{\bs_4}\circ    D_{\bs_1}$   induces   $F^2$   on
$\bX_{\bL_I}(\bs_1,\dz  F)$;  by  \eg, \cite[3.3.16]{DMR} this variety has two
non zero cohomology  groups with compact suport: $H^1_c$ on which $F^2$ acts as $1$ and
$H^2_c$  on  which  $F^2$  acts  as  $q^2$.  Thus  $D_{\bb_2}$  satisfies  the
quadratic  relation  $(D_{\bb_2}-1)(D_{\bb_2}-q^2)=0$  on  the  cohomology  of
$\bX_{\bL_I}(\bs_1,\dz  F)$,  and thus by \cite[2.3.13]{DMR} and the K\"unneth
formula, it also satisfies the same relation on the cohomology of $\bX(w)$.

To   show   that   $D_{\bb_3}$   satisfies   the   same   quadratic   relation
we   use    \cite[3.1.8]{DMR}   applied    to   the    triality   automorphism
of  $\bG$.  
The triality   maps   $\bw$
to   $\bs_1\bs_3\bs_4\bs_3\bs_2\bs_3$   and  $D_{\bb_2}$   to
$D_{\bs_2\bs_4}$  acting  on  $\bX(\bs_1\bs_3\bs_4\bs_3\bs_2\bs_3)$  which  by
the  same  argument  as   above,  with  $\bx=\by=\bs_4$  and  $I=\{s_2,s_4\}$,
satisfies   the    expected   quadratic    relation.   The    conjugation   by
$\bs_2\bs_3\in\Hom_\cD(\bw,\bs_1\bs_3\bs_4\bs_3\bs_2\bs_3)$ composed with the
triality is equal to $\bpsi$; and the conjugation by $\bs_2\bs_3$ maps 
back $D_{\bs_2\bs_4}$  to $D_{\bb_3}$ which thus  satisfies the quadratic
relation.

One      proceeds      similarly      for     $D_{\bb_1}$,      using      the
square of the  triality  automorphism  which  maps    $D_{\bs_1\bs_4}$    to
$D_{\bs_2\bs_1}$     acting   on   $\bX(\bs_4\bs_3\bs_2\bs_3\bs_1\bs_3)$.
One has that $\bpsi^2$ is the square of triality followed by conjugation
by $\bs_2\bs_3\bs_1\bs_3$, and
 conjugating     by      $\bs_2\bs_3\bs_1\bs_3\in
\Hom_\cD(\bw,\bs_4\bs_3\bs_2\bs_3\bs_1\bs_3)$ maps $D_{\bs_2\bs_1}$ to
$D_{\bb_1}$.
\end{proof}

If we denote by $T_1, T_2, T_3$ the images in $\Hcyc$ of $\bb_1, \bb_2, \bb_3$,
we  thus get  a representation  $\rho$  of $\Hcyc$  on $\oplus  H^i_c(\bX(w))$
which  maps $T_i$  to $D_{\bb_i}$.  We note  that triality  commutes with  the
automorphism of $\Hcyc$ given by $T_1\mapsto T_2\mapsto T_3\mapsto T_1$.

The next proposition and \ref{D4} show that conjecture \ref{D} holds.
\begin{proposition}
The  representation  induced  by  $\rho$ to  $\sum  (-1)^i  H^i_c(\bX(w))$  is
special.
\end{proposition}
\begin{proof}
Let $\Phi=T_1T_2T_3$. It  follows from the computation of  the generic degrees
and the character  table of $\Hcyc$ in \cite[5.A]{BrMa}  that a representation
is  special  if   (and  only  if)  the  elements   $\{  (T_i)_{i=1,2,3},  (T_i
T_j)_{i=1,2,3, j=1+(i\bmod 3)}, (\Phi^i)_{i=1,2,3}\}$ have trace zero on it.
The quadratic relations show that $\CB=\{(T_i)_{i=1,2,3}, (T_i\Phi)_{i=1,2,3},
(\Phi^i)_{i=1,2,3}\}$  generates  the  same   subspace  of  $\Hcyc$  as  these
elements. We shall  show that the trace  of any element in  $\CB$ vanishes on
$\sum (-1)^i H^i_c(\bX(w))$. It is enough to check this on just one element in
each orbit of  $\CB$ under triality, \eg, for $T_2,  T_1\Phi, \Phi^i$. We have
$\rho(\Phi)=D_\bw=F$. Thus $D_\bw$  as well as its square and  cube verify the
trace formula, and by \cite[5.2.3]{DMR} they all have trace zero on $\sum (-1)^i
H^i_c(\bX(\bw))$. Thus  it only remains  to compute the traces  of $D_{\bb_2}$
and $D_{\bb_1\bw}$.

To    compute    the    trace    of   $D_{\bb_2}=D_{\bs_1\bs_4}$,   we   apply
\cite[5.2.10]{DMR}  with  $g=1$,  $I=\{s_1,s_4\}$  and $\bx=\bs_1\bs_4$, which
gives  that  the trace of $D_{\bs_1\bs_4}$ on $\sum (-1)^i H^i_c(\bX(\bw))$ is
the  value  at  $1$  of  the  class  function  $$R_{\bL_I^{\dz  F}}^\GF  \left
(l\mapsto\TrH{lD_{\bs_1\bs_4}}{\bX_{\bL_I}(\bs_1,\dz  F)} \right ).$$ By
a  general  result  on  Lusztig  induction  (see \eg, \cite[12.17]{DMb}), this
value  is  a  multiple  of  $$\TrH{D_{\bs_1\bs_4}}{\bX_{\bL_I}(\bs_1,\dz
F)}.$$  But this last trace is zero since $D_{\bs_1\bs_4}$ acts as $F^2$
on $\bX_{\bL_I}(\bs_1,\dz F)$, thus has no fixed points.

By  the  computation  done  when checking the quadratic relations, we see that
the  trace  of $D_{\bb_1\bw}$ on $\sum (-1)^i H^i_c(\bX(\bw))$ is equal to the
trace         of        $D_{\bs_2\bs_1}F$        on        $\sum        (-1)^i
H^i_c(\bX(\bs_3\bs_2\bs_3\bs_1\bs_3\bs_4))$.    By    \cite[2.3.13]{DMR}   and
\cite[5.2.9]{DMR}  applied  with  $I=\{s_1,s_2\}$  this  trace  is,  if we set
$z=s_3s_2s_1s_3s_4$,  equal  to  $$|\bL_I^{\dz F}|\inv\sum_{l\in\bL_I^{\dz F}}
\TrH{lD_{\bs_2\bs_1}F}{\bX_{\bL_I}(\bs_2,\dz                               F)}
\TrH{lF}{\tilde\bX_{(I)}(\dz)}$$  (the action of $F$ decomposes as a product
since  $\bU_I$  is $F$-stable). By the Lefschetz trace formula one has
$\TrH{lF}{\tilde\bX_{(I)}(\dz)}=|\tilde\bX_{(I)}(\dz)^{lF}|= 
|\{g\bU_I\mid  g\inv\lexp  Fg\in(\bU_I\dz\bU_I)\cap(l\bU_I)\}|$. But this last
intersection  is  empty  since  $z\notin  W_I$ so that $\bP_I$ does not intersect
$\bU_I\dz\bU_I$. Thus the trace vanishes.
\end{proof}

We  shall  now  prove  conjecture  \ref{E}  by  giving  a  full description of
$H^*_c(\bX(\bw))$  as  a  $\GF\sdp  F$-module.  To  give  the result, we first
introduce  a  notation  for  the  irreducible  characters of the Weyl group of
$D_4$:  they  are  parametrized  by  the pairs  of  partitions of total sum $4$,
except  that  a  pair  of  equal  partitions corresponds to 2 characters. We
shall          thus          denote          the         characters         by
$$1^2+,1^2-,1.1^3,1^4,1^2.2,1.21,21^2,2+,2-,2^2,1.3,31,4$$  where  a missing
``.''  means  that one of the partitions is empty. If $\lambda$ is a parameter
as  above,  we  shall  denote by $\gamma_\lambda$  the  corresponding  unipotent
character  of  the principal series, except that we denote respectively
by $\St$ and $\Id$ the characters $\gamma_{1^4}$ and $\gamma_4$.
There is one more unipotent
character  of  $D_4$,  a  cuspidal  one,  that  we  shall  denote  by  $\theta$.
We  will  use  the  same  convention  as  in  \ref{principal}  to  describe  the
cohomology  by  a  two  variable   polynomial.  With  these  notation,  we  have
\begin{theorem}\label{D4} The cohomology of $\bX(\bw)$ is given by
$$h^6\St+t^2   h^7(\gamma_{1^2+}+\gamma_{1^2-}+\gamma_{21^2})+2t^3h^7\theta+
2t^3 h^8\gamma_{1.21}+t^4 h^9(\gamma_{2+}+\gamma_{2-}+\gamma_{31})+t^6
h^{12}\Id.$$ \end{theorem} \begin{proof}  We will use the  parabolic subgroup of
type $A_3$  generated by $I=\{s_1,s_3,s_4\}$. Let  $w'=s_3s_1s_3s_4s_3$. We will
need the  value of  the sets  defined in  \ref{calculE}: they  are $E_{W_I}(w')=
\{1,s_1,s_3,s_4\}$ and  $E_W(w)= E_{W_I}(w')\cup\{s_2s_3\}$. The value  for $w'$
is obtained by a direct computation in  the Hecke algebra of $W_I$ and the value
for  $w=s_2 w'$  comes from  lemma \ref{calculE2}  whose assumptions  are easily
checked.

We   first    compute   the   cohomology   of    $\bX(\bw^{\prime   \bw_0^I})=
\bX(\bs_3\bs_4\bs_3\bs_1\bs_3)$  in the  Levi  subgroup $\bL_I$  corresponding
to  $I$.  We use  \cite[3.2.10]{DMR};  with  the  notation  of the  proof  of
\ref{principal} ($\rho_2$  corresponds to  the partition $2,2$,  $\gamma_2$ to
$1,1,2$, and  $\gamma_3$ to $1,3$) we  know by \cite{LuCox} the  cohomology of
the Coxeter variety $\bs_4\bs_3\bs_1$:
$$\sum_i h^i\cdot H^i_c(\bX(\bs_4\bs_3\bs_1))=h^3\St+h^4
t\gamma_2+h^5 t^2\gamma_3+h^6 t^3\Id$$
and since $\bs_3\bs_4\bs_3\bs_1$ is a $3$-rd root
of  $\bpi_I$  we know its  cohomology  by  \ref{principal}
$$\sum_i h^i\cdot H^i_c(\bX(\bs_3\bs_4\bs_3\bs_1))=h^4\St+t^2 h^5\rho_2+h^8
t^4\Id.$$
The  exact sequence  given  by \cite[3.2.10]{DMR}  then completely  determines
$H^i_c(\bX(\bs_3\bs_4\bs_3\bs_1\bs_3))$,   except  for   the  $\Id$-isotypic
and   $\St$-isotypic   parts;  we   get   these   by  \cite[3.3.14]{DMR}   and
\cite[3.3.15]{DMR},  and we  obtain $$\sum_i  h^i\cdot H^i_c(\bX(w'))=h^5\St+t^2
h^6(\gamma_2+\rho_2)+ t^3h^7(\gamma_3+\rho_2)+t^5h^{10}\Id.$$

We now determine the principal series  part of the cohomology of $\bX(\bw)$ by
the same  method as  in the  proof of \ref{principal}.  The value  of $E_W(w)$
shows  that  $$\bX(w)=\bX^{w_0}(w)\coprod\bX^{w_0s_1}(w)\coprod\bX^{w_0s_3}(w)
\coprod\bX^{w_0s_4}(w)\coprod\bX^{w_0s_2s_3}(w).$$      We      may      apply
\ref{union         of         pieces}         with         $s=s_2$.         If
we set $$\bY=\left(\bX^{w_0}(w)\coprod\bX^{w_0s_1}(w)\coprod\bX^{w_0s_3}(w)
\coprod\bX^{w_0s_4}(w)\right)/\bU_{\bP_I}^F,$$   we    have   thus   an    equality   of
$\bL_I^F$-modules:
\begin{multline*}\sum_i  h^i\cdot  H^i_c(\bY)=(th^2+h)\times
(\sum_i            h^i\cdot           H^i_c(\bX(w_0^Iw'w_0^I)))\\
=h^6\St+h^7(t\St+t^2(\gamma_2+\rho_2))+h^8t^3(\gamma_2+\gamma_3+2\rho_2)+
h^9t^4(\gamma_3+\rho_2)+t^5h^{11}\Id+t^6 h^{12}\Id.\end{multline*}
To study the  remaining piece $\bX^{w_0 s_2 s_3}(w)$ we  use \ref{restante D4}
with $v=w_0 s_2 s_3$ and taking for $w'$ the element $w''=s_2s_3s_1s_4s_3$. We
show that $w$ and $w''$ satisfy the assumptions of \ref{restante D4}.

We  first  check  the  assumptions  of \ref{restante2}
for  $w$  and $w''$. 
We check on $E_W(w)$ that $v$   is  the  only  element  of   $W_Iv$  such  that
$\bX^v(w)\neq\emptyset$.
On the  other hand, one  gets by  \ref{calculE2} whose assumptions  are easily
checked  that  $E_W(w'')=E_{W_I}(s_3s_1s_4s_3)\cup \{s_2s_3,  s_2s_3s_1s_4\}$,
where  $E_{W_I}(s_3s_1s_4s_3)=\{1,s_1,s_3,  s_4, s_1s_4,  s_3s_1s_4\}$, which
results from a direct computation in the  Hecke algebra of $W_I$. Thus $v$ is
indeed the only element of $W_Iv$ such that $\bX^v(w'')\neq\emptyset$.

We   check   now   the    assumption   \ref{restante3}  for
$w$    and    $w''$.
If  $\alpha_1,  \alpha_2,  \alpha_3, \alpha_4$  are  the
simple  roots of  $D_4$, we  have $v\inv(\{\alpha_1,\alpha_3,\alpha_4\})=
\{-\alpha_1-\alpha_3,-\alpha_2,-\alpha_3-\alpha_4\}$,    thus    $\lexp{v\inv}
\bU_I\subset\bU^-$ .
The    projection    onto    $\lexp{v\inv}\bU_I$    of
$\lexp{v\inv}\bU\cap(\bB  w\bB\coprod \bB  w''\bB)$  is the  same  as that  of
$\lexp{v\inv}\bU\cap\bU^-\cap(\bB w\bB\coprod \bB w''\bB)$; indeed each double
coset is  invariant by right multiplication  by $\lexp{v\inv}\bU\cap\bU$, thus
its  intersection with  $\lexp{v\inv}\bU$  is the  product  of its  respective
intersections with $\lexp{v\inv}\bU\cap\bU^-$ and $\lexp{v\inv}\bU\cap\bU$. As
$\lexp{v\inv}\bU_I\subset\bU^-$,  the part  in $\lexp{v\inv}\bU\cap\bU$  has a
trivial projection.

Let us compute
$\lexp{v\inv}\bU\cap\bU^-\cap(\bB w\bB\coprod \bB w''\bB)$.
By \ref{B(w) inter U-} we have $\bB w\bB\cap\bU^-=\bU^*_{-\alpha_2}.
(\bB_I w'\bB_I\cap\bU^-_I)$ and
$\bB w''\bB\cap\bU^-=\bU^*_{-\alpha_2}.(\bB_I s_3s_1s_4s_3\bB_I\cap\bU^-_I)$,
whence, from the decomposition
$\lexp{v\inv}\bU\cap\bU^-=(\lexp{v\inv}\bU\cap\bU^-_{\bP_I}).
(\lexp{v\inv}\bU\cap\bU^-_I)$ we get
$$ \lexp{v\inv}\bU\cap\bU^-\cap(\bB w\bB\coprod \bB w''\bB)=
\bU^*_{-\alpha_2}.(\lexp{v\inv}\bU\cap\bU^-_I\cap
(\bB_I w'\bB_I\coprod\bB_I s_3s_1s_4s_3\bB_I)). $$
We   have  $\lexp{v\inv}\bU\cap\bU^-_I=\prod_{\alpha\in\Phi_I^--\{-\alpha_3\}}
\bU_\alpha$. In order to make explicit computations, we may replace $\bL_I$ by
the  group $\GL_4$  as the  Borel subgroup  varieties, as  well as  the groups
$\bU_I$ and  $\bU_\alpha$ depend only  on the  isogeny type together  with the
Frobenius action on the  root system. The variety $\lexp{v\inv}\bU\cap\bU^-_I$
is thus isomorphic to the variety of matrices
$\left(\begin{smallmatrix}1&0&0&0\\ *&1&0&0\\ *&0&1&0\\ *&*&*&1\\
\end{smallmatrix}\right)$.

We   now  determine   the  variety   of  matrices   in  $\GL_4$   representing
elements   of   $\lexp{v\inv}\bU\cap\bU^-_I\cap   (\bB_I   w'\bB_I\coprod\bB_I
s_3s_1s_4s_3\bB_I)$. For this we use the following:
\begin{lemma}\label{BwB  dans  GLn}  Let  $\bB$   be  the  Borel  subgroup  of
$\GL_n$  of upper  triangular matrices;  let $W$  be the  Weyl group  relative
to  the  torus  of  diagonal  matrices;  let  $(w_{ij})$  be  the  matrix  for
$w\in  W$;  then   $(a_{ij})\in\GL_n$  is  in  $\bB  w\bB$  if   and  only  if
the   ranks  of   the  submatrices   $(a_{ij})_{i=k,\dots,n,j=1,\dots,l}$  and
$(w_{ij})_{i=k,\dots,n,j=1,\dots,l}$ coincide for all $k$ and $l$.
\end{lemma}
\begin{proof}
The  condition   on  ranks   is  invariant   by  left   multiplication  (resp.
right  multiplication)  by   $\bB$  since  this  replaces   each  line  (resp.
column)  by  a   non-zero  multiple  of  itself  plus   a  linear  combination
of  the  following  lines  (resp.  columns). This  condition  defines  thus  a
union  of  double  $\bB$-cosets.  It  remains to  see  that  elements  of  $W$
are  determined  by  the  rank  conditions:  but  indeed,  in  line  $k$,  the
position  $l$  of  the  non-zero  coefficient is  the  smallest  integer  such
that  the  rank  of   the  matrices  $(w_{ij})_{i=k,\dots,n,j=1,\dots,l}$  and
$(w_{ij})_{i=k+1,\dots,n,j=1,\dots,l}$ differ. \end{proof}

We thus  obtain, characterizing  $w'$ and  $s_3s_1s_4s_3$ by  rank conditions,
that
$\lexp{v\inv}\bU\cap\bU^-_I\cap
(\bB_I w'\bB_I\coprod\bB_I s_3s_1s_4s_3\bB_I)$
is the variety of matrices
$\left(\begin{smallmatrix}1&0&0&0\\ d&1&0&0\\\alpha&0&1&0\\ 0&\beta
&f&1\\\end{smallmatrix}\right)$, with $\alpha$, $\beta$, $d$ and $f$ in
$\Fqbar$, such that
$\left|\begin{smallmatrix}\alpha&0\\0&\beta\\\end{smallmatrix}\right|\neq 0$;
the open subset corresponding to $\bB_I w'\bB_I$ is given by the condition
$\left|\begin{smallmatrix}d&1&0\\\alpha&0&1\\
0&\beta&f\\\end{smallmatrix}\right|\neq 0$.
These matrices may be written as
$\left(\begin{smallmatrix}1&0&0&0\\ 0&1&0&0\\\alpha&0&1&0\\
0&\beta&0&1\\\end{smallmatrix}\right)
\left(\begin{smallmatrix}1&0&0&0\\ d&1&0&0\\0&0&1&0\\
-d\beta&0&f&1\\\end{smallmatrix}\right)\in\lexp{v\inv}\bU_I\lexp{v\inv}\bU_{\bP_I}$.
The projection on $\lexp{v\inv}\bU_I$ of
$\lexp{v\inv}\bU\cap\bU^-\cap\bB w\bB$ as well as that of
$\lexp{v\inv}\bU\cap\bU^-\cap\bB w''\bB$ is thus
$\bU_{-\alpha_2}^*\bU_{-\alpha_1-\alpha_3}^*\bU_{-\alpha_3-\alpha_4}^*$.
The assumption \ref{restante3} thus holds for $w$ and $w''$.

To  check  the   assumptions  of  \ref{restante  D4},  we   must  compute  the
fibers of  the maps $\overline\pi'$ and  $\overline\pi\coprod\overline\pi'$ of
\ref{restante D4}.

The above  computations show  that for  $y\in\lexp{v\inv}\bU_{\bP_I}\cap\bU^-$ and
$x\in\lexp{v\inv}\bU_I$,  we  have  $y.x\inv\lexp F  x\in\bB  w\bB\coprod  \bB
w''\bB$ if and only if $x\inv.\lexp F x\in
\prod_{\alpha\in\lexp{v\inv}\bU_I}\bU_\alpha^*$ and $\lexp{\lexp F x\inv x}y$
is in $\bU_{-\alpha_1}\bU_{-\alpha_4}\bU_{-\alpha_1-\alpha_4}$ and are such that
the latter element corresponds to the matrix
$\left(\begin{smallmatrix}1&0&0&0\\ d&1&0&0\\ 0&0&1&0\\ -d\beta&0
&f&1\\\end{smallmatrix}\right)$
where the projection of $x\inv\lexp F x$ on $\bU_I^-$ is given by the matrix
$\left(\begin{smallmatrix}1&0&0&0\\ 0&1&0&0\\\alpha&0&1&0\\
0&\beta&0&1\\\end{smallmatrix}\right)$.
The   closed   subset   corresponding   to    $\bB   w''\bB$   is   given   by
$d\beta+f\alpha=0$.   We   see    thus   that   the   fibers    of   the   map
$\overline\pi\coprod\overline\pi'$  of  \ref{restante  D4}  are  2-dimensional
affine spaces corresponding  to the $d$ and $f$ coordinates  of the matrix for
$y$ and  that the fibers of  the map $\overline\pi'$ are  1-dimensional affine
subspaces corresponding to the equation $d\beta+f\alpha=0$.

The assumptions of  \ref{restante D4} thus hold, and we  get the cohomology of
$\bX^v(w)/\bU_{\bP_I}^F$ as an $\LF\genby F$-module  by multiplying by $th^3+t^2h^4$ the
two-variable polynomial  encoding the  cohomology of  the Coxeter  variety for
$\bL_I$. We get
$$\sum_ih^i\cdot H^i_c(\bX^v(w))^{\bU_{\bP_I}^F}=h^6t\St+h^7t^2(\gamma_2+\St)+
h^8t^3(\gamma_2+\gamma_3)+h^9t^4(\gamma_3+\Id)+h^{10}t^5\Id.$$

The long exact sequence of $\bL_I^F\genby F$-modules
$$\ldots\to H^i_c(\bY)\to H^i_c(\bX(w))^{\bU_{\bP_I}^F}\to
H^i_c(\bX^v(w))^{\bU_{\bP_I}^F}\to H^{i+1}_c(\bY)\to\ldots$$
gives
$$\displaylines{\hfil
0\to\St\to H^6_c(\bX(w))^{\bU_{\bP_I}^F}\to t\St\to t^2(\gamma_2+\rho_2)+ t\St
\to H^7_c(\bX(w))^{\bU_{\bP_I}^F}\to t^2(\gamma_2+\St)\to 0\hfil\cr
\hfil H^8_c(\bX(w))^{\bU_{\bP_I}^F}=t^3(2\gamma_2+2\rho_2+2\gamma_3)\hfil\cr
\hfil H^9_c(\bX(w))^{\bU_{\bP_I}^F}=t^4(\rho_2+2\gamma_3+\Id)\hfil\cr
\hfil 0\to H^{10}_c(\bX(w))^{\bU_{\bP_I}^F}\to t^5\Id\to t^5\Id\to
H^{11}_c(\bX(w))^{\bU_{\bP_I}^F}\to 0\hfil\cr
\hfil H^{12}_c(\bX(w))^{\bU_{\bP_I}^F}=t^6\Id\hfil\cr}$$
To obtain the non cuspidal part of $H^*_c(\bX(\bw))$, we first
use \cite[3.3.14]{DMR} and \cite[3.3.15]{DMR} which give the $\Id$ and $\St$
isotypic parts. We then consider for each $b$
the $t^b$-isotypic part of the above exact sequences
arguing as the proof  of \ref{principal} and using the following table
which describes $\scal{\Ind_{A_3} ^{D_4}\chi}\psi{D_4}$:
\halign{\tabskip 1em$#$&&$#$\cr
&1^2+&1^2-&1.1^3&\St&1^2.2&1.21&21^2&2+&2-&2^2&1.3&31&\Id\cr
\Id & 0& 0& 0& 0& 0& 0& 0& 0& 1& 0& 1& 0& 1 \cr
\gamma_2 & 1& 0& 1& 0& 1& 1& 1& 0& 0& 0& 0& 0& 0\cr
\rho_2 & 0& 1& 0& 0& 0& 1& 0& 0& 1& 1& 0& 0& 0 \cr
\gamma_3 & 0& 0& 0& 0& 1& 1& 0& 1& 0& 0& 1& 1& 0 \cr
\St& 0& 1& 1& 1& 0& 0& 0& 0& 0& 0& 0& 0& 0 \cr}

We get:
\begin{multline*}\sum_i h^i\cdot H^i_c(\bX(w))_{\text{non cuspidal}}=\\
h^6\St+t^2 h^7(\gamma_{1^2+}+\gamma_{1^2-}+\gamma_{21^2})+
2t^3 h^8\gamma_{1.21}+
t^4 h^9(\gamma_{2+}+\gamma_{2-}+\gamma_{31})+t^6 h^{12}\Id.\\
\end{multline*}

To study  the $\theta$-part  of the  cohomology of $\bX(w)$,  we will  use the
variety $\overline\bX(w)$. One may  check that all Kazhdan-Lusztig polynomials
$P_{y,w}$  for $y\le  w$ are  1, thus  $\overline\bX(w)$ is  rationally smooth
(\cf\ \cite[3.2.5]{DMR}),  which by  \cite[3.3.8 (ii)]{DMR} allows  to compute
$\sum_i h^i\cdot (H^i_c(\overline\bX(w)))_\theta=t^3h^6\theta$.

By  \cite[3.1.3]{DMR}, if  $y$  lies in  a proper  parabolic
subgroup,  $H^i_c(\bX(y))$  cannot  have  a  cuspidal  part.  The  only  $y<w$
for  which  this  does  not  hold  are  the elements of
$$\CC=\{s_2s_3s_1s_4s_3,  s_2s_1s_3s_4s_3,
s_2s_1s_3s_1s_4, s_2s_3s_1s_4, s_2s_1s_4s_3, s_2s_1s_3s_4\}.$$ Thus
$\overline\bX(w)=\bX(w)\coprod\bX\coprod\bY$ where
$\bX=\coprod_{v\in\CC}\bX(v)$    and   where    $\bY$    is    a   union    of
Deligne-Lusztig  varieties,   closed  in   $\overline\bX(w)$  and   such  that
$H^i_c(\bY)_\theta=0$ for  any $i$. The  long exact sequence  corresponding to
$\overline\bX(w)=\bY\coprod(\overline\bX(w)-\bY)$ shows thus  that for any $i$
we have $H^i_c(\overline\bX(w))_\theta=H^i_c(\bX(w)\cup\bX)_\theta$.

The      varieties      $\bX(s_2s_3s_1s_4)$,      $\bX(s_2s_1s_4s_3)$      and
$\bX(s_2s_1s_3s_4)$,       corresponding       to      Coxeter       elements,
satisfy       (\cf\       \eg,      \cite{LuCox})       $\sum_i       h^i\cdot
(H^i_c)_\theta=t^2h^4$,   and    they   are   connected    components   of
their    union,    thus     $\sum_i    h^i\cdot    H^i_c(\bX(s_2s_3s_1s_4)\cup
\bX(s_2s_1s_4s_3)\cup\bX(s_2s_1s_3s_4))_\theta=3t^2h^4$.    The   elements
$s_2s_3s_1s_4s_3$,  $s_2s_1s_3s_4s_3$ and  $s_2s_1s_3s_1s_4$ are  conjugate by
cyclic  permutation  respectively  to  $s_2s_3s_1s_4s_2$,  $s_4s_2s_1s_3s_4$
and  $s_1s_2s_3s_4s_1$  which,  by  \cite[3.1.6]{DMR}  and
\cite[3.2.10]{DMR}, allows to compute the cuspidal part of their  cohomology 
For each of them
we   get   $\sum_i    h^i\cdot   (H^i_c)_\theta=t^2h^5+t^3h^6$,   thus
for    their    union   (of    which    they    are   connected    components)
$\sum_i   h^i\cdot   (H^i_c)_\theta=   3t^2h^5+3t^3h^6$.   The   union
$$\bX(s_2s_3s_1s_4s_3)\coprod\bX(s_2s_1s_3s_4s_3)\coprod\bX(s_2s_1s_3s_1s_4)$$
is    open    in    $\bX$;    the   corresponding    long    exact    sequence
gives       $$0\to      H^4_c(\bX)_\theta\to3t^2\theta\to3t^2\theta\to
H^5_c(\bX)_\theta\to   0,$$   and  $H^6_c(\bX)_\theta=3t^3\theta$.   There
exists   thus   an   integer   $\eta\le  3$   such   that   $\sum_i   h^i\cdot
H^i_c(\bX)_\theta=\eta t^2(h^4+h^5)\theta+3t^3h^6\theta$.    The   long
exact sequence  corresponding to the  union of  $\bX$ and $\bX(w)$  shows then
that only  the characters $t^2$ and  $t^3$ of $\genby F$  may occur in
the cohomology of $\bX(w)$, and gives  for the corresponding isotypic parts of
the cohomology:
$H^5_c(\bX(w))_{\theta,t^2}=H^6_c(\bX(w))_{\theta,t^2}=
\eta t^2\theta$,
and  $$0\to H^6_c(\bX(w))_{\theta,t^3}
\to t^3\theta\to3t^3\theta\to H^7_c(\bX(w))_{\theta,t^3}\to 0.$$
By \cite[3.3.22]{DMR} we have $H^5_c(\bX(w))=0$, thus also
$H^6_c(\bX(w))_{\theta,t^2}=0$.

To lift  the ambiguity on  $(\theta,t^3)$, we now  use \cite[3.3.21]{DMR};
we  take  for  $w'$  a  Coxeter  element,  which  is  not  conjugate  to  $w$;
thus   any  eigenvalue   of  $F$   on  $H^6_c(\bX(w))_\theta$   must  have   a
module  less  than  $q^3$.  This  shows  that  $H^6_c(\bX(w))_\theta=0$,  thus
$H^7_c(\bX(w))_{\theta,t^3}=2t^3\theta$.
\end{proof}

From the values  of the generic degrees of $\Hcyc$  (see \cite[5.A]{BrMa}), we
see  that  if  we  denote  by  $\rho(x_1,x_2,x_3)$ where $x_i\in\{1,q^2\}$  the
1-dimensional  representation  of  $\Hcyc$  given  by  $T_i\mapsto  x_i$,  and
$\rho^+$ (resp.  $\rho^-$) the  irreducible 2-dimensional  representation where
$T_1T_2T_3$ acts  as the scalar  $q^3$ (resp.  $-q^3$), we have  the following
equalities, if we denote by $m_\rho$ the multiplicity of $\rho$ in $
\sum_i(-1)^iH^i_c(\bX(\bw))$:
$m_{\rho(1,1,1)}=\dim \St$, $m_{\rho(q^2,q^2,q^2)}=\dim\Id$,
$m_{\rho(q^2,q^2,1)}= m_{\rho(1,q^2,q^2)}= m_{\rho(q^2,1,q^2)}=
-\dim\gamma_{1^2+}=-\dim\gamma_{1^2-}=-\dim\gamma_{21^2}$,
$m_{\rho(q^2,1,1)}= m_{\rho(1,q^2,1)}= m_{\rho(1,1,q^2)}=
-\dim\gamma_{2+}=-\dim\gamma_{2-}=-\dim\gamma_{31}$,
$m_{\rho^+}=\dim\gamma_{1.21}$ et $m_{\rho^-}=-\dim\theta$.
Thus we can determine $\oplus_i H^i_c(\bX(\bw))$ as a $\GF\times\Hcyc$-module
up to an ambiguity on the correspondance between characters of $\GF$ and 
of $\Hcyc$ which appear in the same degree.

\vfill\eject

\end{document}